\documentclass[12pt]{article}
\usepackage{style}
\title{HOMOLOGICAL MIRROR SYMMETRY FOR FUNCTORS BETWEEN FUKAYA CATEGORIES OF VERY AFFINE HYPERSURFACES}
\date{\today}

\author{Benjamin Gammage and Maxim Jeffs}

\begin{document}

\maketitle

\begin{abstract}
We prove that homological mirror symmetry for very affine hypersurfaces respects certain natural symplectic operations (as functors between partially wrapped Fukaya categories), verifying conjectures of Auroux. These conjectures concern compatibility between mirror symmetry for a very affine hypersurface and its complement, itself also a very affine hypersurface. We find that the complement of a very affine hypersurface has in fact two natural mirrors, one of which is a derived scheme. These two mirrors are related via a non-geometric equivalence mediated by Kn\"orrer periodicity; Auroux's conjectures require some modification to take this into account. Our proof also introduces new techniques for presenting Liouville manifolds as  gluings of Liouville sectors.
\end{abstract}

\tableofcontents

\section{Introduction}

\subsection{Motivation and Background}

In some sense, homological mirror symmetry (HMS) is an under-specified problem: there may be many possible equivalences between two categories, some of which clearly involve making choices. For instance, many proofs of HMS proceed by matching endomorphism algebras of chosen collections of generators on the $A$-side and $B$-side. How do we nail down a particular distinguished HMS equivalence? One answer is that we should ask that HMS intertwines naturally-defined functors on the A-side and B-side. Such results have previously appeared in the literature: see for example \cite{AndrewJeff, Andrew, Nadler,Catherine, HackingKeating}, and work in progress by  Cannizzo-Azam-Lee-Liu; in addition, it has long been known (starting with \cite{elliptic_curve, Fukaya}) that mirror symmetry for abelian varieties respects monoidal structures. 

The example that will be most relevant to us is mirror symmetry for anticanonical divisors in toric varieties: it is known by work of various authors \cite{Abouzaid_toric1, Abouzaid_toric2, FLTZ, Andrew, Kuwagaki,Zhou-ccc} that for a toric Fano variety $X$ (satisfying some assumptions such as those listed in \S \ref{subsec:terms} below) there is a Landau-Ginzburg mirror $((\CC^{\ast})^n, f)$ so that there is a quasiequivalence
\begin{equation}\label{eq:torichms}
    \scr{W}( (\CC^{\ast})^n, f) \simeq \mathrm{Coh}(X),
\end{equation}
where here $\scr{W}$ denotes the (split-closure of the twisted complexes over the) partially-wrapped Fukaya category, $\mathrm{Coh}$ denotes the (dg-derived) category of coherent sheaves, and $f$ is a Laurent polynomial (see \S \ref{subsec:terms} for details). Moreover, it was shown \cite{BenVivek, Zhou} using microlocal sheaf methods that the very affine hypersurface $H = f\inv(0)$ is mirror to the toric anticanonical divisor $i_D: D \hookrightarrow X$. Combined with the comparison theorem of \cite[Theorem 7.22]{GPS3} between microlocal sheaves and the Fukaya category, this entailed an equivalence 
\begin{equation*}
\scr{W}(H) \simeq \mathrm{Coh}(D).
\end{equation*}
Moreover, there is a commutative diagram \cite[Example 7.25]{GPS3}:
\begin{equation}\label{eqn:diagram1}
\begin{tikzcd}[baseline=(current  bounding  box.center)]
\scr{W}(H) \ar{r}{\sim} \ar{d}{\cup} & \mathrm{Coh}(D) \ar{d}{i_{D\ast}}\\
 \scr{W}( (\CC^{\ast})^n, f) \ar{r}{\sim} & \mathrm{Coh}(X),
\end{tikzcd}
\end{equation}
where $\cup$ denotes the cup functor \cite{sylvanOrlov} (envisioned by Abouzaid-Ganatra \cite{GPS3}) and the bottom equivalence is the equivalence \eqref{eq:torichms} from \cite[Theorem 1.2]{Kuwagaki} combined with the work of \cite{GPS3}. This is the prototype of the theorems we prove: we show that homological mirror symmetry respects certain naturally defined functors on each side of mirror symmetry. 

Auroux observed in \cite{Speculations} that the complement $(\CC^{\ast})^n\setminus H$ is also a very affine hypersurface, mirror to the toric anticanonical divisor $Z$ inside the canonical bundle $K_X$. Our functoriality results show that the HMS equivalences of \cite{BenVivek} intertwine symplectically-defined functors to and from $\scr{W}((\CC^{\ast})^n\setminus H)$, with various pullback and pushforward maps coming from the algebraic geometry of the mirror $Z$. We shall say informally, that such functors \textit{correspond under mirror symmetry} when a commutative diagram like (\ref{eqn:diagram1}) exists (with respect to specified mirror symmetry equivalences).

\subsection{Statement of Results}

On the $A$-side of mirror symmetry, Auroux in \cite{Speculations} sketches constructions of natural functors relating $\scr{W}((\CC^{\ast})^n\setminus H)$ to $\scr{W}(H)$ and $\scr{W}( (\CC^{\ast})^n, f)$:
\begin{itemize}
    \item The ($\ZZ/2$-graded) \textit{restriction functor} $\rho: \scr{W}((\CC^{\ast})^n\setminus H) \to \scr{W}(H)$ which takes a Lagrangian in $\comp$ to its `ends' along the removed fiber $H$ (Definition \ref{defn:rho});
    \item The \textit{lifting functor} $j: \scr{W}(H) \to \scr{W}((\CC^{\ast})^n\setminus H)$ which parallel transports a Lagrangian in $H$ along a ray from $0$ to $\infty$ avoiding all of the critical values of $f$ (Definition \ref{defn:lifting});
    \item The \textit{wrapping at infinity functor} $\alpha_\infty: \scr{W}( (\CC^{\ast})^n, f) \to \scr{W}((\CC^{\ast})^n\setminus H)$ which turns on the wrapping around infinity (Definition \ref{defn:alphas}, (1));
    \item The \textit{wrapping around $H$ functor} $\alpha_0: \scr{W}( (\CC^{\ast})^n, f) \to \scr{W}((\CC^{\ast})^n\setminus H)$ which turns on the wrapping around the removed fiber $H$ (Definition \ref{defn:alphas}, (2)).
\end{itemize}
We construct versions of these functors in detail in the sections below. Auroux conjectures that these functors should satisfy the relations
\begin{enumerate}
    \item $\rho \circ \alpha_0 \simeq \cap$;
    \item $\rho \circ \alpha_\infty  \simeq 0$;
    \item $\rho \circ j \simeq \mathrm{id}$;
\end{enumerate}
as well as the exact triangle
\begin{center}
\begin{tikzcd}
j \cap \ar{rr}{+1} & {} & \alpha_{\infty}\ . \ar{dl}\\
{} & \alpha_0 \ar{ul} & {}
\end{tikzcd}
\end{center}

Auroux then makes four conjectures about the above functors and their mirrors in algebraic geometry, which we recall here, indicating afterward the section of this paper in which each is discussed.

\begin{conjecture}[\cite{Speculations}]\label{conj:main}
The following functors correspond under mirror symmetry:
\begin{enumerate}
    \item (Conjecture 1.2) The functor $\rho$ is mirror to the quotient functor $q: \mathrm{Coh}(Z) \to \mathrm{Sing}(Z)$  composed with the Kn\"orrer periodicity equivalence $\mathrm{Sing}(Z) \to \mathrm{Coh}(D)$ (\S \ref{sec:rho});
    \item (Conjecture 1.3 (1)) The functor $\alpha_0$ is mirror to $i_{X\ast}$ for the inclusion $i_{X}: X \to Z$ (\S \ref{sec:gluing});
    \item (Conjecture 1.3 (2)) The functor $\alpha_\infty$ is mirror to $\pi_{X}^{\ast}$ for the projection $\pi_{X}: Z \to X$ (\S \ref{sec:rho});
    \item (Conjecture, \S 6.2) The functor $j$ is mirror to $j_{D \ast}p_{D}^{\ast}(\scr{K}_{X}\inv|_{D} \otimes \cdot \;)$ for $j_{D}: K_{X}|_D \to Z, p_{D}: K_{X}|_D \to D$ the inclusion and projection respectively (\S \ref{sec:j}).
\end{enumerate}
\end{conjecture}

In this paper, we prove these conjectures, after suitably modifying them to take into account grading data. The main results of this paper can be summarized as follows:

\begin{theorem*}
Under the mirror symmetry equivalences of \cite{Kuwagaki} and \cite{BenVivek} as in Theorem \ref{thm:alphas}, composed with the Kn\"orrer periodicity equivalences of Theorem \ref{prop:orlov} and Lemma \ref{lemma:small_orlov}: 
\begin{enumerate}
    \item (Theorem \ref{thm:alphas}, (2)) The functor $\alpha_0:\scr{W}( (\CC^{\ast})^n, f) \to \scr{W}_{\infty}(\comp)$ is mirror to $i_{X\ast}$ for the inclusion $i_{X}: X \to Z$;
    \item (Theorem \ref{thm:alpha_infinity}) The functor $\alpha_\infty:\scr{W}( (\CC^{\ast})^n, f) \to \scr{W}_{0}(\comp) $ is mirror to $\pi_{X}^{\ast}[1]$ for the projection $\pi_{X}: Z \to X$;
    \item (Theorem \ref{thm:rho}) The functor $\rho: \bar{\scr{W}}_{0}(\comp) \to \bar{\scr{W}}(H)$ is mirror to the quotient functor $q: \bar{\mathrm{Coh}}(Z) \to \bar{\mathrm{Sing}}(Z)$; 
    \item (Theorem \ref{thm:bigraded_rho}) The graded lift of $\rho: \bar{\scr{W}}_{0,\infty}(\comp) \to \bar{\scr{W}}(H)$ is (graded) mirror to the graded lift of $q: \bar{\mathrm{Coh}}_{\GG_m}(Z) \to \bar{\mathrm{Sing}}_{\GG_m}(Z)$;
    \item (Theorem \ref{thm:negative_lift}) The functor $j$ is mirror to $j_{D \ast}p_{D}^{\ast}$ for $j_{D}: K_{X}|_D \to Z, p_{D}: K_{X}|_D \to D$ the inclusion and projection respectively.
\end{enumerate}
and furthermore there is an exact triangle (Corollary \ref{cor:triangle}):
\begin{center}
\begin{tikzcd}
j \cap \ar{rr}{+1} & {} & \alpha_{0} \ar{dl}\\
{} & \alpha_\infty \ar{ul} & {}
\end{tikzcd}
\end{center}
\end{theorem*}

The definitions of the categories $\scr{W}_{\infty},\scr{W}_{0}$, and $\scr{W}_{0,\infty}$ will be provided in \S \ref{subsec:terms}, while the bar over $\overline{\scr{W}}, \overline{\mathrm{Coh}}$ denotes the reduction of the grading $\mathrm{mod}-2$.

\begin{remark}
Our functors differ from those defined by Auroux in \cite{Speculations} by an application of the monodromy: Auroux defines $\alpha_0, \alpha_{\infty}$ so that the Lagrangians both begin parallel transport in the direction of negative real infinity in the base. Because we define $\alpha_0,\alpha_{\infty}$ as sectorial inclusions instead, the Lagrangians go directly towards negative and positive real infinity respectively  (cf. \cite[Fig. 4]{Speculations}); this means that our results sometimes differ from the conjectures in \cite{Speculations} by a tensor product with the canonical bundle $K_X$ (cf. Theorem \ref{thm:negative_lift}) and a grading shift.
\end{remark}

\begin{remark}
    The hypotheses for the conjectures of \cite{Speculations} are somewhat more general than our assumptions in \S \ref{subsec:terms}; in particular, Auroux does not assume that $X$ is Fano. While the work of \cite{BenVivek} does extend to the non-Fano setting, there are interesting and subtle questions concerning different presentations of the skeleta that we do not address in this article. Hence there remain interesting open questions about the conjectures of \cite{Speculations}. 
\end{remark}

\subsection{Notation and Terminology}\label{subsec:terms}

\subsubsection{$A$-side Definitions}

Suppose $H \subset (\CC^{\ast})^n$ is a smooth algebraic hypersurface (a \textit{very affine hypersurface}), and choose a Laurent polynomial $f$ in $n$ variables so that $H = f\inv(0)$. We will write $f$ in the form
\begin{equation*}
    f(z) = \sum_{\alpha \in A} c_{\alpha} T^{\varphi(\alpha)} z^{\alpha}
\end{equation*}
where $A \subset \ZZ^n$ is the set of all multi-indices $\alpha$ of monomials in $f$; the $c_{\alpha}$ are for the moment some arbitrary non-zero complex numbers; $T$ is a small real number; and $\varphi$ is a function $A \to \RR$. Let $P$ be the polytope in $\RR^n$ given by the convex hull of $A$; then we require the function $\varphi$ to be the restriction of a convex piecewise-linear function on $P$, whose maximal domains of linearity define a polyhedral decomposition of $P$ whose vertices are exactly the points in $A$. 

We will further require $f$ to satisfy the following conditions:
\begin{enumerate}
    \item All of the non-zero points of $A$ lie on the boundary of $P$;
    \item The polyhedral decomposition of $P$ induced by $\varphi$ is maximal (all the cells of the polyhedral decomposition are congruent to standard simplices under the action of $\mathrm{GL}_n(\ZZ)$);
    \item The Laurent polynomial $f$ has a constant term $0 \in A$, and it is a vertex of every maximal cell in the polyhedral decomposition.
    % \item The polytope $P$ is perfectly centered (see \cite{BenVivek, Zhou}).
\end{enumerate}

These hypotheses will ensure that the mirror toric variety $X$ is smooth and Fano. Moreover, the complement of the tropical amoeba $\Pi_f \subset \RR^n$ will have exactly one compact component by hypothesis (3), so that the mirror for $z_0 f +1$ can be canonically identified with the canonical bundle $K_X$, and all of the critical points of $f$ can be made, by hypothesis (2), to live over an arbitrarily small neighbourhood of $1$ by taking $T$ sufficiently small \cite[Lemma 5.2]{Speculations}. Hypotheses (1) and (2) may be weakened if one is willing to work in the context of toric stacks as in \cite{BenVivek}: we will remark on this below. We will take $\varphi$ to be a generic small perturbation of the function $\varphi(\alpha)=1$ for every $\alpha \in A\setminus\set{0}$.

The complement of $H$ can be considered as a very affine hypersurface in $(\CC^{\ast})^{n+1}$ in one of two natural ways. It may be defined by either the Laurent polynomial $\tilde{f}_{0} = f + z_{0}\inv$ or by $\tilde{f}_{\infty} = z_{0} f + 1$, where $z_0 \in \CC^{\ast}$ is the first factor in $(\CC^{\ast})^{n+1}$. Of course, these are identical as algebraic varieties, and we shall call them $(\CC^{\ast})^n \setminus H$. Consider $\phi:(\CC^{\ast})^n \to \RR$, $\phi(z) = |\mathrm{Log}(z)|^2$ the standard Stein function on $(\CC^{\ast})^n$; by restriction both $H$ and $(\CC^{\ast})^n \setminus H \subset (\CC^{\ast})^{n+1}$ have the structure of Stein manifolds, and thus Liouville manifolds, allowing us to talk about their wrapped Fukaya categories. Note that the Stein structure obtained by restricting $\phi$ to $(\CC^{\ast})^n \setminus H$ is simple Liouville homotopic (see Definition \ref{def:simple}) to the Stein structures constructed in \cite[p.13]{periodicity} on the complement of a hypersurface in a Stein manifold. 

For us, the difference between the spaces $\{\tilde{f}_0=0\}$ and $\{\tilde{f}_\infty=0\}$ will be in the different {\em grading/Maslov} data --- namely, a choice of trivialization of the canonical line bundle, which determines a $\ZZ$-grading on the Fukaya category --- with which we equip them.
Given a hypersurface $H\subset (\CC^\ast)^n$, a function $f$ for which $H=\{f=0\}$ defines a trivialization $\eta$ of the canonical bundle of $H$: writing
%In addition we obtain a trivialization $\eta$ of the determinant line bundle of the tangent bundle of $H$ via its embedding in $(\CC^{\ast})^n$ by the wedge product with $\dd f$. Consider the standard holomorphic top form on $\cstar{n}$:
\begin{equation*}
    \Omega = \frac{\dd z_1}{z_1} \wedge \cdots \wedge \frac{\dd z_n}{z_n}
\end{equation*}
for the standard holomorphic top form on $H,$ we define 
define $\eta$ by the condition $\eta \wedge \dd f = \Omega$. 
%this gives a $\ZZ$-grading on $\scr{W}(H)$. 
Similarly, $\dd \tilde{f}_0$ and $\dd \tilde{f}_\infty$ give trivializations $\eta_0, \eta_{\infty}$ of the canonical bundle of $\comp$. In concrete terms, $\eta_0$ will be obtained as the restriction $\Omega|_{\comp}$ while $\eta_{\infty}$ will be $f\inv \Omega|_{\comp}$. Note that $\eta_0$ therefore extends across the removed fiber $H$ while $\eta_{\infty}$ does not. Changing the gradings from $\eta_0$ to $\eta_{\infty}$ amounts to dividing $\eta_0$ by $f$. In terms of cohomology, this corresponds to twisting the trivialization $\eta_0$ by the class $[\mathrm{arg}(f)] \in H^1(\comp; \ZZ)$. 

These two trivializations define two \textit{different} $A_{\infty}$ categories, which we shall call  $\scr{W}_0(\comp)$ and $\scr{W}_\infty(\comp)$. Note that these have \textit{different objects}: not every Lagrangian in $\comp$ has a graded lift for both gradings. Alternatively, we can consider a \textit{graded $A_{\infty}$ category} $\scr{W}_{0,\infty}(\comp)$, whose objects are restricted further to those that admit graded lifts for \textit{both} of the gradings. There are of course other grading structures one could consider on $\comp$ but these are less natural from the $B$-side perspective.

\begin{remark}
It will be important to note that choosing a Laurent polynomial $f$ provides $H$ with a framing inside $\cstar{n}$: horizontal lifts of $\partial_x,\partial_y \in T_0 \CC$ provide a trivialization of the normal bundle of $H$ in $\cstar{n}$.
\end{remark}

We shall use the following terminology: a \textit{sectorial decomposition} or \textit{sectorial gluing} of a Liouville manifold $X$ is a sectorial hypersurface $F$ in $X$ so that $X$ is the union of two Liouville sectors $X_1, X_2 \subset X$ intersecting along their boundaries $\partial X_1 = \partial X_2 = F$. This is what is informally referred to as `sector gluing' in \cite[\S 11.2]{GPS2}, analogous to the gluing of Liouville pairs in \cite[(12.17)]{GPS2}. See \S \ref{subsec:gluing} for our setup concerning the Liouville sector $(X,W)$ associated to a holomorphic function $W:X \to \CC$ on a Stein manifold $X$ (Definition \ref{def:general_sector}). 

\subsubsection{$B$-side Definitions}

In the above, the boundary of the polytope $P$ has a polyhedral decomposition; let $\Sigma_X$ be the toric fan of cones on the faces of this polyhedral decomposition. The rays of this fan are generated by the vectors in $A \subset \ZZ^n$. Let $X$ be the associated complex toric variety; let $i_D: D \to X$ be the inclusion of the toric anticanonical divisor, which is the zero-locus of a section $s$ of the anticanonical bundle $\scr{K}_{X}\inv$. We use $\scr{K}_X$ to denote the canonical bundle of $X$ as a sheaf on $X$ and we will denote the total space of the canonical bundle by $K_X$; this is also a toric variety with fan $\tilde{\Sigma}_X \subset \RR^{n+1}$ whose rays are generated by $(1, -\alpha)$ for $\alpha \in A$, and where $(1, -\alpha_i)$ span a cone if and only if the $\alpha_i$ span a cell of the polyhedral decomposition of the polytope $P$ induced by $\varphi$. Moreover let $\tilde{S}$ denote the toric anticanonical divisor inside $X \times \mathbb{A}^1$ and write $\hat{Z}$ for the toric anticanonical divisor inside $K_X \times \mathbb{A}^1$ (the `grand anticanonical divisor' inside the canonical bundle of the canonical bundle). Also, let $X \times \mathbb{A}^1[-1]$ denote the derived scheme where $\mathbb{A}^1[-1] = \mathrm{Spec}\;\CC[t]$ where $t$ has degree $-1$. We will use $\tilde{Z}$ to denote the derived scheme obtained by the gluing of $X \times \set{0}$ along $D$ to $D \times \AA^1[-1]$. Under hypotheses (1) and (2) above, $X$ must be smooth so all non-zero sections of $\scr{K}_{X}\inv$ will be regular.

\begin{remark}\label{rem:unimodular}
Hypothesis (2) can be weakened to require only that cells of the polyhedral decomposition are simplices, not necessarily of minimal volume. Given a fan $\Sigma_X$ and a choice of generator $\alpha$ for each ray, we can associate a \textit{toric (Deligne-Mumford) stack} $X$ \cite{BCS}.
The simplicial hypothesis guarantees that the stack $X$ will be smooth. The homological mirror symmetry theorem of \cite{BenVivek} is stated under this condition, and we expect that all of the results of this paper continue to hold in the stacky case.
%Under this weaker condition, $X$ will be smooth as a Deligne-Mumford stack \cite{BenVivek}. All of the algebraic geometry results we use in this paper should apply equally well in this case.
\end{remark}

In addition to the usual derived categories of coherent sheaves $\mathrm{Coh}(X),\mathrm{Coh}(D)$, we have various categories of sheaves associated to a smooth complex quasi-projective variety $X$ with a regular function $W: X \to \CC$ with a single critical fiber $X_0$ over $0$:
\begin{itemize}
    \item the matrix factorization category $\mathrm{Coh}(X,W)$ (Definition \ref{def:mf_cat});
    \item the singularity category $\mathrm{Sing}(X_0)$ (Definition \ref{def:sing}).
\end{itemize}
If moreover we have an affine algebraic group $G$ acting on $X$, and $\chi$ is a character of $G$ so that $W$ is $\chi$ -semi-invariant, then we also have:
\begin{itemize}
    \item the category of graded matrix factorizations $\mathrm{Coh}_G(X,W)$ (Definition \ref{def:graded_mf});
    \item the equivariant singularity category $\mathrm{Sing}_{G}(X_0)$ (Definition \ref{def:graded_sing}).
\end{itemize}
where the character $\chi$ is left implicit.

Throughout, all dg/$A_{\infty}$ categories will be taken to be $\CC$-linear (equipped as usual with their internal/cohomological $\ZZ$-grading); by a \textit{graded dg category} (or $A_{\infty}$ category), we shall mean a dg category equipped with an \textit{additional} $\ZZ$-grading: in other words, a {\em grading} on a dg category is a $\Coh(B\GG_m)$-linear structure. A \textit{$2$-periodicity structure} on a dg category is a choice of $\cz$-linear structure, where $|\beta|=2$; a \textit{$2$-periodic dg category} will be one that admits a $2$-periodicity structure. Equivalently, localizing this $\cz$-linear dg category at the natural transformation $\beta$ gives a differential $\ZZ/2$-graded category. For an $A_{\infty}$ category, reducing the homological grading mod-$2$ yields a $\cz$-linear $A_{\infty}$ category; we shall use a bar to denote the reduction of the homological grading mod-$2$: for instance, $\overline{\scr{W}}_{0}(\comp)$ and $\overline{\mathrm{Coh}}(Z)$. Equivalences of $A_{\infty}$ categories in particular descend to the level of $\cz$-linear $A_{\infty}$ categories by reducing the grading mod-$2$.

In the following, for a variety $X$, $\mathrm{Coh}(X,0)$ will denote the matrix factorizations of zero, that is, a $\cz$-linear version of the category of coherent sheaves $\mathrm{Coh}(X)$. 
As we will explain in \S \ref{sec:mfs}, we can recover the $\CC$-linear dg category $\mathrm{Coh}(X)$ as a category of {\em graded} matrix factorizations of zero; we will denote this latter category as $\mathrm{Coh}_{\GG_m}(X,0).$ The geometric meaning of this $\GG_m$ action will be explained in \S \ref{sec:mfs}.
%We can restore this grading using an action of $\GG_m$: the category $\mathrm{Coh}_{\GG_m}(X,0)$ will be equivalent the original ($\ZZ$-graded) category of coherent sheaves $\mathrm{Coh}(X)$. See \S \ref{sec:mfs} for more discussion.

\begin{remark}
Henceforth all of the functors below will be $\ZZ$-graded functors between $\ZZ$-graded categories (i.e. dg-functors between dg categories), unless otherwise stated. We will use the definition of the (partially) wrapped Fukaya category from \cite{GPS1}, denoted there as $\scr{W}$. For us, $\scr{W}$ will instead denote the split-closure of the $A_{\infty}$ category of twisted complexes over this (partially) wrapped Fukaya category; $\mathrm{Coh}$ will denote the dg-derived category of coherent sheaves. All functors are derived, and all colimits are to be understood in the appropriate homotopical sense.
\end{remark}

\subsubsection{An Example}

In this short subsection, we consider the example of $H = \set{-1}$. We hope that even though this example is simple it will help the reader keep track of the various spaces and functors involved. On the $B$-side we can identify
\begin{itemize}
    \item $X = \mathrm{Spec}(\CC[x])$;
    \item $K_X = \mathrm{Spec}(\CC[x,y])$;
    \item $\AA^1[-1] = \mathrm{Spec}(\CC[t])$, $|t|=-1$.
\end{itemize}
On the $A$-side, there are two different holomorphic functions, $f=z+1$ and $1/f$, on the pair of pants $\Pi_1 = \CC^{\ast}\setminus H$. In this case, the two corresponding choices of gradings can be easily described:
\begin{align*}
    \eta_0 & = \frac{\dd{z}}{z}, & \eta_{\infty} & = \frac{\dd{z}}{z(z+1)}
\end{align*}
For $\eta_0$, the simple clockwise Reeb orbit around $-1$ will have degree $2$, while those around the other two punctures will have degree $0$. The situation is reversed for $\eta_{\infty}$, where the simple clockwise Reeb orbit around $\infty$ will have degree $2$. We hope that this goes some way towards explaining our choice of notation.

This example is discussed further in \S \ref{subsec:example}.

\subsubsection{Summary of Notation}

We list here for the reader's convenience all the morphisms of schemes we use: first the inclusions:
\begin{itemize}
    \item $i_D:D \to X$; $i_Z:Z \to K_X$; $i_{\tilde{Z}}: \tilde{Z} \to X \times \mathbb{A}^1[-1]$, $i_{\hat{Z}}: \hat{Z} \to K_{X} \times \mathbb{A}^1$, $i_{X}:X \to Z$, $i_{\tilde{X}}: X \to \tilde{Z}$; 
    \item $j_{D}: K_{X}|_D \to Z$, $j_{\tilde{D}}: D \times \mathbb{A}^1[-1] \to \tilde{Z}$, $j_{Z}: Z \times \mathbb{A}^1[-1] \to \hat{Z}$; $j_{\tilde{Z}}: K_{X \times \mathbb{A}^1}|_{\tilde{Z}} \to \hat{Z}$;
    \item $\iota_{\tilde{D}}: D \to D \times \mathbb{A}^1[-1]$; $\iota_{D}: D \to K_{X}|_D$; $\iota_{Z}: Z \to Z \times \mathbb{A}^1[-1]$; $\iota_{\tilde{Z}}: \tilde{Z} \to K_{X \times \mathbb{A}^1}|_{\tilde{Z}}$;
    \item $\iota_{X}: X \to K_{X}$; $\iota_{\tilde{X}}: X \to X \times \mathbb{A}^1[-1]$; $\iota_{K_X}: K_{X} \to K_{X} \times \mathbb{A}^1[-1]$; $\iota_{X \times \mathbb{A}^1[-1]}: X \times \mathbb{A}^1[-1] \to K_{X} \times \mathbb{A}^1[-1]$.
\end{itemize}
and the projections:
\begin{itemize}
    \item $p_{D}: K_{X}|_D \to D$; $p_{\tilde{D}}: D \times \mathbb{A}^1[-1] \to D$; $p_{Z}: Z \times \mathbb{A}^1[-1] \to Z$; $p_{\tilde{Z}}: K_{X \times \mathbb{A}^1}|_{\tilde{Z}} \to \tilde{Z}$;
        \item $\pi_{Z}: K_X \times \mathbb{A}^1[-1] \to K_{X}$; $\pi_{\tilde{Z}}: K_{X} \times \mathbb{A}^1[-1] \to X \times \mathbb{A}^1[-1]$; $\pi_{X}: Z \to X$, $\pi_{\tilde{X}}: \tilde{Z} \to X$.
\end{itemize}
and analogously for the classical scheme $\tilde{S}$.

\subsection{Summary of Methods}

In this section, we will give a somewhat heuristic overview of some of the main techniques used in this paper. The first is the homotopy pushout theorem for wrapped Fukaya categories of \cite[Theorem 1.28]{GPS2}, and the corresponding proper descent theorem for categories of coherent sheaves of \cite[Proposition 4.7.2.2]{GRI}. The two presentations of $\comp$, as $\tilde{f}_0=0$ or $\tilde{f}_{\infty}=0$, naturally give rise to two (singular) symplectic fibrations over $\CC$, given by $f$ and $1/f$, and two corresponding presentations of $\comp$ as the gluing of Liouville sectors; these in turn correspond to the decompositions of the two mirror spaces $\tilde{Z}$ and $Z$ into their irreducible components. This kind of gluing argument was envisioned in \cite[Remark 1.5]{Speculations}. However, this symplectic gluing is difficult to realize.

Given a polynomial function $W:X \to \CC$ on an affine variety $X$ with Stein function $\phi:X \to \RR$, the real hypersurface $F = \set{\mathrm{Re}(W) = C}$ should be the prototypical example of a sectorial hypersurface, presenting $X$ as the sectorial gluing of the two sectors $\set{\mathrm{Re}(W) \leq C}$ and $\set{\mathrm{Re}(W) \geq C}$ along $F$. While $F$ satisfies the second condition of \cite[Definition 1.2]{GPS1} with the function $I = \mathrm{Im}(W)$, implicit in this definition is that $F$ be tangent to the Liouville vector field of $X$ at infinity. This condition is in general impossible to ensure, even after a deformation of the Liouville structure on $X$. Criteria under which such a simple Liouville deformation exists will appear in future work of the second author and collaborators. Here we describe an alternative construction that will yield the same results (\S \ref{sec:gluing}).

Instead of rectifying the Liouville structure along the entire hypersurface $F$, we may apply a result of Sylvan \cite[Prop 2.6]{sylvan_lemma} to produce a neighbourhood of a fiber of $W$ near infinity with a standard Liouville form. If we let $H = W\inv(R_1) \cap \set{\phi \leq R_2 - \epsilon}$ then this proposition gives us for some $\rho>0$ (and sufficiently large $R_1, R_2$) a standard Liouville neighborhood $U$ of $H$ of the form $(\hat{H} \times \set{z \in \CC \;\;: \mathrm{Re}(z) > - \rho}, \lambda_{X}|_H + \lambda_{\CC}^{\mathrm{std}})$, properly embedded inside the completion $\hat{X}$ (see Figure \ref{fig:step1}). The reader is warned that the projection $z: U \to \CC$ need not be compatible with the original map $W:X \to \CC$ (though the fibers of $z$ will be \textit{isotopic} to fibers of $W$). We now have a sectorial hypersurface $\set{\mathrm{Re}(z) = 0} \subset U$ properly embedded inside $X$. In our setting, with $W = f$, this simply presents $\cstar{n}$ as the trivial gluing of $(\cstar{n},f)$ to $(H \times \CC,z)$. Removing the hypersurface $H$ instead inside the neighbourhood $U$ allows us to present $\comp$ as the sectorial gluing of $(\cstar{n},f)$ to $(H \times \CC^{\ast},z)$ (Proposition \ref{prop:sector_gluing}).

The same argument cannot work for the function $1/W:X \to \CC$ as $1/W$ is not defined everywhere on $X$ (and is no longer a polynomial). We may instead build the sectorial decomposition of $X$ induced by $1/W$ by `swapping zero and infinity' in the sectorial decomposition of $X$ induced by $W$. In more precise terms, this means starting with our previous construction and applying an isotopy of the fiber $H$ given by the global monodromy of $f$. By applying uniqueness results of \cite{GPS1}, one may see that it is always possible to find a corresponding family of sectorial hypersurfaces `following' the isotopy. The resulting sectorial decomposition presents $\comp$ as the sectorial gluing of $(\cstar{n}, f)$ to a $H$-bundle over $\CC^{\ast}$ twisted by the global monodromy of $f$ (see Figure \ref{fig:step5}, Proposition \ref{prop:A2}).

The second main technique we use in this paper is derived Kn\"orrer periodicity for gauged LG models (as in \cite{Hirano}). While the sectorial gluing of $\comp$ induced by $1/f$  corresponds under mirror symmetry to the $B$-side gluing of $Z$ from its irreducible components (Theorem \ref{thm:alphas}), on the other hand, the unusual grading $\eta_0$ means the sectorial decomposition of $\comp$ induced by $f$ induces an equivalence between the wrapped Fukaya category of $\comp$ with grading $\eta_0$ and coherent sheaves on the \textit{derived scheme} $\tilde{Z}$. To relate this back to the original mirror symmetry conjectures of \cite{Speculations} for $Z$, we need an equivalence between $\mathrm{Coh}(Z)$ and $\mathrm{Coh}(\tilde{Z})$. These two schemes are related by a combination of Koszul duality and derived Kn\"orrer periodicity: the two schemes have a common Kn\"orrer stabilization, to a category of (graded) matrix factorizations of $s y z$ on $K_X \times \AA^1$, where $y,z$ denote the coordinates on the fiber of $K_X$ and $\AA^1$ respectively (Proposition \ref{prop:orlov}). This equivalence between $\mathrm{Coh}(Z)$ and $\mathrm{Coh}(\tilde{Z})$ has surprising properties that make it the algebraic analog of passing from $f$ to $1/f$. This equivalence takes the functor $i_{\tilde{X}\ast}$ to $\pi_{X}^{\ast}[1]$ (Lemma \ref{lem:locally_free}) and vice-versa, and thus the induced equivalence on the $A$-side swaps the functor $\alpha_0$ with $\alpha_{\infty}$ (Theorem \ref{thm:alpha_infinity}), in some sense `swapping zero and infinity'. The general significance of such intriguing non-geometric equivalences is a subject that will be taken up in future work. 

\subsection*{Acknowledgements}

MJ would first like to thank Denis Auroux for his patience and guidance throughout the project, as well as for reading many drafts of this paper. Thank you also to David Nadler for initial conversations on this subject, and special thanks to Daniel Chupin and Kevin Lin for helpful conversations about descent. MJ also benefited from discussions with Sheel Ganatra and Andrew Hanlon. We would also like to thank the anonymous referee for their detailed suggestions.

MJ was partially supported by the Rutherford Foundation of the Royal Society of New Zealand, NSF grant DMS-1937869 and DMS-2202984, and by Simons Foundation grant \#385573.

BG is supported by an NSF Postdoctoral Research Fellowship, DMS-2001897.

\clearpage

\section{Gradings}

\subsection{Matrix factorization categories}\label{sec:mfs}

Suppose $X$ is a smooth complex quasi-projective variety and $W: X \to \CC$ a regular function with a single critical fiber $X_0$ over $0$. We have two categories naturally associated to the LG model $(X,W)$:

\begin{definition}\label{def:sing}
The singularity category $\mathrm{Sing}(X_0)$ is the dg-quotient category of $\mathrm{Coh}(X_0)$ by the full subcategory $\mathrm{Perf}(X_0)$ of perfect complexes on $X_0$. 
\end{definition}

The other is the category $\mathrm{Coh}(X, W)$ of matrix factorizations of $W$:

\begin{definition}\label{def:mf_cat}
A \textit{matrix factorization} of $W$ on $X$ consists of a pair $F_0,F_1$ of coherent sheaves on $X$, along with morphisms $\phi_0:F_0 \to F_1, \phi_1: F_1 \to F_0$ such that $\phi_0 \phi_1 = W \cdot \mathrm{id}$ and $\phi_1 \phi_0 = W \cdot \mathrm{id}$. For a pair of objects $E = (E_0,E_1,\phi^E), F = (F_0,F_1,\phi^F)$ in this category, the morphism complex is given by
\begin{align*}
    \mathrm{Hom}^{0}(E,F) & = \mathrm{Hom}(E_0,F_0) \osum \mathrm{Hom}(E_1,F_1) \\
    \mathrm{Hom}^{1}(E,F) & = \mathrm{Hom}(E_0,F_1) \osum \mathrm{Hom}(E_1,F_0)
\end{align*}
with differential $\dd^i(f) = \phi^F \circ f - (-1)^i f \circ \phi^E $. The dg-derived category formed by matrix factorizations of $W$ on $X$ is denoted $\mathrm{Coh}(X, W)$ (cf. \cite[Definition 2.10]{Hirano}).
\end{definition}

Note that $\sing(X_0)$ is a $2$-periodic dg category (as every object in $\sing$ has a $2$-periodic infinite resolution), while $\mathrm{Coh}(X,W)$ is a $\cz$-linear dg category. An isomorphism between them involves making a choice of $2$-periodicity structure on $\sing(X_0)$ induced by $W$ (since different choices of $W$ can have the same zero locus $X_0$): localizing at this $2$-periodicity natural transformation gives a differential $\ZZ/2$-graded category. See \cite[\S 2.5]{pascaleff} for further discussion.

%Scaling $W$ by a constant factor may yield a different $\ZZ/2$-graded category.

Let $i:X_0 \to X$ be the inclusion of the closed subscheme $X_0 = W\inv(0)$. Explicitly, this isomorphism works by taking $F$ a coherent sheaf on $X_0$ and pushing it forward to $X$ to get a two-step resolution by locally free sheaves:
\begin{center}
\begin{tikzcd}
 0 \ar{r} & E^{-1} \ar[bend left]{r}{\phi_{-1}} & E^0 \ar[bend left]{l}{\phi_0} \ar{r} & i_{\ast} F \ar{r} & 0
\end{tikzcd}
\end{center}
with $\phi_{-1}\phi_0 = W$ and $\phi_{0}\phi_{-1} = W$: this yields the corresponding matrix factorization. Pulling back to $X_0$, the $E_{i}$s remain locally free, and since $W$ vanishes on $X_0$ we have an exact sequence
\begin{equation*}
    0 \to F \to E^{-1} \to E^{0} \to  F \to 0
\end{equation*}
Iterating this, we have an infinite resolution of $F$ on $X_0$ by 
\begin{equation*}
  \cdots \to E^{-3} \to E^{-2} \to E^{-1} \to E^{0} \to  F \to 0
\end{equation*}
which is $2$-periodic as a complex of sheaves. 

%To have both of these categories be $\ZZ$-graded, we need to take a $\mathbb{G}_m$-action on $X$ so that $W$ is homogeneous.
We may restore the $\ZZ$-grading to these categories by equipping $X$ with a $\mathbb{G}_m$ action for which $W$ is homogeneous.
Suppose $G$ is an affine algebraic group acting on $X$, and $\chi$ is a character of $G$ so that $W$ is $\chi$-semi-invariant.

\begin{definition}\label{def:graded_mf}
The category $\mathrm{Coh}_G(X,W)$ of \textit{graded matrix factorizations} has objects given by pairs $F_0,F_1$ of $G$-equivariant coherent sheaves on $X$, along with $G$-invariant morphisms $\phi_1:F_1 \to F_0$, $\phi_0:F_0 \to F_1(\chi)$ satisfying $\phi_0 \phi_1 = W \cdot \mathrm{id}$ and $\phi_1(\chi) \phi_0 = W \cdot \mathrm{id}$. The morphism complex between two objects is now given by
\begin{align*}
    \mathrm{Hom}^{2n}(E,F) & = \mathrm{Hom}(E_0,F_0(\chi^n)) \osum \mathrm{Hom}(E_1,F_1(\chi^n)) \\
    \mathrm{Hom}^{2n+1}(E,F) & = \mathrm{Hom}(E_0,F_1(\chi^n)) \osum \mathrm{Hom}(E_1,F_0(\chi^{n+1}))
\end{align*}
with the same differential $\dd^i(f) = \phi^F \circ f - (-1)^i f \circ \phi^E $ as above.
\end{definition}
Compare \cite[Definition 2.2]{Hirano}, where this category is denoted $\mathrm{Coh}_G(X,\chi,W)$; we choose to suppress the choice of character $\chi$ in our notation. We refer to \cite{Hirano} for all our other conventions about matrix factorization categories and functors among them. 

For now, we return to the case where $G = \GG_m$, where this means that our matrix factorization category is $\ZZ$-graded. In the resulting triangulated category, the usual homological $\ZZ/2$-grading will be extended to a `sheared' $\ZZ$-grading where $[2] = (\chi)$.

\begin{definition}\label{def:graded_sing}
We define the equivariant singularity category $\mathrm{Sing}_{\GG_m}(X_0)$ to be the dg-quotient category of $\Coh_{\GG_m}(X_0)$, the $\GG_m$-equivariant coherent sheaves on $X_0$, by the full subcategory of equivariant perfect complexes.
\end{definition}

Now, the homological $\ZZ$-grading on $\mathrm{Sing}(X_0)$ is no longer $2$-periodic once we take into account the equivariant structure on the resolutions coming from $\Coh_{\GG_m}(X_0)$: in our previous explicit construction of the $2$-periodic resolution, one sees that in fact $E^{\bullet}[2] = E^{\bullet}(\chi)$. The equivalence between $\mathrm{Coh}(X,W)$ and $\mathrm{Sing}(X_0)$ proved by Orlov \cite{Orlov2} now extends to an equivalence of dg categories; the following theorem follows from \cite[Theorem 3.6]{Hirano} by taking $G = \mathbb{G}_m$.

\begin{theorem}($\ZZ$-graded Orlov's theorem) \label{thm:orlov}
Given a $\mathbb{G}_m$-action on $X$ for which $W$ is quasi-homogenous of weight $1$, we have an equivalence of dg categories
\begin{equation*}
    \mathrm{Sing}_{\GG_m}(X_0) \to \mathrm{Coh}_{\mathbb{G}_m}(X, W)
\end{equation*}
given by the pushforward under the inclusion $X_0 \hookrightarrow X$.
\end{theorem}

Observe that since $W$ is zero on $X_0$, pushing forward a complex of sheaves from $X_0$ yields a matrix factorization on $X$.

\begin{lemma}\label{lem:trivial_action} 
If $Y$ is any smooth algebraic variety with a trivial $\GG_m$ action, then $\mathrm{Coh}(Y)$ is equivalent to $\mathrm{Coh}_{\GG_m}(Y,0)$ as dg categories.
\end{lemma}

See \cite[Proposition 2.14]{Hirano}.

The basic $\ZZ$-graded version of Orlov's Kn\"orrer periodicity theorem \cite[Theorem 2.1]{Orlov3} that we will need is the following. 

\begin{lemma}($\ZZ$-graded Kn\"orrer periodicity)\label{lemma:small_orlov}
If $K_X$ carries a $\mathbb{G}_m$-action of weight 1 in the fibers, then there are equivalences of dg categories:
\begin{equation*}
    \mathrm{Coh}(D) \xrightarrow{j_{D\ast}p_{D}^{\ast}} \mathrm{Sing}_{\mathbb{G}_m}(Z) \xrightarrow{i_{Z\ast}} \mathrm{Coh}_{\mathbb{G}_m}(K_X, s(x)y)
\end{equation*}
where $y$ is the coordinate on the fiber of $K_X$.
\end{lemma}
\begin{proof}
This follows by combining Theorem \ref{thm:orlov} above with Hirano's $\ZZ$-graded enhancement of derived Kn\"orrer periodicity \cite[Theorem 1.2]{Hirano} with $G = \mathbb{G}_m$ and $\scr{E}$ the line bundle $\scr{K}_X$ on $X$ with section $s$.
\end{proof}

We will also need various upgraded versions of this. Applying \cite[Theorem 1.2]{Hirano} with $G = \mathbb{G}_m$ and $\scr{E}$ the trivial line bundle over $K_X$ gives us:
\begin{proposition}\label{prop:knorrer1}
There is an equivalence of dg categories
\begin{equation*}
    j_{Z\ast}p_{Z}^{\ast}: \mathrm{Coh}_{\mathbb{G}_m}(Z, 0) \to \mathrm{Coh}_{\mathbb{G}_m}(K_X \times \mathbb{A}^1, s(x) y z)
\end{equation*}
where here:
\begin{itemize}
    \item $\mathbb{G}_{m}$ acts trivially on $K_X$ and with weight $1$ on $\mathbb{A}^1$;
    \item $s$ is the section of $\scr{K}_{X}\inv$ defining the toric anticanonical divisor:
    \item $x$ is the coordinate on $X$; $y$ is the coordinate on the fiber of $K_X$; and $z$ is the coordinate on $\mathbb{A}^1$;
    \item $p_Z$ is the projection $Z \times \mathbb{A}^1 \to Z$; and $j_Z$ is the inclusion $Z \times \mathbb{A}^1 \to K_X \times \mathbb{A}^1$.
\end{itemize}
\end{proposition}
Similarly, applying Hirano's theorem to $\scr{E} = \scr{K}_{X}$ over $X \times \mathbb{A}^1$ with a $\mathbb{G}_m$-action of weight $1$ on the $\mathbb{A}^1$ factor gives:
\begin{proposition}\label{prop:knorrer2}
There is an equivalence of dg categories:
\begin{equation*}
    j_{\tilde{S}\ast}p_{\tilde{S}}^{\ast}: \mathrm{Coh}_{\mathbb{G}_m}(\tilde{S}, 0) \to \mathrm{Coh}_{\mathbb{G}_m}(K_X \times \mathbb{A}^1, s(x) y z)
\end{equation*}
where $\mathbb{G}_m$ acts trivially on $K_X$ and with weight $1$ on $\mathbb{A}^1$ as before; $p_{\tilde{S}}$ is the projection $\tilde{S} \times \mathbb{A}^1 \to \tilde{S}$ and $j_{\tilde{S}}$ is the inclusion $K_{X}|_{\tilde{S}} \to K_X \times \mathbb{A}^1$.
\end{proposition}
This follows from \cite[Theorem 1.2]{Hirano} since in order for $z s(x)$ to give a $\mathbb{G}_m$-\textit{invariant} section of $\scr{K}_{X}\inv$, the bundle $\scr{E} = \scr{K}_X$ must be given an equivariant structure with character $\chi_{-1}$. Then the fiber of the total space of $\scr{E} \otimes \scr{O}(\chi_1)$ has weight $0$ with respect to the $\mathbb{G}_m$-action.

Consider the derived scheme obtained by the gluing of $X \times \set{0}$ along $D$ to $D \times \AA^1[-1]$, which we denoted by $\tilde{Z}$.

\begin{proposition}\label{prop:compatibility}
We have an equivalence of dg categories:
 $\mathrm{Coh}(\tilde{Z}) \cong \mathrm{Coh}_{\mathbb{G}_m}(\tilde{S}, 0).$
\end{proposition}

For the proof, we proceed in a series of steps.

\begin{lemma}\label{lem:koszul}
(Koszul duality) The category $\mathrm{Coh}(\AA^1[-1])$ is equivalent to the category $\mathrm{Perf}(\CC[z])$ of perfect dg-modules over the graded ring $\CC[z],$ where $z$ is a variable of cohomological degree $2$.
\end{lemma}
\begin{proof}
See for instance \cite[Corollary 5.1.10]{ArinkinGaitsgory} and pass to compact objects.
\end{proof}

\begin{lemma}\label{lem:mf-koszul}
The category $\mathrm{Coh}(\AA^1[-1])$ is equivalent to the category $\mathrm{Coh}_{\GG_m}(\AA^1,0)$ of graded matrix factorizations on $\AA^1$ with a weight-$1$ $\GG_m$-action.
\end{lemma}
\begin{proof} Firstly, by the Koszul duality equivalence above, the category $\mathrm{Coh}(\AA^1[-1])$ is equivalent to $\mathrm{Perf}(\CC[z])$. Objects of this latter category consist of perfect dg-$\CC[z]$-modules $(V^{\bullet}, d^{\bullet})$, i.e. a chain complex where $V^{\bullet}$ are finite-dimensional vector spaces and $z: V^{i} \to V^{i+2}$ are linear maps commuting with $d^{\bullet}$. Objects of the category $\mathrm{Coh}_{\GG_m}(\AA^1,0)$ consist of triples $(F_0,F_1, \varphi)$ where $F_0,F_1$ are finitely-generated $\ZZ$-graded $\CC[x]$-modules (where $|x|=1$) and $\varphi_1: F_1 \to F_0$ and $\varphi_0: F_0 \to F_1[-1]$ are morphisms of graded modules with $\varphi_0 \varphi_1 = 0$ and $\varphi_1[-1] \varphi_0 = 0$. 

Given a dg-$\CC[z]$-module $(V^{\bullet}, d^{\bullet})$ one can produce a matrix factorization $T(V) = (F_0, F_1, \varphi)$ via:
\begin{align*}
    F_0 & = \bigosum_{n \in \ZZ} V^{2n}[n] & F_1 & = \bigosum_{n \in \ZZ} V^{2n+1}[n]
\end{align*}
where $F_i$ have $\CC[x]$-module structures given by the odd and even parts of $z$ considered as a map $z:F_i \to F_i[-1]$; and the maps $\varphi_0, \varphi_1$ come from the differentials $\varphi_0 = \osum_n d^{2n}[n]$ and $\varphi_1 = \osum_n d^{2n+1}[n]$. If $(V^\bullet,d^\bullet)$ is a perfect $\CC[z]$-module then $F_0,F_1$ are coherent $\CC[x]$-modules.

In the reverse direction, given a matrix factorization $(F_0,F_1, \varphi)$, one can produce a perfect dg-$\CC[z]$-module $S(F_0,F_1,\varphi)$ by writing each of $F_0,F_1$ as a sum of graded pieces:
\begin{align*}
    F_0 & = \bigoplus_{m \in \ZZ} A^{m} & F_1 & = \bigoplus_{m \in \ZZ} B^m
\end{align*}
and setting
\begin{align*}
    V^{2n} & = A^n & V^{2n+1} & = B^n
\end{align*}
with $\CC[z]$-module structure where $z$ acts by $x:V^{2n} = A^n \to A^{n+1} = V^{2n+2}$ on the even part, as $x:V^{2n+1}=B^n \to B^{n+1} = V^{2n+3}$ on the odd part, and the differential $d:V^{2n+1} \to V^{2n}$ is given by $\varphi_1:B^n \to A^{n}$ and $d:V^{2n} \to V^{2n-1}$ is given by $\varphi_0:A^n \to B^{n-1}$. Observe that $S \circ T(V^\bullet,d^\bullet)$ is equal to $(V^\bullet,d^\bullet)$ as a dg-$\CC[z]$-module, and that $S(F_0,F_1,d)$ is a perfect $\CC[z]$-module whenever $F_0,F_1$ are coherent $\CC[x]$-modules.

A degree-$2k$ morphism $f^{\bullet}: (V^{\bullet}, d_{V}^{\bullet}) \to (W^{\bullet}, d_{W}^{\bullet})$ between two dg-$\CC[z]$-modules consists of linear maps $f^n: V^n \to W^{n+2k}$ commuting with the action of $z$ and with $d_V,d_W$. Likewise, a degree-$2k$ morphism $g: T(V) \to T(W)$ of $\GG_m$-equivariant matrix factorizations consists of morphisms of graded $\CC[x]$-modules $g_0: T(V)_0 \to T(W)_0[-k]$ and $g_1: T(V)_1 \to T(W)_1[-k]$, commuting with $\varphi_V, \varphi_W$. Thus one obtains a morphism $T(f):T(V) \to T(W)$ simply by restricting $f^\bullet$ to the odd and even degrees of the complex: $f^{2n}:V^{2n}[n] \to W^{2n+2k}[n]$ and $f^{2n+1}:V^{2n+1}[n] \to W^{2n+2k+1}[n]$ give us $T(f)_0 : T(V)_0 \to T(W)_0[-k]$ and $T(f)_1 : T(V)_1 \to T(W)_1[-k]$. The case of odd-degreee morphisms is analogous.

It is evident that $T$ respects composition of morphisms (since it respects compositions of linear maps), so $T$ defines a functor $\mathrm{Perf}(\CC[z]) \to \mathrm{Coh}_{\GG_m}(\AA^1,0)$ (of ordinary categories). We have seen above that $T$ is essentially surjective, and we claim that the functor $T$ is also faithful. Any degree-$2k$ morphism $g: T(V) \to T(W)$ can be restricted to the graded pieces to give linear maps $g_A:A_{V}^m \to A_{W}^{m+k}$ and $g_B:B_{V}^m \to B_{W}^{m+k}$ which together give a linear map $g: V^m \to W^{m+2k}$ that commutes with $z$ and $d_V, d_W$, and hence yields a degree-$2k$ morphism $S(g): (V^\bullet, d_{V}^\bullet) \to (W^\bullet, d_{W}^\bullet)$ of dg-$\CC[z]$-modules. Starting with a morphism $f:(V^\bullet, d_{V}^\bullet) \to (W^\bullet, d_{W}^\bullet)$ of dg-$\CC[z]$-modules, one can see that $S(T(f))$ is equal to $f$ itself. A similar argument applies in the case of odd-degree morphisms.

We claim that $T$ is in fact a dg-functor. The differential on degree-$k$ morphisms between two dg-$\CC[z]$ modules $f^{\bullet}: (V^{\bullet}, d_{V}^{\bullet}) \to (W^{\bullet}, d_{W}^{\bullet})$ simply acts by
\begin{equation*}
    d^k(f) = d_W \circ f - (-1)^k f \circ d_V
\end{equation*}
while $T(d_W \circ f) = \varphi_{T(W)} \circ T(f)$ and $T(f \circ d_V) = T(f) \circ \varphi_{T(V)}$, by our construction of the matrix factorization $T(V)$ which used $d_V$ to define $\varphi_V$. Thus $T(d^k f) = \varphi_{T(W)} \circ T(f) - (-1)^k  T(f) \circ \varphi_{T(V)}$, which is the differential acting on the degree-$k$ morphism $T(f)$ of equivariant matrix factorizations. Hence $T$ gives an equivalence of dg categories.
\end{proof}
 
From this lemma it follows that also:

\begin{corollary}\label{cor:divisor-factorization}
There is an equivalence of dg categories between $\mathrm{Coh}(D \times \AA^1[-1])$ and $\mathrm{Coh}_{\GG_m}(D \times \AA^1,0)$, where $\GG_m$ acts trivially on the $D$ factor and with weight $1$ on $\AA^1$.
\end{corollary}

%Compare \cite[\S 2.2]{pascaleff}. 
The following can also be considered a more general form of Lemma \ref{lem:mf-koszul}:

\begin{proposition}\label{prop:mf_to_coh}
For a scheme $Y$ with an action by an affine algebraic group $G$ and a choice of character $\chi$, there is an equivalence of dg categories between the matrix factorization category $\mathrm{Coh}_G(Y,0)$ and the category of equivariant coherent sheaves $\mathrm{Coh}_{G}(Y \times \AA^1[-1])$ where $G$ acts on $\AA^1[-1]$ by $\chi$.
\end{proposition}
\begin{proof}
An object of $\mathrm{Coh}_{G}(Y \times \AA^1[-1])$ consists of a complex $\scr{F}^{\bullet}$ of $G$-equivariant coherent sheaves on $Y$, along with a degree-$(-1)$ $G$-equivariant map  $t:\scr{F}^{\bullet} \to \scr{F}^{\bullet-1}(\chi)$ such that $t^2 = 0$: this is the action of the weight-$\chi$ coordinate $t$ of $\AA^1[-1] = \mathrm{Spec}\CC[t]$, $|t| = -1$ on $\scr{F}$. Given such an object, we produce a matrix factorization of $0$ with:
\begin{align*}
    F_0 & = \bigosum_{n \in \ZZ} \scr{F}^{2n}(\chi^{-n}) & F_1 & = \bigosum_{m \in \ZZ} \scr{F}^{2m+1}(\chi^{-m})
\end{align*}
and $\varphi_1 = t: F_1 \to F_0$, $\varphi_0 = t: F_0 \to F_1(\chi)$. Similarly, a degree-$2k$ morphism $f^{\bullet}: \scr{F}^{\bullet} \to \scr{G}^{\bullet}$ yields morphisms of $G$-equivariant coherent sheaves on $Y$ via $f_0: F_0 \to G_0(\chi^{k})$ and $f_1: F_1 \to G_1(\chi^{k})$, commuting with $\varphi_F, \varphi_G$, and thus a degree-$2k$ morphism in $\mathrm{Coh}_G(Y,0)$ (and likewise for degree $2k+1$). As above, it is not difficult to check that this gives an  equivalence of dg categories.
\end{proof}

\begin{proposition}\label{prop:descent}
Suppose $G$ is an affine algebraic group, let $\mathrm{Sch}_{G}$ denote the category of finite-type schemes over $\CC$ equipped with a $G$-action, with equivariant morphisms between them. Then the functor $\mathrm{Coh}_{G}(\; \cdot \;, 0): \mathrm{Sch}_{G} \to \mathrm{Cat}_{dg}$ satisfies descent with respect to proper equivariant morphisms of schemes with $G$-action.
\end{proposition}
\begin{proof}
Because of Proposition \ref{prop:mf_to_coh}, we can apply the descent result of \cite[Proposition 4.7.2.2]{GRI} and then pass to compact objects. Though this is stated for aft schemes, the statement that we require for stacks follows by descent.
\end{proof}

% We adapt the Barr-Beck-Lurie descent argument of \cite[Proposition 7.2.2]{GRI}: one may define componentwise an analogous colimit-preserving functor $f^{!}$ on ind-matrix factorizations that is right adjoint to the pushforward functor $f_{\ast}$ of \cite[\S 2.3.1]{Hirano} for $f$ a proper map. When $f$ is smooth, $f^{!}$ preserves compact objects; when $f$ is surjective, $f^{!}$ will also be conservative: checking these properties amounts to checking them componentwise on the sheaves $F_0, F_1$ since $W=0$.

%Compare \cite[\S 3]{pascaleff}. 
Finally, we may prove Proposition \ref{prop:compatibility}:
\begin{proof}[Proof of Proposition \ref{prop:compatibility}]
We will show that $\mathrm{Coh}_{\mathbb{G}_m}(\tilde{S}, 0)$ can be written as the pushout of the same diagram giving $\mathrm{Coh}(\tilde{Z})$ (see Proposition \ref{prop:derived_pushout}). We know from Lemma \ref{lem:trivial_action} and Corollary \ref{cor:divisor-factorization} that we can write each of the three categories in the pushout diagram as matrix factorization categories:
\begin{itemize}
    \item $\mathrm{Coh}(D) \simeq \mathrm{Coh}_{\GG_m}(D,0)$;
    \item $\mathrm{Coh}(D \times \AA^1[-1]) \simeq \mathrm{Coh}_{\GG_m}(D \times \AA^1,0)$;
    \item $\mathrm{Coh}(X) \simeq \mathrm{Coh}_{\GG_m}(X,0)$;
\end{itemize}
with the corresponding functors in the diagram given by the pushforward functors on matrix factorizations. Then by Proposition \ref{prop:descent},
\begin{center}
    \begin{tikzcd}
\mathrm{Coh}_{\GG_m}(D,0) \ar{r} \ar{d} & \mathrm{Coh}_{\GG_m}(X,0) \ar{d}\\
 \mathrm{Coh}_{\GG_m}(D \times \mathbb{A}^1,0) \ar{r} & \mathrm{Coh}_{\GG_m}(\tilde{S},0)
\end{tikzcd}
\end{center}
is a pushout diagram, where the arrows are given by pushforwards.
\end{proof}

Combining Propositions \ref{prop:knorrer1}, \ref{prop:knorrer2}, \ref{prop:compatibility} we have:

\begin{theorem}\label{prop:orlov}(Kn\"orrer periodicity) Suppose $\mathbb{G}_m$ acts trivially on $K_X$ and with weight $1$ on the $\AA^1$ fibers of $K_X \times \AA^1$, then we have equivalences of dg categories:
\begin{center}
\begin{tikzcd}
\mathrm{Coh}(Z) \ar{r} & \mathrm{Coh}_{\mathbb{G}_m}(K_{X} \times \mathbb{A}^1, s(x)yz) & \ar{l} \mathrm{Coh}(\tilde{Z})
\end{tikzcd}
\end{center}
\end{theorem}

\subsection{Bigradings}

In this section, we will consider the extension of the above results to graded dg categories: i.e., categories linear over the category $\Coh(B\mathbb{G}_m)$ of graded vector spaces. We shall refer to grading-preserving functors among graded dg categories as \textit{graded dg functors}.

%These are categories where the morphism complexes have a $\ZZ \times \ZZ$-grading; these gradings do not necessarily play a symmetric role: the differential on the complex need not have the same degree for both.
%I thought this might be helpful to the reader but please take it out if it's not strictly correct

Suppose we consider the $\GG_m \times \GG_m$-action on $K_X \times \AA^1$ where the first $\GG_m$ acts by weight $1$ on the $\AA^1$ factor (and trivially on $K_X$), and the second $\GG_m$ acts with weight $1$ on the fibers of $K_X$ (and trivially on $X \times \AA^1)$. We consider $Z, \tilde{S}$ to carry the induced $\GG_m \times \GG_m$ actions. The weights for the action of the first factor $\GG_m$ will give the homological grading, while the weights for the action of the second second factor will give the auxilliary grading.

\begin{theorem}\label{thm:bigraded_orlov}
There are equivalences of graded dg categories:
\begin{equation*}
    \mathrm{Coh}_{\mathbb{G}_m \times \GG_m}(Z, 0) \xrightarrow{j_{Z\ast}p_{Z}^{\ast}} \mathrm{Coh}_{\mathbb{G}_m \times \GG_m}(K_X \times \mathbb{A}^1, s(x) y z) \xleftarrow{j_{\tilde{S}\ast}p_{\tilde{S}}^{\ast}} \mathrm{Coh}_{\mathbb{G}_m \times \GG_m}(\tilde{S}, 0)
\end{equation*}
induced by the same Kn\"orrer periodicity functors as above.
\end{theorem}
\begin{proof}
Follows from \cite[Theorem 1.2]{Hirano} with $G = \GG_m \times \GG_m$.
\end{proof}

The next result follows similarly from Proposition \ref{prop:descent} by taking $G = \GG_m \times \GG_m$.
\begin{proposition}\label{prop:bigraded_descent}
We have a pushout diagram of graded dg categories:
\begin{center}
    \begin{tikzcd}
\mathrm{Coh}_{\GG_{m}^2}(D,0) \ar{r} \ar{d} & \mathrm{Coh}_{\GG_{m}^2}(X,0) \ar{d}\\
 \mathrm{Coh}_{\GG_{m}^2}(D \times \mathbb{A}^1,0) \ar{r} & \mathrm{Coh}_{\GG_{m}^2}(\tilde{S},0)
\end{tikzcd}
\end{center}
where the second $\GG_m$ factor acts trivially and the first acts with weight $1$ on the $\AA^1$ factor.
\end{proposition}

We will also need an upgraded form of Lemma \ref{lem:trivial_action}:

\begin{lemma}\label{lem:trivial_action2} 
If $Y$ is any smooth algebraic variety with a $\GG_m$ action, then the category $\mathrm{Coh}_{\GG_m}(Y)$ of equivariant coherent sheaves is equivalent to $\mathrm{Coh}_{\GG_m^2}(Y,0)$ (where the second $\GG_m$-factor acts trivially) as a graded dg category.
\end{lemma}

The equivalences between various categories of coherent sheaves constructed in this section (Propositions \ref{prop:orlov} and \ref{thm:bigraded_orlov}) shall be used later in \S 4, where they arise as mirrors of natural automorphisms of wrapped Fukaya categories of very affine hypersurfaces.

\section{Gluing Diagrams}\label{sec:gluing}

\subsection{Liouville sectors for Laurent polynomials}\label{subsec:gluing}

In order to turn the Landau-Ginzburg model $(\cstar{n}, f)$ into a Liouville sector, we shall apply a modified form of the construction from \cite[Proposition 1]{periodicity} to the function $1-f$. The key difference is that in \cite{periodicity}, the holomorphic function in question was assumed to have a single critical value, and the associated Liouville sector lived over a small neighborhood of this critical value, whereas in our setting $1-f$ will in general have more than one critical value.

To outline the general construction, suppose first that $X$ is a smooth affine variety with an embedding $i: X \to \CC^N$; then $X$ becomes a Stein manifold with the Stein function $\phi: X \to \RR$ given by the restriction of $\phi(z) = |z|^2$ on $\CC^N$. Suppose that $W: X \to \CC$ is the restriction of a polynomial on $\CC^N$ to $X$.

\begin{definition}
A point $p \in \RR$ is called a \textbf{non-Malgrange point} of $|W|^2:X \to \RR$ if the Malgrange condition for $|W|^2$ fails at $p$. Likewise, a point $p \in \CC$ is a non-Malgrange point of $W:X \to \CC$ if the (complex) Malgrange condition fails there.
\end{definition}

Recall that the Malgrange condition at $p \in \RR$ says that there exists $R, \epsilon, \eta>0$ so that if $|z|>R$ and $||W|^2(z) - p| \leq \epsilon$ then
\begin{equation*}
    |z| |\grad_X |W|^2| > \eta
.\end{equation*}
and likewise in the complex case. 

It is a well-known result that the number of (real or complex) non-Malgrange points of a polynomial mapping is finite (for instance see \cite[Remark 3]{Malgrange} and take the intersection with the algebraic variety $X$). So choose $R > 0$ sufficiently large so that all of the critical points and non-Malgrange points of $W$ are contained inside $|W|^2 < R$. If we define a (simple Liouville) deformation of the Stein structure $\phi$ on $X$ to
\begin{equation*}
    \psi(z) = \phi(z) + D\phi(z)|W|^{2m}
\end{equation*}
for $D>0$ a positive constant and $m$ a positive integer. Then
\begin{lemma}
The Liouville vector-field of $\psi$ is outward pointing along $|W|^2 = R$ for $R, D, m>0$ sufficiently large.
\end{lemma}
\begin{proof}
The proof follows \cite[Proposition 1]{periodicity}: we need to verify that the inner product
\begin{equation*}
    \langle \grad_X |W|^2, \grad_X \psi \rangle
\end{equation*} 
is strictly positive everywhere on $|W|^2 = R$. Rewriting this inequality and applying the Cauchy-Schwarz inequality gives the condition
\begin{equation}\label{eqn:inequality}
    2\br{ \frac{1}{D |W|^{2m-2}} + |W|^2 } < m |\phi|^{1/2}|\grad_X |W|^2|
.\end{equation} 
When the Malgrange condition for $|W|^2$ holds at $R$ then we have a constant $\eta>0$ such that
\begin{equation*}
    |z| |\grad_X |W|^2| > \eta
\end{equation*}
on $|W|^2 = R$; if $R$ is sufficiently large so that no critical points of $W$ occur on $|W|^2 = R$ then $|\grad_X |W|^2|>0$ there and hence there is a constant $C>0$ (depending on $R$) such that
\begin{equation*}
    |W|^2 < C |z| |\grad_X |W|^2|
.\end{equation*}
and hence the inequality \ref{eqn:inequality} holds for $D, m>0$ sufficiently large.
\end{proof}

Similarly, we shall need to know that the parallel transport of $W$ is well-defined outside of a disk:
\begin{lemma}
Symplectic parallel transport gives exact symplectomorphisms between smooth fibers of $W$ over $|W|^2 = R$ for some $R>0$ sufficiently large.
\end{lemma}
Again the proof is analogous to \cite[Lemma 1]{periodicity}, where we take $R>0$ sufficiently large so that none of the critical points or non-Malgrange points of $W$ lie on $|W|^2 = R$.

We will want to have the freedom to use a deform the Liouville structure to compute wrapped Fukaya categories. We now record here a standard lemma we shall use implicitly throughout:

\begin{definition}\label{def:simple}
Suppose $X$ is a Liouville manifold with Liouville form $\lambda$; a smooth family $\lambda_t$, $t \in [0,1]$ of Liouville forms for $X$ with $\lambda_0=\lambda$ such that the  union of the skeleta of the Liouville structures $\lambda_t$ stays within a compact set, is called a \textbf{simple Liouville deformation}.
\end{definition}

For instance, this condition is satisfied by a family of Weinstein functions whose critical points remain inside a fixed compact set.

\begin{lemma}
Given a Liouville manifold $X$ and a simple Liouville deformation $\lambda_t$, all of the Liouville manifolds $(X, \lambda_t)$ are exact symplectomorphic and the wrapped Fukaya categories $\scr{W}(X,\lambda_t)$ are all quasiequivalent.
\end{lemma}

This means we can deform the Liouville structure when considering only the Fukaya category up to quasi-equivalence (cf. \cite[\S 2]{periodicity}).

\begin{proposition}\label{prop:swappable}
There is a simple Liouville deformation of $(X, \phi)$ to the Stein manifold $(X, \psi)$ for which the subset $\set{|W|^2 \leq R_1} \cap \set{\psi \leq R_2}$ has convex boundary, for suitably large $R_1, R_2>0$. Moreover, all of the fibers $W\inv(z) \cap \set{\psi \leq R_2}$ are exact symplectomorphic via symplectic parallel transport along $|W|^2(z) = R_1$; and they remain exact symplectomorphic after rounding the corners to give a Liouville domain $X_0$.
\end{proposition}

By making the rounding parameter sufficiently small, the contact boundary of the Liouville domain $X_0$ contains the Liouville hypersurface $F = W\inv(R_1) \cap \set{\psi \leq R_2-\epsilon}$ with convex boundary; and we may also take $\epsilon>0$ sufficiently small so that the completion of $F$ is exact symplectomorphic to $W\inv(R_1)$. The pair $(X_0, F)$ gives a sutured Liouville domain that is independent of sufficiently large $R_1, R_2$. By \cite[Lemma 2.32]{GPS1}, there is a (homotopically) unique Liouville sector associated to a sutured Liouville domain, which we will denote by $(X,W)$.

\begin{definition}\label{def:general_sector}
The Liouville sector $(X,W)$ associated to a pair $X, W$ of an affine variety and a polynomial function $W:X \to \CC$ is the Liouville sector associated to the (rounding of the) sutured Liouville domain 
$(\set{|W|^2 \leq R_1} \cap \set{\psi \leq R_2},  W\inv(R_1) \cap \set{\psi \leq R_2-\epsilon})$ for sufficiently large $R_1, R_2>0$ and small $\epsilon>0$.
\end{definition}

We observe here that there are two constructions for turning a sutured Liouville domain into a Liouville sector. One is the `completion' described in \cite[Definition 2.14]{GPS1}. Alternatively, applying \cite[Prop 2.6]{sylvan_lemma} to $F = W\inv(R_1) \cap \set{\psi \leq R_2 - \epsilon}$ gives us for some $\rho>0$ (after possibly increasing $R_1, R_2$) a standard Liouville neighborhood $U$ of $F$ of the form $(\hat{F} \times \set{z \in \CC \;\;: \mathrm{Re}(z) > - \rho}, \lambda_{X}|_F + \lambda_{\CC}^{\mathrm{std}})$, properly embedded inside the completion $\hat{X_0}$. Then the hypersurface $\set{\mathrm{Re}(z) = 0}$ is parallel to the Liouville vector field at infinity and has a defining function $I = \mathrm{Im}(z)$ (see Figure \ref{fig:step1}) in the sense of \cite[Definition 2.4]{GPS1}. This gives $\set{\mathrm{Re}(z) \leq 0}$ the structure of a Liouville sector. By \cite[Lemma 2.32]{GPS1}, there is a (homotopically unique) homotopy of Liouville sectors from $\set{\mathrm{Re}(z) \leq 0} \subset \hat{X}_0$ to the convex completion of $(X_0, F)$ (further inspection of the argument shows that this can be taken to be a \textit{simple} Liouville homotopy). 

Returning to the situation at hand, it is known that by taking $T>0$ sufficiently small, all of the critical values of $f$ can be made to lie in an arbitrarily small neighbourhood of $1$ under our hypothesis (2) (see for instance \cite[Lemma 5.2]{Speculations}). Though \cite{Speculations} assumes hypothesis (2) that the maximal cells are unimodular simplices, the proof of \cite[Lemma 5.2]{Speculations} only uses the weaker hypothesis that they are simplices, so by Remark \ref{rem:unimodular}, we expect that arguments will continue to apply in the setting of toric Deligne-Mumford stacks. Now we may apply the arguments above from \cite{periodicity} to $W=1-f$, to see that $\set{|f-1| \leq R_1} \cap \set{\psi \leq R_2}$ has convex boundary. We summarize this as:

\begin{definition}\label{def:sector}
The Liouville sector $(\cstar{n}, f)$ is defined by applying Definition \ref{def:general_sector} with $R_1>1$ to the function $1-f$ on $(\CC^{\ast})^n$ with its standard Stein structure $\phi$. The hypersurface $F$ we take is a Liouville domain whose completion is isotopic to the reference fiber $H = f\inv(0)$.
\end{definition}

Let $\gamma: [0,1] \to \CC$ be the smooth arc $\gamma(t)  = R_1 \e{i\theta(t)}$ where $\theta(t) = \delta + t (2 \pi - \delta)$ for some small $\delta>0$. We identify each fiber $F_z = W\inv(z)$ with our reference fiber $H = f\inv(0) = W\inv(1)$ via parallel transport along a ray from $1$ to $z$ (or a small perturbation thereof if these rays happen to cross a non-Malgrange point of $W$). Under this identification, the counterclockwise global monodromy $\mu: H \to H$ is defined to be the parallel transport symplectomorphism $\Phi_{\gamma}: F_{R_1\e{i\delta}} \to F_{R_1}$ for $W$ with $\delta>0$ sufficiently small.

\subsection{Compatibility of Liouville structures}\label{sec:comparisons}

In this section, we verify that several (a priori different) methods for obtaining a Liouville structure on the complement $(\CC^{\ast})^n \setminus H$ yield (simple) homotopic results. 

Firstly, note that the Stein structure on $(\CC^{\ast})^n \setminus H$, considered as a very affine hypersurface, comes from the embedding into $(\CC^{\ast})^{n+1}$ via $\iota(z) = (-1/f,z) \in (\CC^{\ast})^n \times \CC^{\ast}$; if $\phi$ denotes the standard Stein function on $(\CC^{\ast})^n$ then the induced Stein function on $(\CC^{\ast})^n \setminus H \subset (\CC^{\ast})^n$ is
\begin{equation*}
    \psi_{0}(z) = \phi(z) + C\phi(z) (\log|f(z)|)^{2n}
\end{equation*}
for some sufficiently large positive constant $C$ and positive integer $n$: the same argument as in \cite[Proposition 3]{periodicity} implies that this gives a Liouville structure on $(\CC^{\ast})^n \setminus H$. Under the deformation of the Stein structure on $(\CC^{\ast})^n$ in \S \ref{sec:gluing} this Stein function becomes
\begin{equation*}
    \psi_{1}(z) = \phi(z) + C\phi(z)(\log|f(z)|)^{2n} + D \phi(z) |1-f|^{2m}
\end{equation*}
for $D>0$ a sufficiently large positive constant and $m$ a sufficiently large positive integer. We call the Liouville form induced by $\psi_1$ the \textit{hypersurface Liouville structure} on $(\CC^{\ast})^n \setminus H$. 

Observe that by setting $f_t(z) = f(z) - t$, we have a family of Stein functions
\begin{equation*}
\psi_{1,t}(z) = \phi(z) + C\phi(z)(\log|f_t(z)|)^{2n} + D \phi(z) |1-f|^{2m} 
\end{equation*}
inducing a simple Liouville homotopy between $(\CC^{\ast})^n \setminus H$ with $H = f^{-1}(0)$ and $(\CC^{\ast})^n \setminus \hat{F}$ with $\hat{F} = W^{-1}(R_1)$ with their hypersurface Liouville structures (so long as $t$ does not pass through non-Malgrange points of $f$).

On the other hand, applying Sylvan's construction \cite[Proposition 2.6]{sylvan_lemma} to the sutured Liouville domain obtained from $X = (\CC^{\ast})^n, W = 1-f$ as in Definition \ref{def:general_sector} gives us for some $\rho>0$ a standard Liouville neighborhood $U$ of $F = W^{-1}(R_1) \cap \{\psi \leq R_2\}$ of the form $\hat{F} \times \set{z \in \CC \;\;: \mathrm{Re}(z) > - \rho}$, properly embedded inside $(\CC^{\ast})^n$. This means that there is a smooth fibration $g: U \to \CC$ such that $g^{-1}(0) = F$ and $\lambda_{X}|_{U} = \lambda_{X}|_{F} + g^{\ast} \lambda_{\CC} + dh$ for some compactly supported function $h$. Sylvan then defines a family of cutoff functions $\kappa_t$ and a family of Liouville forms $\lambda_{t} = \lambda_{F} + g^{\ast} \lambda_{\CC} + d(\kappa_t h)$ and shows this is a simple Liouville deformation to $\lambda_0 = \lambda_F + g^{\ast} \lambda_{\CC}$. Then the hypersurface $\{\mathrm{Re}(g) = 0 \} \subset U$ is parallel to the Liouville vector field of $\lambda_0$ at infinity and is a sectorial hypersurface with $I$-function given by $I = \mathrm{Im}(g)$.

Now we wish to do this in the punctured case. Let $\beta_{s}:\RR \to [0,1]$ be a smooth cutoff function with:
\begin{align*}
    \beta_{s}(r) = \left\{ \begin{matrix}
        1 && \text{for $|r| \leq s/2$} \\
        0 && \text{for $|r| > s$}
    \end{matrix}\right.
\end{align*}

\begin{lemma}\label{lem:one-dimensional}
Consider $\CC^{\ast}$ with equipped with the function
\begin{equation*}
    \psi = C |z|^2 + \beta_{s}(|z|) (\log|z|)^2
\end{equation*}
where $C>0$ is a constant and $s>0$; then for $C$ sufficiently large (depending on $s$), $\psi$ is a Stein function for $\CC^{\ast}$ that is simple Liouville deformation equivalent to the standard Liouville structure $\lambda_0$ on $\CC^{\ast}$. 
\end{lemma}
\begin{proof}
Clearly $\psi$ is exhausting for $s>0$ so to show $\psi$ is a Stein function with respect to the standard complex structure on $\CC^{\ast}$ it suffices to show that $\Delta \psi > 0$: since $\psi$ only depends on the radial coordinate $r$ we can write
\begin{align*}
    \Delta \psi & = r^{-1} \partial_r \psi  + \partial^{2}_r \psi\\
    & = 4 C + 2 r^{-2} \beta_{s}(r) + \text{terms involving $\beta_{s}\dash(r), \beta_{s}\ddash(r)$}
    % & = r^{-1}(2 C_1 r + C_2 \beta_{s}\dash(r) \log^2(r) + 2 C_2 r^{-1} \beta_{s}(r) \log r)  \\ & + (2 C_1 + C_2 \beta_{s}\ddash(r) \log^2(r) + 4 C_2 r^{-1} \beta_{s}\dash(r) \log(r)  - 2 C_2 r^{-2} \beta_{s}(r) \log r + 2 C_2 r^{-2} \beta_{s}(r) \\
    % & = 4 C_1 + 2 C_2 r^{-2} \beta_{s}(r) + C_2 \beta_{s}\dash(r) \left( r^{-1} \log^2 r + 4 r^{-1} \log(r) \right) + C_2 \beta_{s}\ddash(r) \log^2(r)
\end{align*}
since
\begin{equation*}
    (r^{-1}\partial_r + \partial_{r}^2)\log^2(r) = 2 r^{-2} \log(r) - 2 r^{-2} \log(r) + 2 r^{-2} = 2 r^{-2}
\end{equation*}
Since $\beta_{s}$ is a cutoff function, we have $\beta_{s}(r) = 0, \beta_{s}\dash(r) = \beta_{s}\ddash(r) = 0$ for all $r$ sufficiently large, and $\beta_{s}(r) = 1, \beta_{s}\dash(r) = \beta_{s}\ddash(r) = 0$ for all $r$ sufficiently small. Therefore, by taking $C > 0$ sufficiently large, $\Delta \psi$ will be strictly positive. Moreover, as $s \to \infty$ the value of this constant $C$ required does not increase. To see that this deformation is simple it suffices to observe that
\begin{equation*}
    \partial_r \psi = 2 C r + \beta_{s}\dash(r) \log^2(r) + 2 r^{-1} \beta_{s}(r) \log(r) = 0
\end{equation*}
is only possible if $r<1$ or $s/2 < r < s$, otherwise the sum must be strictly positive; by taking $C$ larger if necessary we can assume that for all $s>1$, we have $2 C r > \log^2(r)$ for $s/2<r<s$. Therefore the only critical points of $\psi$ can occur where $r<1$ and hence the deformation is simple as $s \to \infty$.
\end{proof}

Let $\lambda_s = -d^c \psi$ denote the Liouville form on $\CC^{\ast}$ in Lemma \ref{lem:one-dimensional} above. We now define  a Liouville form on $(\CC^{\ast})^n \setminus F$ by extending the $1$-form on $U \setminus g^{-1}(0)  = \hat{F} \times \set{z \in \CC^{\ast} \;\;: \mathrm{Re}(z) > - \rho}$ given by
\begin{equation*}
    \lambda_1 = \lambda_{F} + g^{\ast} \lambda_{s} + dh
\end{equation*}
for $0 < s < \rho$ sufficiently small. This gives a new Liouville structure on $(\CC^{\ast})^n \setminus F$ which we call the \textit{punctured Liouville structure}.

\begin{lemma}\label{lem:deformation1}
The hypersurface Liouville structure on $(\CC^{\ast})^n \setminus F$ (via its embedding in $(\CC^{\ast})^{n+1}$ as a hypersurface) is simple Liouville homotopic to the punctured Liouville structure on $(\CC^{\ast})^n \setminus F$ (obtained by applying Sylvan's construction to $(\CC^{\ast})^n$).
\end{lemma}

\begin{proof}
This follows from Lemma \ref{lem:one-dimensional} by noting that there is an isotopy between $g: U \to \CC$ and $f: U \to \CC$ on $U$. 
\end{proof}

\begin{lemma}
The hypersurface Liouville structure on $H = f^{-1}(0)$ obtained via its embedding as a hypersurface in $(\CC^{\ast})^n$ is simple Liouville homotopic to the Liouville structure on the semitropicalization $H_{s,t} = f_{s,t}^{-1}(0)$ used by \cite{BenVivek,Zhou}.
\end{lemma}

See \cite[Theorem 6.3]{vandenbergh}. We may then apply this to $(\CC^{\ast})^n \setminus H$ considered as a hypersurface in $(\CC^{\ast})^{n+1}$ to see that:

\begin{corollary}
The hypersurface Liouville structure obtained on $(\CC^{\ast})^n \setminus H$ via its embedding in $(\CC^{\ast})^{n+1}$ as a hypersurface, is simple Liouville homotopic to the semitropical Liouville structure on the semitropicalizations $\tilde{f}_{\infty,s,t}^{-1}(0)$ and $\tilde{f}_{0,s,t}^{-1}(0)$.
\end{corollary}

\begin{corollary}
The Liouville sector associated to $((\CC^{\ast})^n, f)$ in Definition \ref{def:sector} is simple Liouville homotopic to the Liouville sector associated to the pair $((\CC^{\ast})^n, f)$ in \cite[Corollary 6.2.6]{BenVivek}. 
\end{corollary}

\subsection{Sectorial Gluing}

\begin{figure} %Step 1
\begin{center}
    \begin{tikzpicture}[scale=0.8]
%\draw[help lines] (-8,-8) grid (8,8);
\fill [fill=gray!5] (0,0) circle (5);
\fill [fill=gray!10] (0,0) circle (4);
\draw (0,0) circle (4);

\fill[fill=green!10] (-23:9) to[out=175, in=-100] (-20:4.5) arc[start angle=-20, end angle=20,radius=4.5] (20:4.5)  to[out=100,in=177] (23:9) to (-23:9);

\draw[dashed] (0,0) circle (5);
\fill[black] (5,0) circle (0.1);
\draw[thick] (-10:5) arc[start angle=-10, end angle=10,radius=5];
\draw[thick] (11:5) to[out=100,in=180] (15:8);
\draw[thick] (-10:5) to[out=-100, in=180] (-15:8);

\draw[dashed] (-20:4.5) arc[start angle=-20, end angle=20,radius=4.5];
\draw[dashed] (21:4.5) to[out=100,in=180] (25:8);
\draw[dashed] (-20:4.5) to[out=-100, in=180] (-25:8);

\fill[black] (1.5,0) circle (0.1);
\node[anchor=north west] at (1.5,0) {$H$};
\draw (-0.5,0) node[left] {$\huge{\times}$};
    \draw (0.5,-0.5) node[left] {$\huge{\times}$};
    \draw (0.5,0.5) node[left] {$\huge{\times}$};
    \node[anchor=south east] at (-0.5,0.5) {$\mathrm{crit}(f)$};
    \node[anchor=west] at (5,0) {$F_{R_1}$};
\node at (7,1) {$U$};
\node[anchor=north] at (7,-2) {$\mathrm{Re}(z)=0$};
\node at (-4,4) {$\cstar{n}$};
\node[anchor=north] at (-4,-4) {$|W| = R_1$};
    \end{tikzpicture}
\end{center}
\caption{From proof of Proposition \ref{prop:sector_gluing}: $U$ (in green) is a product neighborhood of the fiber $F_{R_1}$ (black dot)}
\label{fig:step1}
\end{figure}
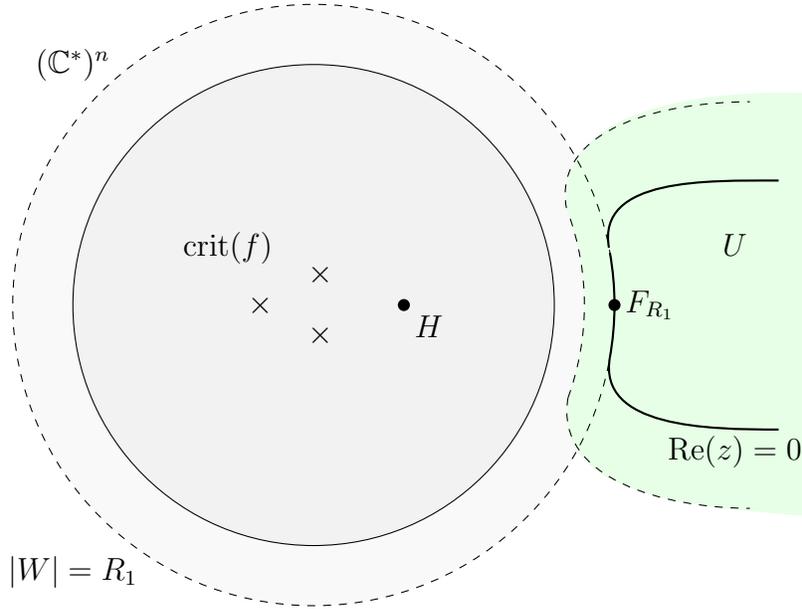

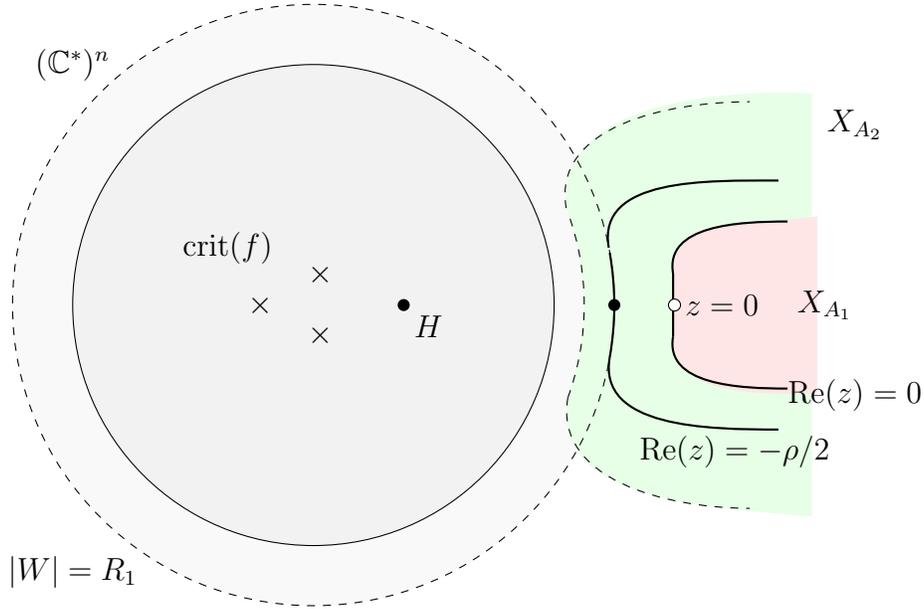
\begin{figure}%Step 2
\begin{center}
    \begin{tikzpicture}[scale=0.8]
%\draw[help lines] (-8,-8) grid (8,8);
\fill [fill=gray!5] (0,0) circle (5);
\fill [fill=gray!10] (0,0) circle (4);
\draw (0,0) circle (4);

\fill[fill=green!10] (-23:9) to[out=175, in=-100] (-20:4.5) arc[start angle=-20, end angle=20,radius=4.5] (20:4.5)  to[out=100,in=177] (23:9) to (-23:9);

\fill[fill=red!10] (5:6) to (-5:6) to[out=-90, in=180] (-10:8.5) to (10:8.5) to[out=190,in=90] (5:6);

\draw[dashed] (0,0) circle (5);
\fill[black] (5,0) circle (0.1);
\draw[thick] (-10:5) arc[start angle=-10, end angle=10,radius=5];
\draw[thick] (11:5) to[out=100,in=180] (15:8);
\draw[thick] (-10:5) to[out=-100, in=180] (-15:8);

\draw[dashed] (-20:4.5) arc[start angle=-20, end angle=20,radius=4.5];
\draw[dashed] (21:4.5) to[out=100,in=180] (25:8);
\draw[dashed] (-20:4.5) to[out=-100, in=180] (-25:8);

\draw[thick] (-5:6) to (5:6);
\draw[thick] (5:6) to[out=100,in=180] (10:8);
\draw[thick] (-5:6) to[out=-100, in=180] (-10:8);
\draw[fill=white] (0:6) circle (0.1);

\fill[black] (1.5,0) circle (0.1);
\node[anchor=north west] at (1.5,0) {$H$};
\draw (-0.5,0) node[left] {$\huge{\times}$};
    \draw (0.5,-0.5) node[left] {$\huge{\times}$};
    \draw (0.5,0.5) node[left] {$\huge{\times}$};
    \node[anchor=south east] at (-0.5,0.5) {$\mathrm{crit}(f)$};
    \node[anchor=west] at (6,0) {$z=0$};
\node at (9,-1.5) {$\mathrm{Re}(z)=0$};
\node at (9,3) {$X_{A_2}$};
\node at (8.5,0) {$X_{A_1}$};
\node[anchor=north] at (7,-2) {$\mathrm{Re}(z)=-\rho/2$};
\node at (-4,4) {$\cstar{n}$};
\node[anchor=north] at (-4,-4) {$|W| = R_1$};
    \end{tikzpicture}
    \caption{From proof of Proposition \ref{prop:A2}: $X_{A_2}$ is in green; $X_{A_1}$ is in red.}
    \label{fig:step2}
\end{center}
\end{figure}

\begin{proposition}\label{prop:sector_gluing}
The Liouville manifold $\comp$ is simple Liouville homotopic to the sectorial gluing of $(H \times \CC^{\ast},z)$ along $H$ to $( (\CC^{\ast})^n, f)$ 
\end{proposition}
\begin{proof}
As above, applying \cite[Prop 2.6]{sylvan_lemma} to $F_{R_1}\cap\set{\phi \leq R_2}$ gives us (after a small deformation explained in \cite{sylvan_lemma}) a standard Liouville neighborhood $U$ of the form:
$$(H \times \set{z \in \CC \;\;: \mathrm{Re}(z) > - \rho}, \lambda_{\cstar{n}}|_F + \lambda_{\CC}^{\mathrm{std}})$$
properly embedded inside the completion $\cstar{n}$. Here, we reiterate our choice of identification of $F_{R_1}$ with $H$ via parallel transport along a ray. The sectorial hypersurface $\set{\mathrm{Re}(z) = 0}$ thus exhibits $(\CC^{\ast})^n$ as the sectorial gluing of $(H \times \CC,z)$ along $H$ to $((\CC^{\ast})^n, f)$ (see Figure \ref{fig:step1}).

Next, we introduce a Liouville form on $H \times \set{z \in \CC \;\;: \mathrm{Re}(z) > - \rho}\setminus\set{\rho}$ so that the difference with the original Liouville form is supported on a small open neighbourhood of $H \times \set{\rho}$, as in \cite[Proposition 6]{periodicity}. Then, the hypersurface $\set{\mathrm{Re}(z)=0}$ is still parallel to the Liouville vector field at infinity, and so we may glue $H \times \set{z \in \CC \;\;: \mathrm{Re}(z) > - \rho}\setminus\set{\rho}$ back to (a sector canonically identified with) $((\CC^{\ast})^n, f)$ along $H$. But by the discussion from \S \ref{sec:comparisons}, this is simple Liouville homotopic to the hypersurface Liouville form on $\comp$. This demonstrates that the gluing of $(H \times \CC^{\ast},z)$ along $H$ to $( (\CC^{\ast})^n, f)$ yields $\comp$ (after homotopy): see Figure \ref{fig:step2}.
\end{proof}

\begin{remark}
We should think of this as the sectorial decomposition of $\comp$ associated to the presentation as $\tilde{f}_0=0$.
\end{remark}

Recall the $A_2$ sector from \cite[Remark 13.1]{GPS2}.

\begin{figure}
\begin{center}
\begin{tikzpicture}[scale=1.2]
\begin{axis}[
    xmin = -2, xmax = 2,
    ymin = -2, ymax = 2,
    zmin = 0, zmax = 1,
    axis equal image,
    yticklabels={,,},
    xticklabels={,,},
    view = {0}{90},
    axis line style={draw=none},
]
    \addplot3[
        quiver = {
            u = {(x*(-4 + 4*x^2 - 2/(x^2 + y^2+0.01)^2))*(x^2+y^2+0.1)/100},
            v = {(4*y^3 - (2*y)/(x^2 + y^2+0.01)^2)*(x^2+y^2+0.1)/100},
        },
        -stealth,
        domain = -2:2,
        domain y = -2:2,
    ] {0};
\end{axis}
\draw[thick,dashed,blue] (2.85,0) to (2.85,5.7);
\draw[fill=white] (2.85,2.85) circle (0.1);
\fill[black] (2.85,4.25) circle (0.08);
\fill[black] (2.85,1.45) circle (0.08);
\draw[red,thick,red] (1.3,0) to (1.3,5.7);
\node[anchor=north] at (1.3,0) {$ - \rho/2$};
\node[anchor=north] at (2.85,0) {$0$};
\node[anchor=east] at (0,2.6) {$Z$};
\end{tikzpicture}
    \caption{A sketch of the deformed Liouville vector field from the proof of Proposition \ref{prop:A2}: the hypersurface $\mathrm{Re}(z) = -\rho/2$ is in red. The two black points are zeroes of the Liouville vector field.}
    \label{fig:vector_field}
\end{center}
\end{figure}
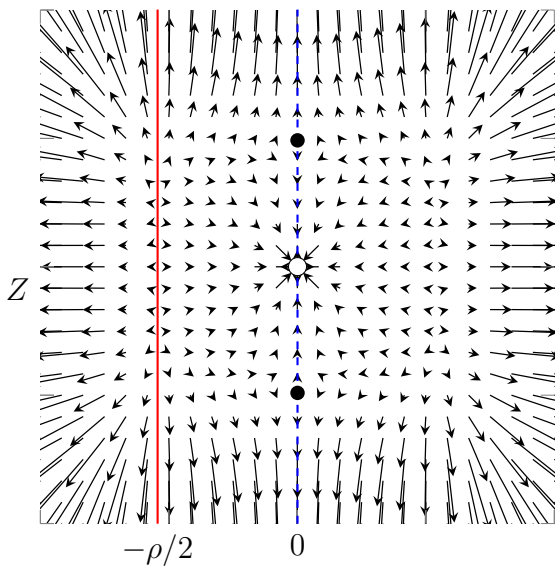

\begin{proposition}\label{prop:A2}
The Liouville manifold $\comp$ is simple Liouville homotopic to the sectorial gluing of $( (\CC^{\ast})^n, f)$ along a fiber $H$ to a sector $(E,H)$ given by gluing of $H \times A_2$ and $H \times A_1$ along two boundary components via the identifications $\mathrm{id}:H \to H$ and $\mu: H \to H$, the counterclockwise global monodromy.
\end{proposition}
\begin{proof}
We break up this proof into four steps. We observe first that this argument could be made substantially simpler by assuming the existence of an open book decomposition of the boundary of $\cstar{n}$ induced by $f$ as in \cite[Example 2.20]{GPS1}. The reader may also wish to compare with \cite[Proposition 3.7]{MilnorFibers}.

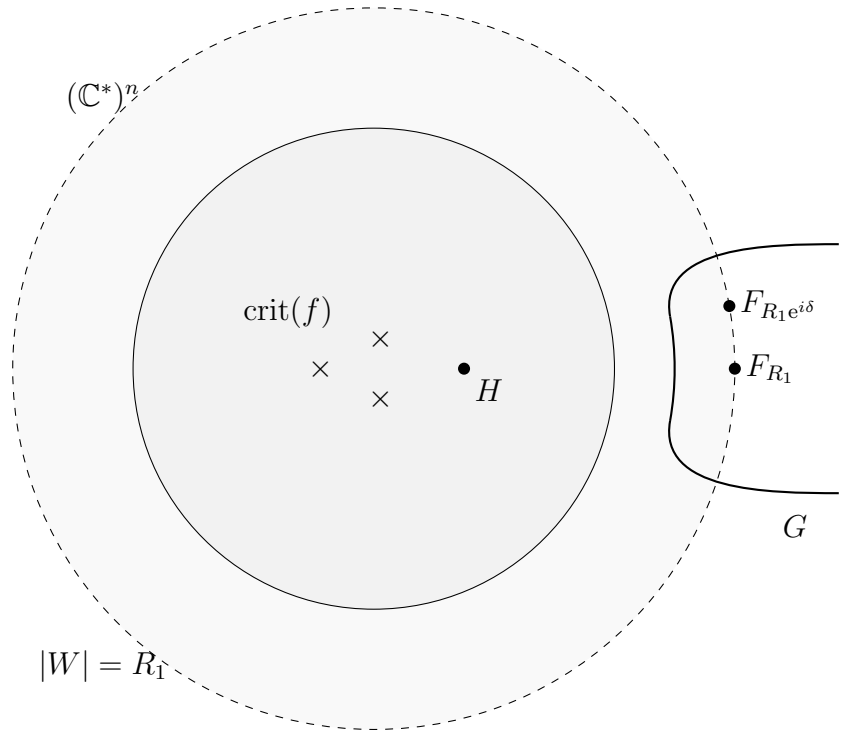
\begin{figure}%Step 3
\begin{center}
    \begin{tikzpicture}[scale=0.8]
%\draw[help lines] (-8,-8) grid (8,8);
\fill [fill=gray!5] (0,0) circle (6);
\fill [fill=gray!10] (0,0) circle (4);
\draw (0,0) circle (4);
\draw[dashed] (0,0) circle (6);
\draw[thick] (-10:5) arc[start angle=-10, end angle=10,radius=5];
\draw[thick] (10:5) to[out=100,in=180] (15:8);
\draw[thick] (-10:5) to[out=-100, in=180] (-15:8);

\fill[black] (0:6) circle (0.1);
\fill[black] (10:6) circle (0.1);

\fill[black] (1.5,0) circle (0.1);
\node[anchor=north west] at (1.5,0) {$H$};
\draw (-0.5,0) node[left] {$\huge{\times}$};
    \draw (0.5,-0.5) node[left] {$\huge{\times}$};
    \draw (0.5,0.5) node[left] {$\huge{\times}$};
    \node[anchor=south east] at (-0.5,0.5) {$\mathrm{crit}(f)$};
    \node[anchor=west] at (6,0) {$F_{R_1}$};
    \node[anchor=west] at (10:6) {$F_{R_1 \e{i\delta}}$};
\node[anchor=north] at (7,-2.25) {$G$};
\node at (-4.5,4.5) {$\cstar{n}$};
\node[anchor=north] at (-4.5,-4.5) {$|W| = R_1$};
    \end{tikzpicture}
        \caption{From proof of Proposition \ref{prop:A2}: the hypersurfaces $F_{R_1 \e{i \delta}}$ and $F_{R_1}$ and the sectorial hypersurface $G$.}
        \label{fig:step3}
\end{center}
\end{figure}

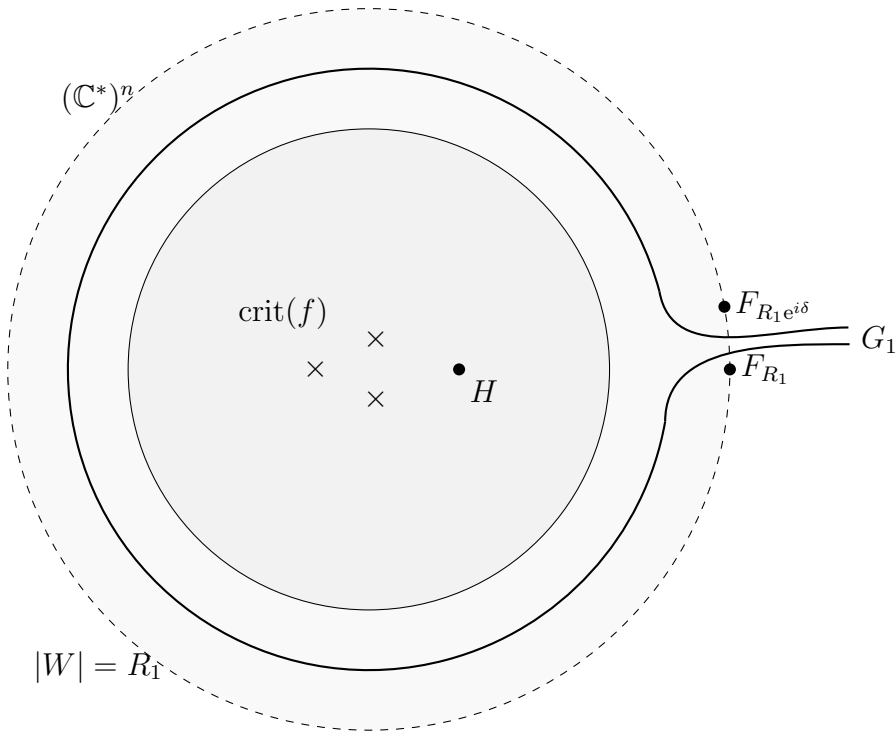
\begin{figure}%Step 4
\begin{center}
    \begin{tikzpicture}[scale=0.8]
%\draw[help lines] (-8,-8) grid (8,8);
\fill [fill=gray!5] (0,0) circle (6);
\fill [fill=gray!10] (0,0) circle (4);
\draw (0,0) circle (4);
\draw[dashed] (0,0) circle (6);
\draw[thick] (15:5) arc[start angle=15, end angle=350,radius=5];
\draw[thick] (15:5) to[out=-80,in=180] (5:8);
\draw[thick] (-10:5) to[out=90, in=180] (3:8);

\fill[black] (10:6) circle (0.1);
\fill[black] (0:6) circle (0.1);

\fill[black] (1.5,0) circle (0.1);
\node[anchor=north west] at (1.5,0) {$H$};
\draw (-0.5,0) node[left] {$\huge{\times}$};
    \draw (0.5,-0.5) node[left] {$\huge{\times}$};
    \draw (0.5,0.5) node[left] {$\huge{\times}$};
    \node[anchor=south east] at (-0.5,0.5) {$\mathrm{crit}(f)$};
    \node[anchor=west] at (6,0) {$F_{R_1}$};
    \node[anchor=west] at (10:6) {$F_{R_1 \e{i\delta}}$};
\node[anchor=west] at (8,0.5) {$G_1$};
\node at (-4.5,4.5) {$\cstar{n}$};
\node[anchor=north] at (-4.5,-4.5) {$|W| = R_1$};
    \end{tikzpicture}
\end{center}
\caption{From proof of Proposition \ref{prop:A2}: the image $G_1$ of the sectorial hypersurface $G$ after applying the isotopy that exchanges $F_{R_1 \e{i \delta}}$ and $F_{R_1}$}
\label{fig:step4}
\end{figure}

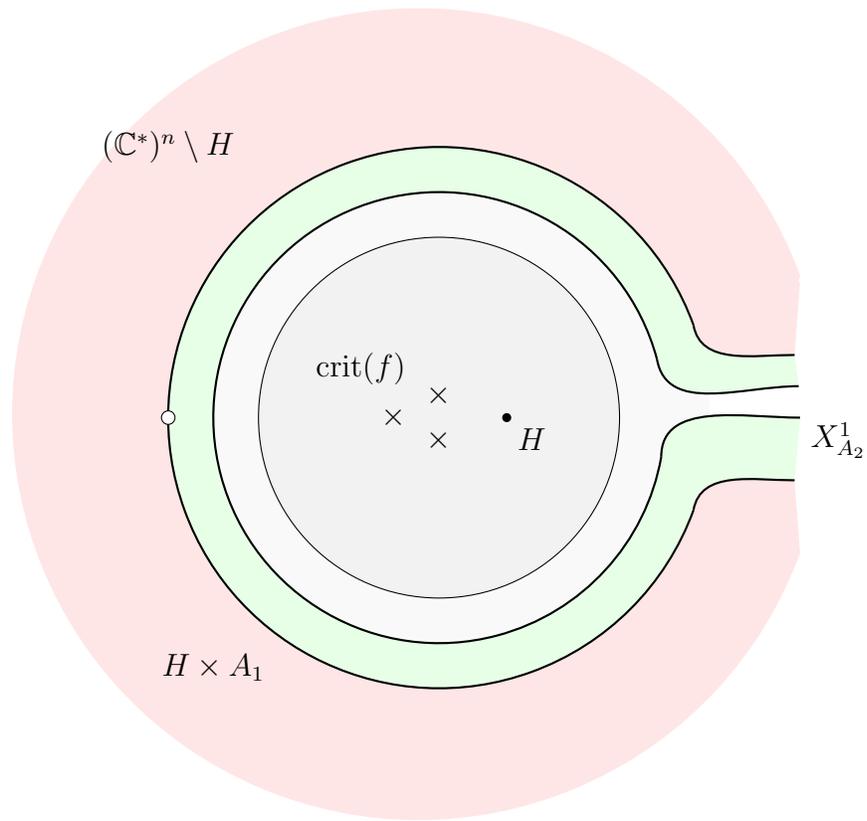
\begin{figure}
\begin{center}
    \begin{tikzpicture}[scale=0.6]
%\draw[help lines] (-8,-8) grid (8,8);
\fill [fill=gray!5] (0,0) circle (6);
\fill [fill=gray!10] (0,0) circle (4);
\draw (0,0) circle (4);

\fill[fill=green!10] (15:5) arc[start angle=15, end angle=350,radius=5] (-10:5) to[out=90, in=180] (0:8) to (-10:8) to[out=180, in=80] (-20:6) arc[start angle=-20, end angle=- 340,radius=6] to[out=-80,in=180] (10:8) to (5:8) to[out=180,in=-80] (15:5);

\fill[fill=red!10] (20:6) arc[start angle=20, end angle=340,radius=6] (-20:6) to[out=90, in=180] (-10:8) to (8
,-3) arc[start angle=-20, end angle=- 340,radius=9] to[out=-80,in=180] (8,3) to (10:8) to[out=180,in=-80] (20:6);

\draw[thick] (15:5) arc[start angle=15, end angle=350,radius=5];
\draw[thick] (15:5) to[out=-80,in=180] (5:8);
\draw[thick] (-10:5) to[out=90, in=180] (0:8);

\draw[thick] (20:6) arc[start angle=20, end angle=340,radius=6];
\draw[thick] (20:6) to[out=-80,in=180] (10:8);
\draw[thick] (-20:6) to[out=80, in=180] (-10:8);
%\draw[thick,dashed] (-180:6) to (-180:9);
%\draw[thick,dashed,->] (20:6.5) arc[start angle=20, end angle=170,radius=6.5];
\draw[fill=white] (-180:6) circle (0.15);
\fill[black] (1.5,0) circle (0.1);
\node[anchor=north west] at (1.5,0) {$H$};
\draw (-0.5,0) node[left] {$\huge{\times}$};
    \draw (0.5,-0.5) node[left] {$\huge{\times}$};
    \draw (0.5,0.5) node[left] {$\huge{\times}$};
    \node[anchor=south east] at (-0.5,0.5) {$\mathrm{crit}(f)$};
\node[anchor=west] at (8,-0.5) {$X_{A_2}^1$};
\node at (-6,6) {$\cstar{n}\setminus H$};
\node[anchor=north] at (-5,-5) {$H \times A_1$};
    \end{tikzpicture}
    \caption{From proof of Proposition \ref{prop:A2}: the sector $X_{A_2}^1$ is in green, the sector $H \times A_1$ is in red.}
    \label{fig:step5}
\end{center}
\end{figure}

\begin{enumerate}
    \item As in the proof of Proposition \ref{prop:sector_gluing} above, after a small deformation we have a standard Liouville neighborhood $U$ of the form:
$$(H \times \set{z \in \CC \;\;: \mathrm{Re}(z) > - \rho}, \lambda_{\cstar{n}}|_F + \lambda_{\CC}^{\mathrm{std}})$$
properly embedded inside $\cstar{n}$. There is then an explicit simple homotopy of Liouville structures on the neighborhood $H \times \set{z \in \CC \;\;: \mathrm{Re}(z) > - \rho}\setminus\set{0}$ that makes both $\set{\mathrm{Re}(z) = 0}$ and $\set{\mathrm{Re}(z) = -\rho/2}$ parallel to the Liouville vector field at infinity (as well as still having $I$-function $I = \mathrm{Im}(z)$). See Figure \ref{fig:vector_field} for a sketch of the resulting Liouville vector field $Z$: it has positive divergence everywhere in $\CC^{\ast}$ and so defines a Liouville form there. This presents $\comp$ as the sector gluing of $(\cstar{n}, f)$ to two sectors $X_{A_2} \cong H \times A_2$ and $X_{A_1} \cong H \times A_1$ along the boundary components $\partial(X_{A_2})$, isotopic to $H \times \RR$ (see Figure \ref{fig:step2}). 
    \item After an analogous argument to \cite[Proposition 6]{periodicity}, applying an isotopy that moves the removed fiber off to infinity, the uniqueness result of \cite[Lemma 2.32]{GPS1}, applies to show that the Liouville sector $X_0 = \cstar{n}\cap\set{\mathrm{Re}(z) \leq 0}$ is Liouville homotopic to the result of applying \cite[Proposition 2.6]{sylvan_lemma} to the Liouville domain $\cstar{n}_{0} = \set{|W| \leq R_1} \cap \set{\psi \leq R_2}$ along the pair of Liouville hypersurfaces $F = |W|\inv(\set{R_1, R_1\e{i\delta}}) \cap \set{\psi \leq R_2-\epsilon}$. Here we identify $F_{R_1\e{i \delta}} \to F_{R_1}$ via $\Phi_\eta$. Equivalently, let $q:\CC \to \CC$ be a double cover so that the two roots of $q(z) = c$ trace out the curves $\gamma$ and $\eta(t) = \e{i t \delta}$ as $c$ travels around the unit circle. Letting $\tilde{f}(x) = q(1-f(x))$, the Liouville sector $\cstar{n}\cap\set{\mathrm{Re}(z) \leq 0}\setminus \set{0}$ is Liouville homotopic to \cite[Proposition 2.6]{sylvan_lemma} now applied to the Liouville domain $\cstar{n}_0$ with hypersurface $F = \tilde{f}\inv(1) \cap \set{\phi \leq R_2-\epsilon}$. 
    \item By the homotopical uniqueness of \cite[Lemma 2.32]{GPS1}, taking the isotopy $\tilde{F}_t = \tilde{f}\inv(\e{2\pi it})$ of Liouville hypersurfaces for $t \in [0,1]$ gives a unique Liouville homotopy of the associated Liouville sectors $X_t$. The two components of the sectorial boundary $\partial X_1$ will be canonically identified with $F_{R_1 \e{i\delta}}\times \RR$ and $F_{R_1} \times \RR$ by retracting the sectorial boundary under the Liouville flow (corresponding to the two black points labelled in Figure \ref{fig:vector_field}). Furthermore, the sectorial hypersurface $G = \set{\mathrm{Re}(z) = -\rho/2} \subset X_0$ will be carried through this homotopy to a sectorial hypersurface $G_t \subset X_t$ and remain parallel to the Liouville vector field at infinity: see Figure \ref{fig:step3}. This sectorial hypersurface presents $X_t$ as the gluing of two sectors: by the uniqueness Lemma, the first is Liouville homotopic to $(\cstar{n},f)$, and the latter, $X_{A_2}^t$, to $H \times A_2$ (see Figure \ref{fig:step5}). 
    \item Lastly, the result of gluing the sector $H \times A_1$ to $X_1$ along $\partial X_1$ yields a Liouville manifold homotopic to $\comp$ (see Figure \ref{fig:step5}). Let us understand the identifications involved in this gluing. The Liouville isotopy $\tilde{F}_t$ gives symplectomorphisms $\Phi_{\gamma}: F_{R_1} \to H$ and $\Phi_{\eta}: F_{R_1 \e{i\delta}} \to H$ which present the boundary $\partial X_1$ as $(H \times \RR) \sqcup (H \times \RR)$. In our original presentation $X_0$, we obtain $\comp$ by gluing $X_{A_1}$ to $X_{A_2}$. Recall that the fiber $F_{R_1}$ was identified with the reference fiber $H$ via parallel transport along a ray. This fiber $F_{R_1}$ is the source of the Liouville flow after the deformation of \cite[Proposition 2.6]{sylvan_lemma} and so each component of the sectorial boundaries of $(\cstar{n}, f)$, $X_{A_2}$ and $X_{A_1}$ is canonically identified with $F_{R_1} \times \RR$. Therefore gluing back $X_{A_1}$ to $X_1$ via the original gluing maps $\mathrm{id}:F_{R_1} \to F_{R_1}$ and $\Phi_{\eta}: F_{R_1} \to F_{R_1 \e{i\delta}}$, and composing with our given identifications with $H$ gives the gluing maps $\mu: H \to H$ and $\Phi_{\eta}:H \to H$ (which is isotopic to the identity). This completes the proof of our claim. \qedhere
\end{enumerate}
\end{proof}

\begin{remark}
We should think of this as the sectorial decomposition of $\comp$ associated to the presentation as $\tilde{f}_\infty=0$. We now use the projection of $\comp$ given by $-1/f: \comp \to \CC$. This has the effect of swapping $0$ and $\infty$: all of the monodromy of $f$ now occurs around $0 \in \CC^{\ast}$. Correspondingly, what one would like to do is use the hypersurface $\mathrm{Re}(1/f)=1/2$ to cut $\comp$ into two sectors. Technical details make this difficult, so we use instead the argument presented in Proposition \ref{prop:A2} to give an equivalent presentation as a gluing of sectors. Future work of the second author and collaborators will remove the need for this argument.
\end{remark}

\begin{definition}\label{defn:alphas}
Given the sectorial decompositions described above, we define:
\begin{enumerate}
    \item The functor $\alpha_{\infty}: \scr{W}( (\CC^{\ast})^n, f) \to \scr{W}(\comp)$ is the pushforward functor from \cite[(3.61)]{GPS1} induced by the sectorial inclusion in Proposition \ref{prop:sector_gluing};
    \item The functor $\alpha_{0}: \scr{W}( (\CC^{\ast})^n, f) \to \scr{W}(\comp)$ is the pushforward functor from \cite[(3.61)]{GPS1} induced by the sectorial  inclusion in Proposition \ref{prop:A2}.
\end{enumerate}
\end{definition}

\begin{corollary}\label{cor:pushouts}
We have pushout diagrams of wrapped Fukaya categories:
\begin{center}
\begin{tikzcd}
\scr{W}(H) \ar{r}{\cup} \ar{d} & \scr{W}((\CC^{\ast})^n, f) \ar{d}{\alpha_\infty}\\
 \scr{W}(H \times \CC^{\ast}(2), z) \ar{r} & \scr{W}_{0}( (\CC^{\ast})^n \setminus H)
\end{tikzcd}
\begin{tikzcd}
\scr{W}(H) \ar{r}{\cup} \ar{d} & \scr{W}((\CC^{\ast})^n, f) \ar{d}{\alpha_0}\\
 \scr{W}(E, H) \ar{r} & \scr{W}_{\infty}( (\CC^{\ast})^n \setminus H)
\end{tikzcd}
\end{center}
where $\scr{W}(H \times \CC^{\ast}(2), z)$ denotes the wrapped Fukaya category where simple counterclockwise Reeb orbits around $0$ are given degree $2$.
\end{corollary}
\begin{proof}
See \cite[Theorem 1.28]{GPS2}: inspection of the arguments show that they extend to $\ZZ$-graded categories (see \cite{GPS3}), provided that the grading data glues appropriately. We observe that since $\eta_0$ is obtained from $\eta_{\infty}$ by twisting by the cohomology class $[\arg(f)]$, whose restriction to the subsector $((\CC^{\ast})^n, f)$ is trivial (as $f$ is trivial as a Cartier divisor there), we may use the same $\ZZ$-graded category $\scr{W}((\CC^{\ast})^n, f)$ in both diagrams.
\end{proof}

% In fact they are already written in the setting of $\ZZ$-graded categories, since they are insensitive to whether we pick grading data or not.

\begin{remark}
Note that these geometric gluing results hold independently of which grading structure we choose; hence we could write both of the above diagrams with either the $\eta_0$-grading or the $\eta_\infty$-grading
\end{remark}

\subsection{Mirror Symmetry}

\begin{proposition} \label{prop:schober_stuff}
We have a commutative diagram:
\begin{center}
\begin{tikzcd}
\scr{W}(H) \ar{r} \ar[bend right=20,swap]{d}{\cup} & \mathrm{Coh}(D) \ar[bend right=20,swap]{d}{i_{D \ast}}\\
 \scr{W}( (\CC^{\ast})^n, f) \ar{r} \ar[bend right=20,swap]{u}{\cap} & \mathrm{Coh}(X) \ar[bend right=20,swap]{u}{i_{D}^{\ast}}
\end{tikzcd}
\end{center}
where the horizontal arrows are the equivalences given by \cite{BenVivek,Kuwagaki,GPS3}. Moreover, the counterclockwise monodromy functor $\mu: \scr{W}(H) \to \scr{W}(H)$ of $f$ around $\infty$ is mirror to tensoring by the anticanonical bundle $ \scr{K}_{X}\inv|_D$ on $\mathrm{Coh}(D)$.
\end{proposition}
\begin{proof}
We know that $\cup$ and $i_{D \ast}$ correspond under the mirror symmetry equivalence of  \cite{BenVivek} (see Figure 6). Therefore their left adjoints $\cap^L$ and $i_{D}^{\ast}$ must also be mirror. Since the sector $((\CC^{\ast})^n,f)$ was defined as the suturing of a stopped Liouville manifold (Definition \ref{def:general_sector}) whose stop is swappable \cite[Definition 1.1]{sylvanOrlov} by Proposition \ref{prop:swappable}, then by \cite[\S 4.2]{sylvanOrlov} the functor $\cup$ is spherical, and the image of the left adjoint $\cap^L$ is contained in $\scr{W}(H)$. In the language of perverse schobers, the dual cotwist for a spherical functor $S$ is given by the cone of the counit map $S^{L}S \to \mathrm{id}$ (where $S^{L}$ denotes the left adjoint). In \cite[Theorem 1.3]{sylvanOrlov} it is proved that $\mu$, the geometric counterclockwise monodromy of $f$ around $\infty$, acts on $\scr{W}(H)$ by the dual cotwist of $\cup$, shifted by $[2]$. From \cite[p.12]{Nadler}, we see that the dual cotwist for $i_{D\ast}$ is given by $\otimes \scr{K}_{X}\inv|_{D}[-2]$. Cancelling the $2$-shifts, we see that the functor $\mu$ corresponds to $\otimes \scr{K}_{X}\inv|_D$.
\end{proof}

\begin{remark}
    The cup functor $\cup:\scr{W}(H)\to \scr{W}((\CC^\ast)^n,f)$ is (restriction to compact objects of) left adjoint to the microlocalization functor studied in \cite[\S 7.2]{BenVivek}. It is shown there that when $X$ and $D$ are smooth, this microlocalization is mirror to the pullback functor $i^*:\QCoh(X)\to \QCoh(D).$ However, this equivalence uses an equivalence between $\QCoh$ and $\IndCoh$ for smooth schemes which intertwines the quasicoherent pullback $i^*$ with the Ind-coherent pullback $i^!$. (See \cite[Chapter 6, \S 3.2]{GRI}.) Since $i$ is proper, the left adjoint to $i^!$ is $i_*,$ as claimed above. The same equivalence is used implicitly in Corollary \ref{cor:compatibility} below.
\end{remark}

%Note that this convention for the monodromy is the *reverse* of the one from my periodicity paper: the natural transformation goes in the opposite direction.

\begin{proposition}\label{prop:twisted_schober}
There is a mirror symmetry equivalence between $(E,H)$ and $K_X|_D$ such that the following diagram commutes:
\begin{center}
\begin{tikzcd}
\scr{W}(H) \ar{r} \ar{d}{\cup} & \mathrm{Coh}(D) \ar{d}{\iota_{D \ast}}\\
 \scr{W}( E,z) \ar{r} & \mathrm{Coh}(K_X|_D)
\end{tikzcd}
\end{center}
where the top horizontal arrow is the equivalence of \cite{BenVivek}.
\end{proposition}
\begin{proof} 
The total space of $K_{X}|_{D}$ is the relative Spec of the sheaf of algebras $\mathrm{Sym}(\scr{K}_{X}\inv|_D)$ on $D$. The category of coherent sheaves on the total space is therefore given by the category of quasicoherent sheaves $\scr{F}$ on $D$ with a coherent $\mathrm{Sym}(\scr{K}_{X}\inv|_D)$-module structure. By the universal property of $\mathrm{Sym}$, this data is the same as a morphism  $\scr{F} \otimes \scr{K}_{X}\inv|_D \to \scr{F}$ describing how the generators of $\mathrm{Sym}(\scr{K}_{X}\inv|_D)$ act on sections of $\scr{F}$. We can thus describe $\mathrm{Coh}(K_X|_D)$ as the category of pairs $(\scr{F}, \varphi: \scr{F} \otimes \scr{K}_{X}\inv|_D \to \scr{F})$ of a quasicoherent sheaf $\scr{F}$ on $D$ along with a morphism $\scr{F} \otimes \scr{K}_{X}\inv|_D \to \scr{F}$ giving $\scr{F}$ the structure of a finitely generated module over the $\cdot \otimes \scr{K}_{X}\inv|_D$ monad. Compare \cite[Proposition 5.8]{MilnorFibers}. The pushforward functor $\Coh(D) \to \Coh(K_{X}|_{D})$ sends a sheaf $\scr{F}$ to the trivial action $\scr{F} \otimes \scr{K}_{X}\inv|_D \to \scr{F}$.

We have a similar description of $\scr{W}(E,H)$: the sector decomposition in Lemma \ref{prop:A2} presents $\scr{W}(E,H)$ as a homotopy colimit of wrapped Fukaya categories (note that even though we used $1-f$ rather than $f$, the parallel transport $\mu$ is naturally isomorphic to the geometric counterclockwise monodromy functor). This can be realized by an explicit homotopy colimit  construction as in \cite[A.4]{GPS2}, where $\scr{W}(E,H)$ is given by a category of pairs $(L, s: \mu L \to L)$ where $L$ is a module over $\scr{W}(H)$ and $s: \mu L \to L$ is a morphism giving $L$ a finitely-generated module structure over the $\mu$ monad. Moreover, the cup functor $\scr{W}(H) \to \scr{W}(E,H)$ sends a Lagrangian $L$ to the zero morphism $\mu L \to L$. 

Combining these presentations of the categories $\scr{W}(E,H)$ and $\Coh(K_{X}|_{D})$ with Proposition \ref{prop:schober_stuff}, the result follows.
\end{proof}

\begin{proposition}\label{prop:derived_pushout}
We have pushout diagrams of categories:
\begin{center}
\begin{tikzcd}
\mathrm{Coh}(D) \ar{r}{i_{D\ast}} \ar{d}{\iota_{\tilde{D}\ast}} & \mathrm{Coh}(X) \ar{d}\\
 \mathrm{Coh}(D \times \mathbb{A}^1[-1]) \ar{r} & \mathrm{Coh}(\tilde{Z})
\end{tikzcd}
\begin{tikzcd}
\mathrm{Coh}(D) \ar{r}{i_{D\ast}} \ar{d}{\iota_{D\ast}} & \mathrm{Coh}(X) \ar{d}\\
 \mathrm{Coh}(K_{X}|_{D}) \ar{r} & \mathrm{Coh}(Z)
\end{tikzcd}
\end{center}
where $\mathrm{Coh}(\tilde{Z})$ denotes the category of coherent sheaves on the derived scheme $\tilde{Z}$.
\end{proposition}
\begin{proof}
Follows from proper descent for ind-coherent sheaves on derived schemes (from \cite[Proposition 7.2.2]{GRI}) and passing to compact objects.
\end{proof}

\begin{theorem}\label{thm:alphas}
After choosing appropriate mirror equivalences induced by the pushout diagrams above:
\begin{enumerate}
    \item There is an equivalence of categories between $\scr{W}_{0}(\comp)$ and $\mathrm{Coh}(\tilde{Z})$ so that the functor $\alpha_\infty$ is mirror to $i_{\tilde{X}\ast}$;
    \item There is an equivalence of categories between $\scr{W}_{\infty}(\comp)$ and $\mathrm{Coh}(Z)$ so that the functor $\alpha_0$ is mirror to $i_{X\ast}$.
\end{enumerate}
\end{theorem}
\begin{proof}
By Propositions \ref{prop:schober_stuff} and \ref{prop:twisted_schober}, the mirror equivalences of \cite{BenVivek} intertwine the cup functors of Corollary \ref{cor:pushouts} with the pushforward functors on coherent sheaves in Proposition \ref{prop:derived_pushout}. Therefore, there exist equivalences between the pushout categories $\scr{W}_{0}(\comp)$ and $\mathrm{Coh}(\tilde{Z})$, and between $\scr{W}_{\infty}(\comp)$ and $\mathrm{Coh}(Z)$ making the diagrams below commute.
\end{proof}

\[\begin{tikzcd}
	& {\mathrm{Coh}(D)} && {\mathrm{Coh}(D \times \mathbb{A}^1[-1])} \\
	{\mathcal{W}(H)} && {\mathcal{W}(H \times \CC^{\ast}(2),z)} \\
	& {\mathrm{Coh}(X)} && {\mathrm{Coh}(\tilde{Z})} \\
	{\mathcal{W}((\mathbb{C}^{\ast})^n, f)} && {\mathcal{W}_0((\mathbb{C}^{\ast})^n \setminus H)}
	\arrow["\cup"{description, pos=0.3}, from=2-1, to=2-3]
	\arrow[""{name=0, anchor=center, inner sep=0}, "\cup"', from=2-1, to=4-1]
	\arrow[""{name=1, anchor=center, inner sep=0}, "{\alpha_{\infty}}"', from=4-1, to=4-3]
	\arrow["\sim", from=2-1, to=1-2]
	\arrow["\sim", from=2-3, to=1-4]
	\arrow["\sim", from=4-3, to=3-4]
	\arrow["\sim", from=4-1, to=3-2]
	\arrow["{\iota_{\tilde{D}\ast}}"{description}, from=1-2, to=1-4]
	\arrow["{i_{D\ast}}"{description, pos=0.3}, from=1-2, to=3-2]
	\arrow["{j_{\tilde{D}\ast}}"{description}, from=1-4, to=3-4]
	\arrow["{i_{\tilde{X}\ast}}"{description, pos=0.6}, from=3-2, to=3-4]
	\arrow[from=2-3, to=4-3]
\end{tikzcd}\]

\[\begin{tikzcd}
	& {\mathrm{Coh}(D)} && {\mathrm{Coh}(K_{X}|_D)} \\
	{\mathcal{W}(H)} && {\mathcal{W}(E,H)} \\
	& {\mathrm{Coh}(X)} && {\mathrm{Coh}(Z)} \\
	{\mathcal{W}((\mathbb{C}^{\ast})^n, f)} && {\mathcal{W}_\infty((\mathbb{C}^{\ast})^n \setminus H)}
	\arrow["\cup"', from=2-1, to=4-1]
	\arrow["{\alpha_{0}}"', from=4-1, to=4-3]
	\arrow["\sim", from=2-1, to=1-2]
	\arrow["\sim", from=2-3, to=1-4]
	\arrow["\sim", from=4-3, to=3-4]
	\arrow["\sim", from=4-1, to=3-2]
	\arrow["{\iota_{D\ast}}"{description}, from=1-2, to=1-4]
	\arrow["{i_{D\ast}}"{description, pos=0.3}, from=1-2, to=3-2]
	\arrow["{j_{D\ast}}"{description}, from=1-4, to=3-4]
	\arrow["{i_{X\ast}}"{description, pos=0.6}, from=3-2, to=3-4]
	\arrow[from=2-3, to=4-3]
	\arrow["\cup"{description, pos=0.3}, from=2-1, to=2-3]
\end{tikzcd}\]

\begin{remark}
Theorem \ref{thm:alphas} is stated in this fashion since the argument given does not show that the mirror symmetry equivalence between $\scr{W}_{\infty}(\comp)$ and $\Coh(Z)$ induced by these pushout diagrams is necessarily the same as those coming from \cite{BenVivek} (and likewise for $\scr{W}_{0}(\comp)$ and $\Coh(\tilde{Z})$ ), though we shall show that this is indeed the case in \S \ref{sec:compatibility2} below. 
\end{remark}

\begin{remark}
One cannot construct a functor $\scr{W}_{0}(\comp) \to \scr{W}_{\infty}(\comp)$ by sending Lagrangians to themselves, since Lagrangians may not necessarily have a graded lift for both gradings. Even if we restrict to the subcategory of Lagrangians graded for both, a functor that is the identity on objects and morphisms will not be $\ZZ$-graded, but only a functor on the $\ZZ/2$-graded categories. Similarly, it may seem that there is an obvious functor $\mathrm{Coh}(Z) \to \mathrm{Coh}(\tilde{Z})$, but the correct $\ZZ$-graded equivalence to goes via a Kn\"orrer periodicity equivalence. The corresponding equivalence on the A-side can be thought of informally as applying the transformation $z \mapsto 1/z$ in the base.
\end{remark}

\subsection{Compatibility of Equivalences} \label{sec:compatibility2}

In this section we verify that the sectorial decomposition of $(\CC^{\ast})^n \setminus H$ coming from applying Sylvan's construction to the Liouville structure on $(\CC^{\ast})^n \setminus H$ induced by $f$ respects the combinatorial decomposition of the skeleton coming from the semitropicalization of $f$, under the isotopy relating the two (a priori different) Liouville structures. We will do this by examining the image of the different skeleta of the same hypersurface $(\CC^{\ast})^n \setminus H = \{\tilde{f}_0  = 0\}=\{\tilde{f}_{\infty}=0 \}$ under $f:(\CC^{\ast})^n \to \CC$. Though they have the same zero-locus, the tropicalizations of $\tilde{f}_0$ and $\tilde{f}_{\infty}$ differ. If $A$ denotes the set of multi-indices of monomials in the Laurent polynomial $f$, then we have
$$A_{0} = \{ (0, \alpha) \; : \; \alpha \in A\} \cup \{(-1, 0, \dots, 0)\}$$ while $$A_{\infty} = \{ (1, \alpha) \; : \; \alpha \in A\} \cup \{(0, 0, \dots, 0)\}$$ since $\tilde{f}_{\infty} = z_{0}f + 1$ and $\tilde{f}_0 = f + z_{0}^{-1}$. If $\phi:A \to \RR$ is a convex piecewise-linear function inducing a polyhedral decomposition of $P$, we also obtain a natural convex piecewise-linear function $\tilde{\phi}_{\infty}: A_{\infty} \to \RR$ by setting $\tilde{\phi}_{\infty}(0) = 0$ and $\tilde{\phi}(1, \alpha) = \phi(\alpha)$; and $\tilde{\phi}_0: A_0 \to \RR$ by setting $\tilde{\phi}_0(0, \alpha) = \phi(\alpha)$ and $\tilde{\phi}(-1,0) = 0$. Recall that, using this convex function $\phi$, one defines the tropicalization $L_{\phi}:\RR^n \to \RR$ of $f$ by
\begin{equation*}
	L_{\phi}(u) = \max_{\beta \in A}\{ \langle \beta, u \rangle - \phi(\beta) \}
\end{equation*}
The tropicalization of $f^{-1}(0)$ is given by $\Pi \subset \RR^n$, the locus of discontinuity of $L_{\phi}$. The complement $\Pi^{c} \subset \RR^n$ is given by the locus where one term in $L_{\phi}$ strictly dominates the others; the closure of the region where the constant term $0 \in A$ dominates is denoted by $Q$. Moreover, if $A$ is the set of rays generating the cones in a fan $\Sigma$ of the toric variety $X$ from \S \ref{subsec:terms}, then we likewise have fans $\Sigma_0, \Sigma_{\infty}$ generated by rays in $A_{0}, A_{\infty}$ respectively, determined by the polyhedral decomposition associated to $\phi$. Write $\Sigma_{0}^{0}$ for the subset of cones in $\Sigma_0$ generated only by rays $(0,\alpha)$ for $\alpha \in A$, and write $\Sigma_{0}^1$ for those cones in $\Sigma_0$ whose generators include $(-1,0,\dots, 0)$; likewise for $\Sigma_{\infty}^0, \Sigma_{\infty}^1$.  

Therefore the tropicalizations $L_{\infty}, L_0: \RR \times \RR^n \to \RR$ of $f_{\infty}, f_0$ are, respectively,
\begin{equation*}
	L_{\infty}(u_0, u) = \max\{ 0, u_0+L_{\phi}(u) \}
\end{equation*}
and
\begin{equation*}
	L_{0}(u_0, u) = \max\{ -u_0, L_{\phi}(u) \}
\end{equation*}
Let $\Pi_{\infty}, \Pi_0 \subset \RR^{n+1}$ be the the tropicalizations of $f_{\infty}^{-1}(0), f_{0}^{-1}(0)$ respectively. 

The facets of $\Pi_{\infty}$ are given by the locus where two terms of $L_{\infty}$ become equal. These occur in two kinds: firstly, those where $u_0 = - L_{\phi}(u)$ (those facets corresponding to the boundary of the region in which the constant term dominates); secondly, regions where the terms associated to $(1, \alpha_1), (1, \alpha_2) \in A_{\infty}$ codominate: this occurs when $$u_0 + \langle \alpha_1, u \rangle - \phi(\alpha_1) = u_0 + \langle \alpha_2, u \rangle - \phi(\alpha_2)$$
Since this equality is $u_0$-independent, this facet is simply given by
\begin{equation*}
	\{u_0 \geq - L_{\phi}(u) \} \cap \{ \langle \alpha_1, u \rangle - \phi(\alpha_1) =  \langle \alpha_2, u \rangle - \phi(\alpha_2) \}
\end{equation*}
that is, a half-line $[c, \infty)$ times a facet of $\Pi$. 

For $\Pi_0$, the region corresponding to where the constant term $\alpha = 0 \in A$ dominates in $L_{0}$ will be $(-\infty, 0] \times Q$; we also have facets where $-u_0  = L_{\phi}(u)$ but these facets are not in the boundary of $Q$. 

The key result from \cite{BenVivek, Zhou} is that, for appropriate choice of Liouville structure on $(\CC^{\ast})^n$, there is a cover $U_{\sigma}$ of the Lagrangian skeleton of $H$ corresponding to cones $\sigma \in \Sigma$ in the fan of $X$. The perfectly centered hypothesis from \cite{BenVivek} means that every non-zero cone $\sigma \in \Sigma$ intersects its dual face $F \subset Q$ in the relative interior. For each top-dimensional cone $\sigma \in \Sigma$, let $V_{\sigma}$ be a neighborhood of $\sigma \subset \RR^n$, and let $U_{\sigma} = \mathrm{Log}^{-1}(V_{\sigma}) \cap H$. This is then an open subset of $H$, the union of which contains its entire skeleton. Moreover, there is an equivalence between $\mathrm{\mu Sh}(U_{\sigma})$ and $\mathrm{Coh}(\overline{O}(\sigma))$, where $\overline{O}(\sigma)$ the closure of the torus orbit corresponding to $\sigma$, and this is functorial under inclusions of cones $\sigma \subset \tau$ in the sense that restriction on coherent sheaves is mirror to the microlocalization functor on $\mathrm{\mu Sh}$.
 
Note that on $(\CC^{\ast})^n \setminus H \subset (\CC^{\ast})^{n+1}$ we have $z_0 = 1/f$, so under the log map $(\CC^{\ast})^{n+1} \to \RR^{n+1}$ we have $u_0 = - \log |f|$. Let $\Delta_0 \subset \RR^{n}$ be the subset of the complement of the tropical amoeba $\Pi$ for $H$ where the constant term dominates. Let $U$ be the inverse image under the log map of a slight enlargement of $\Delta_0$. 

% \begin{proposition}
% For $t$ sufficiently small, the subsets $U_{\sigma}$ of the open cover of the Liouville skeleton of $(\CC^{\ast})^n \setminus H$ corresponding to cones $\sigma \in \tilde{\Sigma} \setminus \{(1,0,\dots, 0)\}$ are all contained in the set $|f_{t} - 1|>\delta$ for some $\delta>0$. 
% \end{proposition}

\begin{lemma}
	\cite[Lemma 5.2]{Speculations} If $T>0$ is sufficiently small, the critical points of $f$ all lie inside $U$ and the critical values converge to $1$. 
\end{lemma} 

\begin{corollary}\label{cor:compatibility}
    There is an equivalence of categories between $\mathrm{colim}_{\sigma \in \Sigma_0 \cap \Sigma_1} \mu Sh(U_{\sigma})$ and $\scr{W}(H)$ such that the following diagram commutes:
    \[\begin{tikzcd}
	\scr{W}(H) & \scr{W}(\cstar{n},f) \\
	\scr{W}(E,H) & \scr{W}_{\infty}(\cstar{n}\setminus H) && \colim_{\sigma \in \Sigma_{\infty}^0 \cap \Sigma_{\infty}^1}\mu\mathrm{Sh}(U_{\sigma}) & \colim_{\sigma \in \Sigma_{\infty}^0}\mu \mathrm{Sh}(U_{\sigma}) \\
	\Coh(D) & \Coh(X) && \colim_{\sigma \in \Sigma_{\infty}^1} \mu \mathrm{Sh}(U_{\sigma}) & \colim_{\sigma \in \Sigma}\mu \mathrm{Sh}(U_{\sigma} )\\
	\Coh(K_{X}|_{D}) & \Coh(Z)
	\arrow[from=1-1, to=1-2]
	\arrow[from=1-1, to=2-1]
	\arrow[from=1-2, to=2-2]
	\arrow[from=2-1, to=2-2]
	\arrow[from=2-4, to=1-1]
	\arrow[from=2-4, to=2-5]
	\arrow[from=2-4, to=3-1]
	\arrow[from=2-4, to=3-4]
	\arrow[from=2-5, to=1-2]
	\arrow[from=2-5, to=3-2]
	\arrow[from=2-5, to=3-5]
	\arrow[bend left=70, dashed, from=3-1, to=1-1]
	\arrow[from=3-1, to=3-2]
	\arrow[from=3-1, to=4-1]
	\arrow[bend left=70, dashed, from=3-2, to=1-2]
	\arrow[from=3-2, to=4-2]
	\arrow[from=3-4, to=2-1]
	\arrow[from=3-4, to=3-5]
	\arrow[from=3-4, to=4-1]
	\arrow[from=3-5, to=2-2]
	\arrow[from=3-5, to=4-2]
	\arrow[bend left=70, dashed, from=4-1, to=2-1]
	\arrow[from=4-1, to=4-2]
	\arrow[bend left=70, dashed, from=4-2, to=2-2]
\end{tikzcd}\]
as well as an analogous diagram for $\scr{W}_{0}((\CC^{\ast})^n\setminus H)$, where the arrows in the left hand squares are those from Theorem \ref{thm:alphas}, the south-west-pointing arrows are the equivalences from \cite{BenVivek}, and the north-west-pointing arrows are those from \cite{GPS3}.  
\end{corollary}
\begin{proof}
After applying a semitropicalization to either $\tilde{f}_{0}$ or $\tilde{f}_{\infty}$ we get a (simple) homotopic Liouville structure on $(\CC^{\ast})^n \setminus H$ whose skeleton is close to either $\Pi_{0} \subset \RR^{n+1}$ or $\Pi_{\infty} \subset \RR^{n+1}$ under the $\mathrm{Log}$ map on $(\CC^{\ast})^{n+1}$. There exists some $\epsilon>0$ so that all of the critical points of the quadratic Stein function $\psi$ from \cite{Zhou} restricted to $(\CC^{\ast})^n \setminus H$ corresponding to cones $\sigma \in \Sigma_{0}^1$ have $|u_0| > \epsilon$. But recall that $u_0 = \mathrm{Log}|z_0| = \mathrm{Log}|f_{s}|^{-1}$, or in other words  this means that $|f_{s}(z) - 1| > \delta$ for some $\delta>0$. Since all of the critical points of $f_{s}$ live near $1$ for $T>0$ sufficiently small, the sectorial gluing decompositions in Propositions \ref{prop:sector_gluing} and \ref{prop:A2} also respect this decomposition of the skeleton. The commutative diagram above then follows from the combinatorics of $\Sigma_{\infty}$ and general properties of colimits.
\end{proof}

\section{The Restriction Functor}\label{sec:rho}

\subsection{Kn\"orrer Periodicity}

To relate the two mirror symmetry equivalences in Theorem \ref{thm:alphas} above, we will first need to understand the equivalence between $\mathrm{Coh}(Z)$ and $\mathrm{Coh}(\tilde{Z})$ induced by the Kn\"orrer periodicity equivalences from Theorem \ref{prop:orlov}. Firstly, we will need to understand the effect on coherent sheaves coming from $D$ (Lemma \ref{lemma:orlov}) and secondly on coherent sheaves coming from $X$ (Lemma \ref{lem:locally_free}).

\begin{lemma}\label{lemma:orlov}
Kn\"orrer periodicity gives a commutative diagram of dg functors:
\begin{center}
\begin{tikzcd}
{} & \mathrm{Coh}_{\mathbb{G}_m}(K_{X} \times \mathbb{A}^1, s(x)yz) & {} \\
\mathrm{Coh}(Z) \ar{ur} & {}  & \ar{ul} \mathrm{Coh}(\tilde{Z}) \\
{} & \mathrm{Coh}(D) \ar{ur} \ar{ul}& {}
\end{tikzcd}
\end{center}
\end{lemma}

\begin{proof}
The lower two maps are given by the pull-push functors $j_{D\ast}p^{\ast}_{D}$ (for $p_{D}: K_{X}|_D \to D$, $j_D:K_{X}|_D \to Z$) and $j_{\tilde{D} \ast}p_{\tilde{D}}^{\ast}$ (for $p_{\tilde{D}}: D \times \mathbb{A}^1[-1] \to D$, $j_{\tilde{D}}:D \times \mathbb{A}^1[-1] \to \tilde{Z}$). The upper left functor too is given by $j_{Z \ast}p^{\ast}_{Z}$ (for $p_Z$ the projection $Z \times \mathbb{A}^1 \to K_X$; and $j_{Z}$ the inclusion $Z \times \mathbb{A}^1 \to K_X \times \mathbb{A}^1$). The upper right functor is the chain of equivalences from Propositions \ref{prop:knorrer2} and \ref{prop:compatibility}:
\begin{equation*}
    \mathrm{Coh}(\tilde{Z}) \to \mathrm{Coh}_{\mathbb{G}_m}(\tilde{Z}, 0) \to \mathrm{Coh}_{\mathbb{G}_m}(K_{X} \times \mathbb{A}^1, s(x)yz)
\end{equation*}
where the second is also given by a pull-push functor $j_{\tilde{S}\ast}p^{\ast}_{\tilde{S}}$. One can check from our explicit description of the first equivalence in the proof of Proposition \ref{prop:knorrer2} that a coherent sheaf on $X$ supported on $D$ that is a free $z$-module gives a matrix factorization that is also supported on $D$ and free in the $z$-variable. Finally, one can simply check on affines that on both sides of the diagram $$j_{Z\ast}p_{Z}^{\ast} j_{D\ast}p_{D}^{\ast} = j_{Z\ast}(p_{D} \circ p_{Z}|_{K_{X}|_D})^{\ast} = j_{\hat{D}\ast}p_{\hat{D}}^{\ast} = j_{\tilde{S}\ast}p_{\tilde{S}}^{\ast} j_{\tilde{D}\ast}p_{\tilde{D}}^{\ast}$$ for $p_{\hat{D}}: K_{X}|_D \times \mathbb{A}^1 \to D$ and $j_{\hat{D}}: K_{X}|_D \times \mathbb{A}^1 \to K_X \times \mathbb{A}^1$ which means that the diagram commutes.
\end{proof}

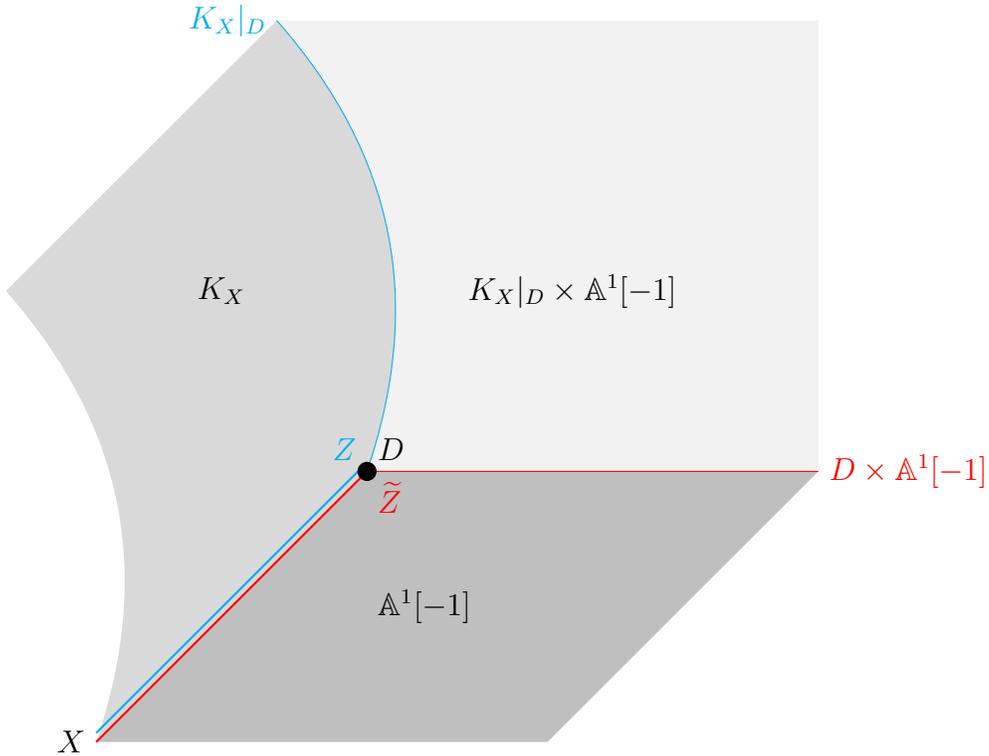
\begin{figure}[H]
\begin{center}
    \begin{tikzpicture}[scale=1.2]
     \fill [fill=gray!30] (0,0) to[bend right] (-1,5) to (-4,2) to[bend left] (-3,-3) to (0,0);
     \fill [fill=gray!10] (0,0) to[bend right] (-1,5) to (5,5) to (5,0) to (0,0);
     \fill [fill=gray!50] (0,0) to (5,0) to (2,-3) to (-3,-3) to (0,0);
 \draw[red] (0,0) -- (5,0);
 \draw (0,-1.5) node[right] {$\mathbb{A}^1[-1]$};
 \draw (-3,-3) node[left] {$X$};
 \draw[cyan] (0,0) node[above left] {$Z$};
  \draw[red] (0,0) node[below right] {$\tilde{Z}$};
 \draw[cyan] (-1,5) node[left] {$K_{X}|_D$};
 \draw (1,2) node[right] {$K_{X}|_{D} \times \mathbb{A}^1[-1]$};
 \draw[red] (5,0) node [right] {$D \times \mathbb{A}^1[-1]$};
 \draw[cyan] (0,0) to[bend right] (-1,5);
  \draw (-2,2) node[right] {$K_X$};
 \draw[thick, red] (0,0) to (-3,-3);
  \draw[thick, cyan] (0,0.1) to (-3,-2.9);
 \filldraw[black] (0,0) circle (0.1) node[above right] {$D$};
     %\draw[help lines] (-5,-5) grid (5,5);
    \end{tikzpicture}
\end{center}
\caption{Schematic picture of $K_X \times \mathbb{A}^1[-1]$ in the case $X = \mathbb{A}^1$.}
\label{fig:schematic}
\end{figure}

\begin{lemma}\label{lem:locally_free}
Kn\"orrer periodicity gives a commutative diagram of dg functors:
\begin{center}
\begin{tikzcd}
{} & \mathrm{Coh}_{\mathbb{G}_m}(K_{X} \times \mathbb{A}^1, s(x)yz) & {} \\
\mathrm{Coh}(Z) \ar{ur} & {}  & \ar{ul} \mathrm{Coh}(\tilde{Z}) \\
{} & \mathrm{Coh}(X) \ar{ur}{i_{\tilde{X}\ast}} \ar{ul}{\pi_{X}^{\ast}[1]}& {}
\end{tikzcd}
\end{center}
\end{lemma}
\begin{proof}
Suppose $\scr{E} = \pi^{\ast}_{X}\scr{F}$ is the pullback of a locally free sheaf $\scr{F}$ on $X$, and let $p_{Z}$ denote the projection $Z \times \mathbb{A}^1 \to Z$. The upper left equivalence $\mathrm{Coh}(Z) \to \mathrm{Coh}_{\mathbb{G}_m}(K_X \times \mathbb{C}, s(x)y z)$ sends $\scr{E}$ to the pushforward of $p^{\ast}_{Z}\scr{E}$ under the inclusion $j_{Z}:Z \times \mathbb{A}^1 \to K_X \times \mathbb{A}^1$. Note that we have an exact sequence of coherent sheaves on $\hat{Z}$:
\begin{equation*}
   0 \to \scr{O}_{Z \times \mathbb{A}^1} \to \scr{O}_{\hat{Z}} \to  \scr{O}_{K_{X} \times \set{0}} \to 0
\end{equation*}
given by restriction of functions $\scr{O}_{\hat{Z}} \to \scr{O}_{K_{X} \times \set{0}}$ (there is no normal-bundle correction since $K_{X} \times \set{0} = \set{z=0}$ is the zero-locus of a section of the trivial line bundle). For the projection $p_{\hat{Z}}: \hat{Z} \to X$, the pullback $p^{\ast}_{\hat{Z}}\scr{F}$ is a locally free sheaf. Tensoring this with the above exact sequence gives an exact triangle:
\begin{center}
\begin{tikzcd}
j_{Z\ast}p^{\ast}_{Z} \scr{E} \ar{rr} & {} & p^{\ast}_{\hat{Z}}\scr{F} \ar{dl}\\
{} &  \scr{O}_{K_{X} \times \set{0}} \otimes p^{\ast}_{\hat{Z}}\scr{F}  \ar{ul}{+1} & {}
\end{tikzcd}
\end{center}
Observe that passing to the quotient category $\mathrm{Sing}(\hat{Z})$ sends all locally free sheaves to zero. Thus $j_{Z\ast}p^{\ast}_{Z}\scr{E}[1]$ is equivalent in $\mathrm{Sing}(\hat{Z})$ to the sheaf $ \scr{O}_{K_{X} \times \set{0}} \otimes p^{\ast}_{\hat{Z}}\scr{F}$. But this sheaf is exactly the image in $\mathrm{Coh}(K_X \times \mathbb{C})$ of $i_{\tilde{X}\ast}\scr{E}|_X$under the functor $j_{\tilde{Z}\ast}p_{\tilde{Z}}^{\ast}:\mathrm{Coh}(\tilde{Z}) \to \mathrm{Coh}_{\mathbb{G}_m}(K_X \times \mathbb{C}, s(x)yz)$. Since $X$ is smooth, the diagram also commutes for all coherent sheaves.

Next we need to check that the diagram commutes on the level of morphisms. Without loss of generality, suppose $h: \scr{F} \to \scr{G}$ is a morphism of locally free sheaves on $X$. Following the arrows up the left hand side, we get the morphism of matrix factorizations given by $\tilde{h}: \scr{O}_{Z \times \AA^1} \otimes p^{\ast}_{\hat{Z}}\scr{F} [1] \to \scr{O}_{Z \times \mathbb{A}^1} \otimes p^{\ast}_{\hat{Z}}\scr{G}[1]$ while on the right hand side we have the morphism $\tilde{h}: \scr{O}_{K_{X} \times \set{0}} \otimes p^{\ast}_{\hat{Z}}\scr{F} \to \scr{O}_{K_{X} \times \set{0}} \otimes p_{\hat{Z}}^{\ast}\scr{G}$, where $\tilde{h}(x,y,z) = h(x)$. These are identified under our choice of isomorphism $\scr{O}_{K_X \times \set{0}} \cong \scr{O}_{Z \times \AA^1}[1]$. 
\end{proof}

\begin{definition}
The functor $q: \bar{\mathrm{Coh}}(Z) \to \bar{\mathrm{Sing}}(Z)$ is the dg quotient functor, sending all perfect complexes to zero.
\end{definition}

\begin{definition}
The localization functor $\mathrm{loc}$ maps $\bar{\mathrm{Coh}}(D \times \CC^{\times}[-1]) \to \bar{\mathrm{Coh}}(D)$ by localization at the degree-$2$ natural transformation given by the action of the generator $z$ of $\CC[z^{\pm 1}] \cong \mathrm{End}(\scr{O}_{\CC^{\times}[-1]})$ (see Lemma \ref{lem:koszul}).
\end{definition}

Compare this definition to Proposition \ref{prop:localization}. Note that both of these functors (as presented here) are $\ZZ/2$-graded only.

\begin{lemma}
The category $\mathrm{Coh}(Z)$ is generated by
\begin{itemize}
    \item the pullbacks $\pi_{X}^{\ast} \scr{E}$ for $\scr{E}$ locally free sheaves on $X$;
    \item the pushforwards $j_{D\ast}p_{D}^{\ast} \scr{F}$ for $\scr{F}$ coherent sheaves on $D$. 
\end{itemize} 
\end{lemma}
\begin{proof}
See \cite[Proposition 2.10]{Orlov3}: though the proposition there is stated in terms of projective bundles, the proof is identical. 
\end{proof}

\begin{proposition}\label{prop:sheafstuff}
Under the Kn\"orrer equivalences in Theorem \ref{prop:orlov}, we have a commutative diagram of $\cz$-linear dg categories:
\begin{center}
\begin{tikzcd}
	{\bar{\mathrm{Coh}}(\tilde{Z})} \ar{rr}{j_{\tilde{D}}^{\ast}} \ar{d}{\sim} && {\bar{\mathrm{Coh}}(D \times \mathbb{C}^{\times}[-1])} \ar{rr}{\mathrm{loc}} && {\bar{\mathrm{Coh}}(D)} \ar{dd}{\sim} \\
	{\bar{\mathrm{Coh}}_{\mathbb{G}_m}(K_X \times \mathbb{A}^1, s(x)yz)} \\
	{\bar{\mathrm{Coh}}(Z)} \ar{u}{\sim} \ar{rrrr}{q} && {} && {\bar{\mathrm{Sing}}(Z)} \\
\end{tikzcd}
\end{center}
\end{proposition}
\begin{proof}
Suppose $\scr{E}$ is a locally free sheaf on $Z$ pulled back from a locally free sheaf on $X$. Then under the sequence of equivalences $\Coh(Z) \to \Coh(\tilde{Z})$ in Lemma \ref{lem:locally_free}, $\scr{E}$ is sent to $i_{\tilde{X}\ast}\scr{E}|_X[-1]$, which is evidently sent to zero under the pullback map to $D \times \CC^{\ast}$. Hence going either way around the diagram sends these locally free sheaves to zero. 

By Lemma \ref{lemma:small_orlov} the $\cz$-linear equivalence $\bar{\mathrm{Coh}}(D) \to \bar{\mathrm{Sing}}(Z)$ factors through the quotient map $q$; so $q$ is an equivalence on those sheaves coming from $D$ via $j_{D\ast}p^{\ast}_{D}$ (for $p_{D}: K_{X}|_D \to D$, $j_D:K_{X}|_D \to Z$ the projection and inclusion respectively). Therefore to complete the proof that the diagram commutes, it suffices to do so for $\scr{F}$ a coherent sheaf on $Z$ supported on $K_{X}|_D$ coming from $D$ via $\scr{F}= j_{D\ast}p_{D}^{\ast} \scr{F}_{D}$. Now by Lemma \ref{lemma:orlov}, this sheaf is sent to $\tilde{\scr{F}}$ in $\mathrm{Coh}(\tilde{Z})$, where $\tilde{\scr{F}} = j_{\tilde{D}\ast}p_{\tilde{D}}^{\ast}\scr{F}_D$. Now note that the top-row composition $\mathrm{Coh}(\tilde{Z}) \to \mathrm{Coh}(D)$ sends $\tilde{\scr{F}}$ to $\scr{F}_D$. Under the equivalence of Lemma \ref{lemma:small_orlov}, $\scr{F}_D$ is sent to the matrix factorization $F$:
\begin{align*}
    F_0 & = \bigosum_{n \in \ZZ}i_{Z\ast} j_{D\ast}p_{D}^{\ast}\scr{F}_{D}^{2n}(\chi^{-n}), & F_1 & = \bigosum_{m \in \ZZ} i_{Z\ast}j_{D\ast}p_{D}^{\ast}\scr{F}_{D}^{2m+1}(\chi^{-m})
\end{align*}
in $\mathrm{Coh}_{\mathbb{G}_m}(K_X, s(x)y)$. While on the other hand, $\scr{F}$ in $\mathrm{Coh}(Z)$ maps to the sheaf $i_{Z\ast}\scr{F}$, which gives a matrix factorization in $\mathrm{Coh}_{\mathbb{G}_m}(K_X, s(x)y)$. As $\ZZ/2$-graded matrix factorizations, this is the same as $F$ above, since $\scr{F} = j_{D\ast}p_{D}^{\ast}\scr{F}_D$.

Similarly, suppose $h_D: \scr{F}_D \to \scr{G}_D$ is a degree-$k$ morphism of complexes of coherent sheaves on $D$; by Lemma \ref{lemma:small_orlov} the corresponding morphism of coherent sheaves $\tilde{h}: \tilde{\scr{F}} \to \tilde{\scr{G}}$ on $\tilde{Z}$ also arises as $j_{\tilde{D}\ast}p_{\tilde{D}}^{\ast}(h_D)$ and has degree $k$. Since $\tilde{h}$ is independent of $z$, under the localization functor $\tilde{h}$ is simply sent to $h_D$ in $\mathrm{Hom}^{k \;(\mathrm{mod}\; 2)}(\scr{F}_D, \scr{G}_D)$. This means the diagram commutes when the degree of the morphisms is taken mod $2$.
\end{proof}

\subsection{Mirror Symmetry}

Firstly, we may establish our version of Conjecture \ref{conj:main}, part 3:

\begin{theorem}\label{thm:alpha_infinity}
There is a mirror symmetry equivalence between $\scr{W}_{0}(\comp)$ and $\mathrm{Coh}(Z)$ given by composing the mirror equivalence of $\scr{W}_{0}(\comp)$ and $\mathrm{Coh}(\tilde{Z})$ from Theorem \ref{thm:alphas} with the sequence of equivalences $\mathrm{Coh}(\tilde{Z}) \to \mathrm{Coh}(Z)$ in Theorem \ref{prop:orlov}.

Under this equivalence, the functor $\alpha_{\infty}: \scr{W}((\CC^{\ast})^n,f) \to \scr{W}_{0}(\comp)$ is mirror to $\pi_{X}^{\ast}[1]:\mathrm{Coh}(X) \to \mathrm{Coh}(Z)$.
\end{theorem}

\begin{proof}
This follows by applying Lemma \ref{lem:locally_free}.
\end{proof}

The \textit{restriction functor} $\rho: \scr{W}_0(\comp) \to \scr{W}(H)$ from \cite[\S 4]{Speculations} is somewhat unique, so we shall spend some time describing it here. It is far from obvious that such a functor even exists: the choice of grading is particularly important.

The first step in this definition of the restriction functor is given by Viterbo restriction \cite[\S 11.1]{GPS2} from $\comp$ to the subdomain $H \times \CC^{\ast}$ in Proposition \ref{prop:sector_gluing}, part (1): this can be done on the level of $\ZZ$-graded categories. Moreover, since $H \times \CC^{\ast}$ is a Weinstein manifold, the Viterbo restriction functor is defined on categories of twisted complexes \cite[(11.5)]{GPS2}. 

By the K\"unneth theorem of \cite[Corollary 1.18]{GPS2}, the wrapped Fukaya category $\scr{W}(H \times \CC^{\ast}(2))$ is generated by (cylindrizations of) objects of the form $\ell \times \RR_{>0}$ where $\ell$ is a Lagrangian in $H$. On the level of objects, the functor $\rho$ should take $\ell \times \RR_{>0}$ to $\ell$: it is not clear, however, that this is well-defined on morphisms, as there is no choice of reference fiber $H \times \set{p}$ in which we can take the intersections of two Lagrangians in $H \times \CC^{\ast}$. However, this will be well-defined if we use the $0$-grading on $H \times \CC^{\ast}$: this is mirror to the algebraic geometry construction of taking the fiber over the `generic point' in $\AA^1[-1]$. 

To illustrate this, let us begin with the simplest example where $H$ is simply a point:
\begin{proposition}\label{prop:localization}
There is a natural functor from the category of finitely-generated $\CC[z^{\pm 1}]$-modules to the category of $\ZZ/2$-graded finite-dimensional $\CC$ vector spaces when $|z|=2$, taking $\CC[z^{\pm}]$ to $\CC$ in degree $0$.
Moreover, this gives a graded functor from a \textit{graded} dg category of $\CC[z^{\pm 1}]$-modules (modules carrying an additional $\ZZ$-grading for which $z$ has degree $0$), to the \textit{graded} $\cz$-periodic dg category of $\CC$-vector spaces.
\end{proposition}
\begin{proof}
Given $F$ a finitely-generated $\CC[z^{\pm 1}]$-module, it corresponds to a sequence $F^{\bullet}$ of finite-dimensional $\CC$-vector spaces along with a degree-$2$ map $z: F^{k} \to F^{k+2}$. Then we may define $\rho(F)$ by $V_0 = F^{2k}, V_1 = F^{2k+1}$ for any $k \in \ZZ$. Given a morphism $f: F \to G$ of $\CC[z^{\pm 1}]$-modules, that is, a linear map $f^{\bullet}:F^{\bullet} \to G^{\bullet}$ commuting with $z$, we get a morphism $\rho(F) \to \rho(G)$ via $f^0 = f^{2k}: F^{2k} \to G^{2k}$ and $f^{1}= f^{2k+1}: F^{2k+1} \to G^{2k+1}$. Importantly, this definition is independent of our choice of $k \in \ZZ$ and yields a functor from finitely-generated $\CC[z^{\pm 1}]$-modules to $\ZZ/2$-graded $\scr{\CC}$ vector spaces.
\end{proof}

Note that no such construction exists if instead we were to take the degree of $z$ be zero: given a Laurent polynomial in a degree-$0$ variable $z$, simply choosing one coefficient does not define a functor.

In the general situation, morphisms in $\scr{W}(H \times \CC^{\ast}(2))$ between two objects $L_1, L_2$ always have an action by $\CC[z^{\pm 1}]$: if $L_1 = \ell_1 \times \RR_{>0}$ and $L_2 = \ell_2 \times \RR_{>0}$, then the morphisms are given by $\mathrm{Hom}(L_1,L_2) \cong \mathrm{Hom}_{\scr{W}(H)}(\ell_1, \ell_2) \otimes_{\CC} \CC[z^{\pm 1}]$, which has an evident action by $\CC[z^{\pm 1}]$. Given a morphism $p$ between $L_1$ and $L_2$, we have a map
\begin{equation*}
    \mathrm{Hom}^{2k}_{\scr{W}(H \times \CC^{\ast}(2))}(L_1, L_2) \overset{\sim}{\to}  \mathrm{Hom}^{0}_{\scr{W}(H \times \CC^{\ast}(2))}(L_1, L_2) \to \mathrm{Hom}_{\scr{W}(H)}^{\mathrm{even}}(\ell_1, \ell_2)
\end{equation*}
given taking (the degree-$2k$ part of) $p$ to $z^{-k}p$. Similarly, taking (the degree-$2k+1$ part of) $p$ to $z^{-k}p$ gives a map
\begin{equation*}
    \mathrm{Hom}^{2k+1}_{\scr{W}(H \times \CC^{\ast}(2))}(L_1, L_2) \overset{\sim}{\to}  \mathrm{Hom}^{1}_{\scr{W}(H \times \CC^{\ast}(2))}(L_1, L_2) \to \mathrm{Hom}_{\scr{W}(H)}^{\mathrm{odd}}(\ell_1, \ell_2)
\end{equation*}
for any $k \in \ZZ$; one can check that this map commutes with the differentials on each morphism complex and so yields a well-defined functor from the cohomology category $H\scr{W}(H \times \CC^{\ast}(2))$ to the $\ZZ/2$-graded cohomlogy category $H\bar{\scr{W}}(H)$. 

Following \cite[\S 4]{Speculations}, to define this more precisely on the $A_{\infty}$ level, one defines $\rho$ formally on $A_{\infty}$-modules as the pullback via the inclusion functor $i: \scr{W}(H) \to \scr{W}(H \times \CC^{\ast}(2))$ that takes a Lagrangian $\ell$ to $\ell \times \RR_{>0}$. Observe that this construction uses a choice of framing of $H$ inside $\cstar{n}$. Given $L$ in $\scr{W}(H \times \CC^{\ast}(2))$, the pullback module $i^{\ast}(L)$ is never representable by an object of $\scr{W}(H)$ since it has unbounded cohomological support. The key observation of Auroux is that the action of $\CC[z^{\pm 1}]$ induces a $2$-periodic structure on $i^{\ast}(L)$ \cite[\S 4]{Speculations}, so that $i^{\ast}(L)$ \textit{is} representable by an object in the $\cz$-linear category $\bar{\scr{W}}(H)$ where we reduce the grading modulo $2$. In some sense, the category $\scr{W}(H \times \CC^{\ast}(2))$ \textit{is} just a $2$-periodization of $\scr{W}(H)$. The functor $\rho$ simply corresponds to the localization $\scr{W}(H \times \CC^{\ast}(2)) \to \bar{\scr{W}}(H)$ at the degree-$2$ natural transformation $z: \mathrm{id} \to [2]$ given by rotating $\CC^{\ast}$ around $0$ \cite[(4.3)]{Speculations}) to give a $\ZZ/2$-graded category. In summary:

\begin{definition}\label{defn:rho}
The restriction functor $\rho: \bar{\scr{W}}_0(\comp) \to \bar{\scr{W}}(H)$ is defined by Viterbo restriction to the subdomain $H \times \CC^{\ast}$ in Proposition \ref{prop:sector_gluing}, part (1); followed by the localization functor $\bar{\scr{W}}(H \times \CC^{\ast}(2)) \to \bar{\scr{W}}(H)$ at the degree-$2$ natural transformation given by multiplication by the generator $z$ of $\mathrm{End}_{\scr{W}(\CC^{\ast}(2))}(\RR_{>0}) \cong \CC[z^{\pm 1}]$ (see \cite[(4.3)]{Speculations}).
\end{definition}

Recall that we use a bar to denote the reduction of the homological grading of an $A_{\infty}$ category mod-$2$:  $\overline{\scr{W}}_{0}(\comp)$ and $\overline{\mathrm{Coh}}(Z)$. Note that the $\rho$ functor is only a functor of $\cz$-linear $A_{\infty}$ categories, since we are localizing at a natural transformation that has degree $2$.

\begin{proposition}\label{prop:viterbo}
Under the mirror symmetry equivalences above for the function $\tilde{f}_0$, the Viterbo restriction from $\comp$ to $H \times \CC^{\ast}$ is mirror to the pullback from $\tilde{Z}$ to $D \times \CC^{\ast}[-1]$; that is, the following diagram of $A_{\infty}$ categories commutes:
\begin{center}
\begin{tikzcd}
\scr{W}_{0}(\comp) \ar{r} \ar{d}{\mathrm{Vit}} & \mathrm{Coh}(\tilde{Z}) \ar{d}{j_{\tilde{D}}^\ast}\\
 \scr{W}(H \times \CC^{\ast}(2)) \ar{r} & \mathrm{Coh}(D \times \CC^{\ast}[-1])
\end{tikzcd}
\end{center}
\end{proposition}
\begin{proof}
%This follows from \cite[Corollary 5.8]{LVL}.
The top horizontal homological mirror symmetry equivalence in the diagram is obtained by gluing together mirror symmetry equivalences for the sectors $(H\times \CC^{\ast}(2),z)$ and $((\CC^{\ast})^n, f)$ along $H$, corresponding on the B-side to a gluing of the derived schemes $D\times \AA^1[-1]$ and $X$ along $D$. \cite[Corollary 5-8]{LVL} shows that for such a gluing, the Viterbo restriction to $H\times \CC^{\ast}(2)$ is mirror to the pullback to $D\times \CC^\times[-1].$ 
Note that although the result in \cite{LVL} is stated for gluings of (non-derived) toric varieties, it works equally well if we allow some directions in those toric varieties to be derived, corresponding to the grading mirror to the choice of $\tilde{f}_0$; the only input is the description of Viterbo restriction \cite[Proposition 11.2]{GPS2} as the quotient by cocores outside the restricted domain, together with the fact that the above descent statement allows us to identify how these cocores behave under mirror symmetry.
\end{proof}

\begin{theorem}\label{thm:rho}
The restriction functor $\rho:\bar{\scr{W}}_0(\comp) \to \bar{\scr{W}}(H)$ is mirror to the quotient functor $q: \bar{\mathrm{Coh}}(Z) \to \bar{\mathrm{Sing}}(Z)$ composed with the Kn\"orrer periodicity equivalence $\bar{\mathrm{Sing}}(Z) \to \bar{\mathrm{Coh}}(D)$ as functors of $\cz$-linear $A_{\infty}$ categories (Lemma \ref{lemma:small_orlov}), under the mirror symmetry equivalence from Theorem \ref{thm:alpha_infinity}.
\end{theorem}

\begin{proof}
In the diagram below, the top left-hand square commutes by Proposition \ref{prop:viterbo}, and the top right-hand square clearly commutes also since we are localizing at the same natural transformation. It remains to show that the lower square commutes, which is Proposition \ref{prop:sheafstuff}.
\end{proof}

\[\begin{tikzcd}
	{\bar{\mathcal{W}}_0((\mathbb{C}^{\ast})^n\setminus H)} && {\bar{\mathcal{W}}(H \times \CC^{\ast}(2))} && {\bar{\mathcal{W}}(H)} \\
	\\
	{\bar{\mathrm{Coh}}(\tilde{Z})} && {\bar{\mathrm{Coh}}(D \times \mathbb{C}^{\times}[-1])} && {\bar{\mathrm{Coh}}(D)} \\
	{\bar{\mathrm{Coh}}_{\mathbb{G}_m}(K_X \times \mathbb{A}^1, s(x) y z)} \\
	{\bar{\mathrm{Coh}}(Z)} && {} && {\bar{\mathrm{Sing}}(Z)}
	\arrow["{\mathrm{Vit}}"{description}, from=1-1, to=1-3]
	\arrow["{\mathrm{loc}}"{description}, from=1-3, to=1-5]
	\arrow["{j_{\tilde{D}}^{\ast}}"{description}, from=3-1, to=3-3]
	\arrow["{\mathrm{loc}}"{description}, from=3-3, to=3-5]
	\arrow["\sim"{description}, from=1-3, to=3-3]
	\arrow["\sim"{description}, from=1-1, to=3-1]
	\arrow["\sim", from=3-1, to=4-1]
	\arrow["\sim",from=5-1, to=4-1]
	\arrow["\sim"{description}, from=3-5, to=5-5]
	\arrow["\sim"{description}, from=1-5, to=3-5]
	\arrow["\rho"{description}, bend right=-18pt, from=1-1, to=1-5]
	\arrow["q", from=5-1, to=5-5]
	%\arrow["\sim"{description}, from=3-3, to=5-3]
\end{tikzcd}\]

\subsection{The $\ZZ$-graded enhancement}

%We can upgrade both $q$ and $\rho$ to $\ZZ$-graded functors if we are willing to work with graded dg/$A_{\infty}$ categories. 
We have seen so far that the functors $q$ and $\rho$ corresponds under a $\ZZ/2$-graded homological mirror symmetry equivalence --- i.e., an equivalence of $\cz$-linear categories. In this section, we show that the categories involved admit an auxiliary $\ZZ$-grading, and that $q$ and $\rho$ can be made to correspond under an equivalence of graded $\cz$-linear categories.

\begin{definition} We introduce the following terminology:
\begin{itemize}
    \item We shall call the grading that comes from giving $\AA^1$ weight $1$ under the $\GG_m$-action, and the form $\eta_0$ on $\comp$, the \textbf{$0$-grading}; we shall consider this the underlying cohomological grading structure on our categories.
    \item We shall call the grading that comes from giving the fibers of $K_X$ weight $1$ under the $\GG_m$-action, and the form $\eta_{\infty}$ on $\comp$, the \textbf{$\infty$-grading}; we shall consider this the additional $\ZZ$-grading structure on our categories.
\end{itemize}
\end{definition}

Using the $\infty$-grading gives an (additional) grading on the $A_{\infty}$ categories $\scr{W}_{0,\infty}(\comp)$ and $\mathrm{Coh}_{\GG_m}(Z)$: reducing the homological grading mod-$2$ yields graded $\cz$-linear $A_{\infty}$ categories; we shall continue to use a bar to denote the reduction of the homological grading mod-$2$: $\overline{\scr{W}}_{0,\infty}(\comp)$ and $\overline{\mathrm{Coh}}_{\GG_m}(Z)$.

We may now define a graded version of the restriction functor $\rho: \bar{\scr{W}}_{0,\infty}(\comp) \to \bar{\scr{W}}(H)$ by
\begin{itemize}
    \item Viterbo restriction to the subdomain $H \times \CC^{\ast}$ (the definition from \cite[\S 11.1]{GPS2} can be carried out in a graded manner), where $H \times \CC^{\ast}$ carries the restriction of the trivialization of the determinant bundle from $\comp$; then,
    \item Localization $\bar{\scr{W}}_{0,\infty}(H \times \CC^{\ast}) \to \bar{\scr{W}}(H)$ at multiplication by the degree-$2$ natural transformation $z$. Since multiplication by $z$ has degree $0$ for the $\infty$-grading, the latter functor is a graded $A_{\infty}$ functor by Proposition \ref{prop:localization}.
\end{itemize}
Hence $\rho$ is a graded functor between graded $\cz$-linear $A_{\infty}$ categories.

Likewise, the quotient map $q: \mathrm{Coh}_{\GG_m}(Z) \to \mathrm{Sing}_{\GG_m}(Z)$ on equivariant coherent sheaves does not preserve the cohomological grading, but preserves the grading coming from the $\GG_m$-equivariant structure (cf. Lemma \ref{lemma:small_orlov}), which is the $\infty$-grading. Hence it gives a graded functor $q: \bar{\mathrm{Coh}}_{\GG_m}(Z) \to \bar{\mathrm{Sing}}_{\GG_m}(Z)$ of graded $\cz$-linear dg categories (where the $\cz$-linear structure on $\mathrm{Sing}_{\GG_m}(Z)$ comes from the equivalence with the category of matrix factorizations in Theorem \ref{thm:orlov}).

Thus both $\rho$ and $q$ have lifts that are graded functors. Next we show that they correspond under mirror symmetry in a graded sense.

\begin{lemma}\label{lem:bigraded_cup}
We have a mirror symmetry equivalence of graded $A_{\infty}$ categories between $\scr{W}_{0,\infty}((\CC^{\ast})^n, f)$ and $\mathrm{Coh}_{\GG_{m}^2}(X,0)$, along with a commutative diagram
\begin{center}
\begin{tikzcd}
\scr{W}(H) \ar{r} \ar{d}{\cup} & \mathrm{Coh}(D) \ar{d}{i_{D \ast}}\\
 \scr{W}_{0,\infty}( (\CC^{\ast})^n, f) \ar{r} & \mathrm{Coh}_{\GG_m^2}(X,0)
\end{tikzcd}
\end{center}
where the top horizontal arrow is the equivalence of \cite{BenVivek}.
\end{lemma}
\begin{proof}
This is simply the observation that both trivializations of the determinant bundle on any generic fiber $H$ and on $((\CC^{\ast})^n,f)$ agree up to homotopy; combined with the fact that both $\GG_m$-actions are trivial on $X$. The former can be seen from the fact that each gluing presentation in Propositions \ref{prop:sector_gluing} and \ref{prop:A2} leaves the sector $((\CC^{\ast})^n,f)$ unchanged. 
\end{proof}

\begin{theorem}\label{thm:bigraded_ms}
There is a mirror symmetry equivalence of graded $A_{\infty}$ categories between $\scr{W}_{0,\infty}(\comp)$ and $\mathrm{Coh}_{\GG_m^2}(\tilde{S},0)$.
\end{theorem}
\begin{proof}
By inspection of the arguments of \cite[\S 11.2]{GPS2}, one can see that the sectorial decomposition from Proposition \ref{prop:sector_gluing} yields a pushout diagram of graded $A_{\infty}$ categories (cf. \cite[\S 2.2]{GPS2}):
\begin{center}
    \begin{tikzcd}
\scr{W}(H) \ar{r}{\cup} \ar{d} & \scr{W}_{0,\infty}((\CC^{\ast})^n, f) \ar{d}\\
 \scr{W}_{0,\infty}(H \times \CC^{\ast}, z) \ar{r} & \scr{W}_{0,\infty}( (\CC^{\ast})^n \setminus H)
\end{tikzcd}
\end{center}
We compare this gluing diagram with the gluing diagram from Proposition \ref{prop:bigraded_descent}:
\[\begin{tikzcd}
	& {\mathrm{Coh}_{\GG_{m}^2}(D,0)} && {\mathrm{Coh}_{\GG_{m}^2}(D \times \mathbb{A}^1,0)} \\
	{\mathcal{W}(H)} && {\mathcal{W}_{0,\infty}(H \times \CC^\ast,z)} \\
	& {\mathrm{Coh}_{\GG_{m}^2}(X,0)} && {\mathrm{Coh}_{\GG_{m}^2}(\tilde{S},0)} \\
	{\mathcal{W}_{0,\infty}((\mathbb{C}^{\ast})^n, f)} && {\mathcal{W}_{0,\infty}((\mathbb{C}^{\ast})^n \setminus H)}
	\arrow["\cup"{description, pos=0.3}, from=2-1, to=2-3]
	\arrow[""{name=0, anchor=center, inner sep=0}, "\cup"', from=2-1, to=4-1]
	\arrow[""{name=1, anchor=center, inner sep=0}, "{\alpha_{\infty}}"', from=4-1, to=4-3]
	\arrow["\sim", from=2-1, to=1-2]
	\arrow["\sim", from=2-3, to=1-4]
	\arrow["\sim", from=4-3, to=3-4]
	\arrow["\sim", from=4-1, to=3-2]
	\arrow["{\iota_{\tilde{D}\ast}}"{description}, from=1-2, to=1-4]
	\arrow["{i_{D\ast}}"{description, pos=0.3}, from=1-2, to=3-2]
	\arrow["{j_{\tilde{D}\ast}}"{description}, from=1-4, to=3-4]
	\arrow["{i_{\tilde{X}\ast}}"{description, pos=0.6}, from=3-2, to=3-4]
	\arrow[from=2-3, to=4-3]
\end{tikzcd}\]
using Lemma \ref{lem:bigraded_cup} to give the desired equivalence.
\end{proof}

\begin{proposition}\label{prop:bigraded_viterbo}
We have a commutative diagram of graded $A_{\infty}$ categories:
\begin{center}
\begin{tikzcd}
\scr{W}_{0,\infty}(\comp) \ar{r} \ar{d}{\mathrm{Vit}} & \mathrm{Coh}_{\GG_m}(\tilde{Z}) \ar{d}{j_{\tilde{D}}^\ast}\\
 \scr{W}_{0,\infty}(H \times \CC^{\ast}) \ar{r} & \mathrm{Coh}_{\GG_m}(D \times \CC^{\ast}[-1])
\end{tikzcd}
\end{center}
\end{proposition}

\begin{proof}
As in \Cref{prop:viterbo}, the arguments of \cite[Corollary 5-8]{LVL}, applied to the above descent diagram (now of graded $A_\infty$ categories) for $\scr{W}_{0,\infty}((\CC^{\ast})^n\setminus H),$ show that Viterbo restriction is mirror to the indicated pullback.
%By inspection of the arguments of \cite[Corollary 5.8]{LVL} and \cite[\S 3]{FLTZ}, one can see that the results extend to graded $A_{\infty}$ categories.
\end{proof}

Homological mirror symmetry equivalences at the level of graded $A_{\infty}$ categories in particular descends to the level of graded $\cz$-linear $A_{\infty}$ categories by forgetting the grading mod-$2$.

\begin{proposition}\label{prop:bigraded_sheafstuff}
Under the Kn\"orrer equivalences in Theorem \ref{thm:bigraded_orlov}, we have a commutative diagram of graded $\cz$-linear dg categories:
\begin{center}
\begin{tikzcd}
	{\bar{\mathrm{Coh}}_{\GG_{m}^2}(\tilde{S},0)} \ar{rr}{j_{\tilde{D}}^{\ast}} \ar{d}{\sim} && {\bar{\mathrm{Coh}}_{\GG_m^2}(D \times \mathbb{C}^{\times},0)} \ar{rr}{\mathrm{loc}} && {\bar{\mathrm{Coh}}(D)} \ar{dd}{\sim} \\
	{\bar{\mathrm{Coh}}_{\mathbb{G}_{m}^2}(K_X \times \mathbb{A}^1, s(x)yz)} \\
	{\bar{\mathrm{Coh}}_{\GG_{m}^2}(Z,0)} \ar{u}{\sim} \ar{rrrr}{q} && {} && {\bar{\mathrm{Sing}}_{\GG_m}(Z)} \\
\end{tikzcd}
\end{center}
where $\mathrm{Coh}_{\GG_{m}^2}(Z,0)$ is identified with $\mathrm{Coh}_{\GG_m}(Z)$ as in Lemma \ref{lem:trivial_action2}.
\end{proposition}
\begin{proof}
We inspect the proof of Proposition \ref{prop:sheafstuff}. Firstly, it is not difficult to see that Lemmas \ref{lemma:orlov} and \ref{lem:locally_free} extend also to graded dg versions. Thus as before, going either way around the diagram sends all locally free equivariant coherent sheaves on $Z$ to zero. Importantly, the equivalence in Lemma \ref{lemma:small_orlov} that factors through $\mathrm{Coh}_{\GG_m}(Z)$ is a graded dg functor. Now, the vertical left equivalences are graded dg equivalences by Theorem \ref{thm:bigraded_orlov}; while the top arrows are only graded $\cz$-linear dg functors, as the localization $2$-periodizes the homological grading. Combining these observations with the proof of Proposition \ref{prop:sheafstuff} gives the result.
\end{proof}

\begin{theorem}\label{thm:bigraded_rho}
The restriction functor $\rho:\bar{\scr{W}}_{0,\infty}(\comp) \to \bar{\scr{W}}(H)$ is mirror to the quotient $q: \bar{\mathrm{Coh}}_{\GG_m}(Z) \to \bar{\mathrm{Sing}}_{\GG_m}(Z)$ composed with the Kn\"orrer periodicity equivalence $\mathrm{Sing}_{\GG_m}(Z) \to \mathrm{Coh}(D)$ from Lemma \ref{lemma:small_orlov} as graded $\cz$-linear $A_{\infty}$ functors, under the mirror symmetry equivalences from Theorem \ref{thm:bigraded_ms} composed with the Orlov equivalences from Theorem \ref{thm:bigraded_orlov}.
\end{theorem}

\begin{proof}
As before, this follows by combining Propositions \ref{prop:bigraded_viterbo} and \ref{prop:bigraded_sheafstuff}.
\end{proof}

\[\begin{tikzcd}
	{\bar{\mathcal{W}}_{0,\infty}((\mathbb{C}^{\ast})^n\setminus H)} && {\bar{\mathcal{W}}_{0,\infty}(H \times \CC^{\ast})} && {\bar{\mathcal{W}}(H)} \\
	\\
	{\bar{\mathrm{Coh}}_{\GG_{m}^2}(\tilde{S},0)} && {\bar{\mathrm{Coh}}_{\GG_{m}^2}(D \times \mathbb{C}^{\times},0)} && {\bar{\mathrm{Coh}}(D)} \\
	{\bar{\mathrm{Coh}}_{\mathbb{G}_{m}^2}(K_X \times \mathbb{A}^1, s(x) y z)} \\
	{\bar{\mathrm{Coh}}_{\GG_{m}^2}(Z,0)} && {} && {\bar{\mathrm{Sing}}_{\GG_m}(Z)}
	\arrow["{\mathrm{Vit}}"{description}, from=1-1, to=1-3]
	\arrow["{\mathrm{loc}}"{description}, from=1-3, to=1-5]
	\arrow["{j_{\tilde{D}}^{\ast}}"{description}, from=3-1, to=3-3]
	\arrow["{\mathrm{loc}}"{description}, from=3-3, to=3-5]
	\arrow["\sim"{description}, from=1-3, to=3-3]
	\arrow["\sim"{description}, from=1-1, to=3-1]
	\arrow["\sim", from=3-1, to=4-1]
	\arrow[from=5-1, to=4-1]
	\arrow["\sim"{description}, from=3-5, to=5-5]
	\arrow["\sim"{description}, from=1-5, to=3-5]
	\arrow["\rho"{description}, bend right=-18pt, from=1-1, to=1-5]
	\arrow["q", from=5-1, to=5-5]
	%\arrow["\sim"{description}, from=3-3, to=5-3]
\end{tikzcd}\]

\section{The Lifting Functor}\label{sec:j}

The hypersurface $H$ carries a natural framing given by the function $f$. We will first need to verify that this agrees with the framing given by the product neighbourhood $H \times \set{\mathrm{Re}(z) \geq -\rho} \subset \cstar{n}$ in \cite[Proposition 2.6]{sylvan_lemma}. Note that over the unperturbed part of the fiber (called $f\inv(1/2)$ in Sylvan's notation), the vector fields uniquely specified by Sylvan's prescription are $X = (g \partial_r, 0), Y = (0,R_\lambda)$ where $R_\lambda$ is the Reeb vector field of $\partial X$ and the smooth function $g: \partial X \to \RR$ is chosen so that $\omega_{X}(g \partial_r, R_{\lambda}) = 1$. These are positive multiples of $\grad |f|$ and $\grad \mathrm{arg}(f)$ respectively, which agrees with the framing obtained from $f$ over the interior region of the fiber. Since the full fiber $\hat{F}$ is obtained by attaching an infinite cylinder to $\partial F$ (by \cite[Proposition 2.6, condition (3)]{sylvan_lemma}), the result of extending this framing over the whole fiber is homotopic to Sylvan's construction. 

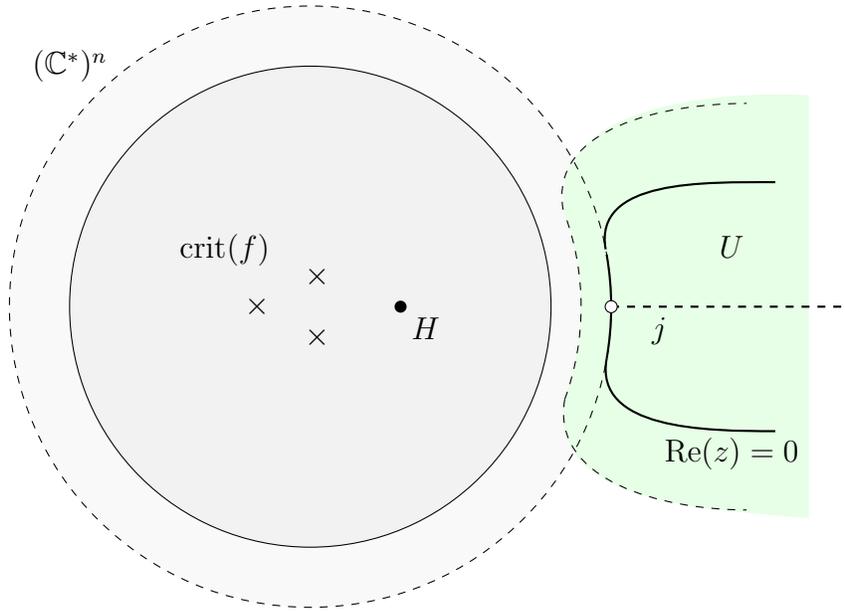
\begin{figure}%Step 2
\begin{center}
    \begin{tikzpicture}[scale=0.8]
%\draw[help lines] (-8,-8) grid (8,8);
\fill [fill=gray!5] (0,0) circle (5);
\fill [fill=gray!10] (0,0) circle (4);
\draw (0,0) circle (4);

\fill[fill=green!10] (-23:9) to[out=175, in=-100] (-20:4.5) arc[start angle=-20, end angle=20,radius=4.5] (20:4.5)  to[out=100,in=177] (23:9) to (-23:9);

% \fill[fill=red!10] (5:6) to (-5:6) to[out=-100, in=180] (-10:8) to (10:8) to[out=180,in=90] (5:6);

\draw[dashed] (0,0) circle (5);
\draw[thick] (-10:5) arc[start angle=-10, end angle=10,radius=5];
\draw[thick] (11:5) to[out=100,in=180] (15:8);
\draw[thick] (-10:5) to[out=-100, in=180] (-15:8);

\draw[dashed] (-20:4.5) arc[start angle=-20, end angle=20,radius=4.5];
\draw[dashed] (21:4.5) to[out=100,in=180] (25:8);
\draw[dashed] (-20:4.5) to[out=-100, in=180] (-25:8);

\draw[thick,dashed] (0:5) to (0:9);
\draw[fill=white] (0:5) circle (0.1);

\fill[black] (1.5,0) circle (0.1);
\node[anchor=north west] at (1.5,0) {$H$};
\draw (-0.5,0) node[left] {$\huge{\times}$};
    \draw (0.5,-0.5) node[left] {$\huge{\times}$};
    \draw (0.5,0.5) node[left] {$\huge{\times}$};
    \node[anchor=south east] at (-0.5,0.5) {$\mathrm{crit}(f)$};
    \node[anchor=north west] at (5.5,0) {$j$};
\node at (7,1) {$U$};
\node[anchor=north] at (7,-2) {$\mathrm{Re}(z)=0$};
\node at (-4,4) {$\cstar{n}$};
    \end{tikzpicture}
    \caption{The lifting functor $j$ corresponds to transporting a Lagrangian along the dashed horizontal line running from the puncture to $\infty$ along the positive real axis.}
    \label{fig:negative_lifting}
\end{center}
\end{figure}

\begin{definition}\label{defn:lifting} We define the lifting functor $j: \scr{W}(H) \to \scr{W}(H \times \CC^{\ast},z) \to \scr{W}(\comp)$ using a sector decomposition as in Proposition \ref{prop:sector_gluing}, by taking the cylindrization \cite[\S 7.2]{GPS2} of the product with $\RR_{\geq 0}$ followed by sector inclusion (see Figure \ref{fig:negative_lifting}).
\end{definition}

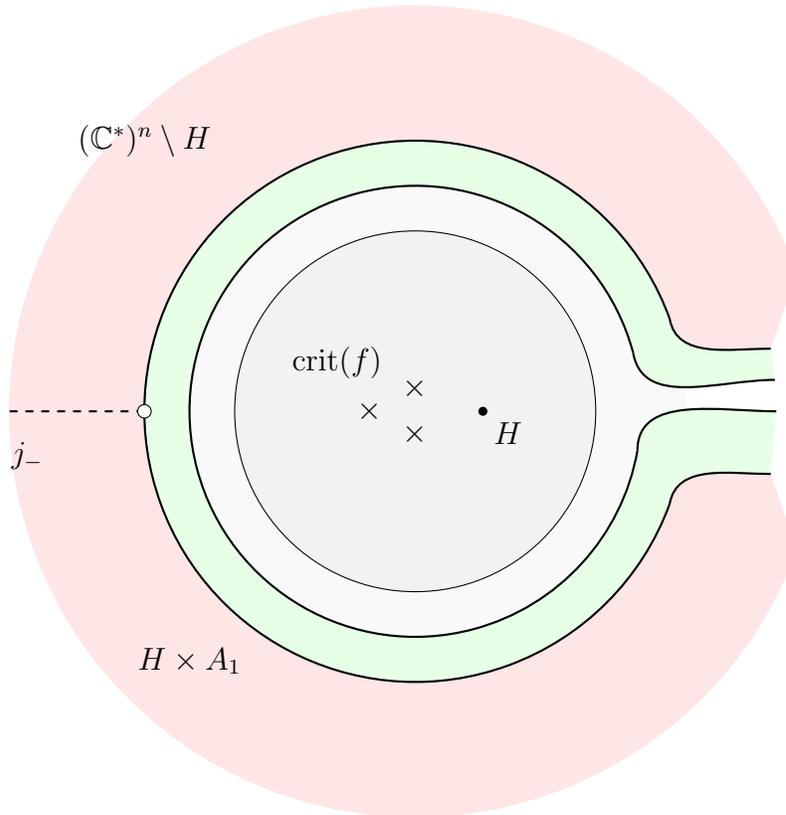
\begin{figure}%Step 5
\begin{center}
    \begin{tikzpicture}[scale=0.6]
%\draw[help lines] (-8,-8) grid (8,8);
\fill [fill=gray!5] (0,0) circle (6);
\fill [fill=gray!10] (0,0) circle (4);
\draw (0,0) circle (4);

\fill[fill=green!10] (15:5) arc[start angle=15, end angle=350,radius=5] (-10:5) to[out=90, in=180] (0:8) to (-10:8) to[out=180, in=80] (-20:6) arc[start angle=-20, end angle=- 340,radius=6] to[out=-80,in=180] (10:8) to (5:8) to[out=180,in=-80] (15:5);

\fill[fill=red!10] (20:6) arc[start angle=20, end angle=340,radius=6] (-20:6) to[out=90, in=180] (-10:8) to (-20:9) arc[start angle=-20, end angle=- 340,radius=9] to[out=-80,in=180] (20:9) to (10:8) to[out=180,in=-80] (20:6);

\draw[thick] (15:5) arc[start angle=15, end angle=350,radius=5];
\draw[thick] (15:5) to[out=-80,in=180] (5:8);
\draw[thick] (-10:5) to[out=90, in=180] (0:8);

\draw[thick] (20:6) arc[start angle=20, end angle=340,radius=6];
\draw[thick] (20:6) to[out=-80,in=180] (10:8);
\draw[thick] (-20:6) to[out=80, in=180] (-10:8);
\draw[thick,dashed] (-180:6) to (-180:9);
\draw[fill=white] (-180:6) circle (0.15);
\fill[black] (1.5,0) circle (0.1);
\node[anchor=north west] at (1.5,0) {$H$};
\draw (-0.5,0) node[left] {$\huge{\times}$};
    \draw (0.5,-0.5) node[left] {$\huge{\times}$};
    \draw (0.5,0.5) node[left] {$\huge{\times}$};
    \node[anchor=south east] at (-0.5,0.5) {$\mathrm{crit}(f)$};
\node[anchor=east] at (-8,-1) {$j_{-}$};
\node at (-6,6) {$\cstar{n}\setminus H$};
\node[anchor=north] at (-5,-5) {$H \times A_1$};
    \end{tikzpicture}
        \caption{The negative lifting functor $j_{-}$ corresponds to transporting a Lagrangian along the dashed horizontal line running from the puncture to $-\infty$.}
    \label{fig:positive_lifting}
\end{center}
\end{figure}

\begin{remark}
One could also define a negative lifting functor $j_{-}: \scr{W}(H) \to \scr{W}(H \times A_1) \to \scr{W}(\comp)$ by cylindrization of the product with $\RR_{\leq 0}$, followed by sector inclusion (see Figure \ref{fig:positive_lifting}). This negative lifting functor $j_{-}$ would differ from the lifting functor $j$ by a twist by the counterclockwise monodromy $\mu$ of $f$ around $\infty$.
\end{remark}

\begin{theorem}\label{thm:negative_lift}
The lifting functor $j:\scr{W}(H) \to \scr{W}_{0}(\comp)$ is mirror to $j_{D \ast}p_{D}^{\ast}: \mathrm{Coh}(D) \to \mathrm{Coh}(Z)$ for the equivalence $\scr{W}_{0}(\comp) \to \mathrm{Coh}(Z)$ as in Theorem \ref{thm:alpha_infinity}, where $j_{D}: K_{X}|_D \to Z, p_{D}: K_{X}|_D \to D$ are the inclusion and projection respectively.
\end{theorem}
\begin{proof}
We first consider the diagram below: 
\[\begin{tikzcd}
	{\mathcal{W}(H)} & {\mathcal{W}(H \times\CC^{\ast}(2),z)} & {\mathcal{W}_{0}((\mathbb{C}^{\ast})^n\setminus H)} \\
	{\mathrm{Coh}(D)} & {\mathrm{Coh}(D \times \mathbb{A}^1[-1])} & {\mathrm{Coh}(\tilde{Z})}
	\arrow["{\times \mathbb{R}_{\leq 0}}", from=1-1, to=1-2]
	\arrow[from=1-2, to=1-3]
	\arrow["{p_{\tilde{D}}^{\ast}}", from=2-1, to=2-2]
	\arrow["{j_{\tilde{D} \ast}}", from=2-2, to=2-3]
	\arrow["\sim", from=1-3, to=2-3]
	\arrow["\sim", from=1-2, to=2-2]
	\arrow["\sim"', from=1-1, to=2-1]
\end{tikzcd}\]
the right-hand square commutes by the diagram from Theorem \ref{thm:alphas}, and the left-hand square commutes since both functors correspond to tensoring with $\CC[z]$ with $|z|=2$. Thus $j:\scr{W}(H) \to \scr{W}_{0}(\comp)$ is mirror to $j_{\tilde{D}\ast}p_{\tilde{D}}^{\ast}: \mathrm{Coh}(D) \to \mathrm{Coh}(\tilde{Z})$. Next, by Lemma \ref{lemma:orlov}, under the Kn\"orrer equivalence $\mathrm{Coh}(\tilde{Z}) \to \mathrm{Coh}(Z)$, $j_{\tilde{D}\ast}p_{\tilde{D}}^{\ast}$ is sent to $j_{D\ast}p_{D}^{\ast}$ again.
\end{proof}

\begin{remark}
Because of the twist by the monodromy, the negative lifting functor $j_{-}$ would be mirror to $j_{D \ast}p_{D}^{\ast}(\scr{K}_{X}\inv|_{D} \otimes \cdot\;)$ for the mirror equivalence of $\scr{W}_{0}(\comp)$ and $\mathrm{Coh}(Z)$, since by Proposition \ref{prop:schober_stuff}, the counterclockwise monodromy $\mu$ is mirror to tensoring by $\scr{K}_{X}\inv|_{D}$. This would correspond to Conjecture \ref{conj:main}, (4).
\end{remark}

\begin{lemma}\label{prop:triangle}
There is an exact triangle of functors $\mathrm{Coh}(X) \to \mathrm{Coh}(Z)$:
    \begin{center}
\begin{tikzcd}
j_{D \ast}p_{D}^{\ast}i_{D}^{\ast} \ar{rr}{+1} & {} & i_{X \ast} \ar{dl}\\
{} & \pi_{X}^{\ast} \ar{ul} & {}
\end{tikzcd}
\end{center}
\end{lemma}
\begin{proof}
The map $i_{X \ast}\scr{F} \to \pi_{X}^{\ast} \scr{F}$ is the inclusion of a subsheaf: the cokernel is exactly the restriction $\scr{F}|_{D}$ pulled back under the projection $K_{X}|_D \to D$.
\end{proof}

\begin{corollary} \label{cor:triangle}
There is an exact triangle relating $\alpha_0, \alpha_{\infty}$ and $j \cap$ as functors $\scr{W}((\CC^{\ast})^n, f) \to \scr{W}_{\infty}(\comp)$:
\begin{center}
\begin{tikzcd}
j \cap \ar{rr}{+1} & {} & \alpha_{0} \ar{dl}\\
{} & \alpha_\infty \ar{ul} & {}
\end{tikzcd}
\end{center}
and there is a corresponding exact triangle for $\scr{W}_{0}(\comp)$ and $j_{-}$.
\end{corollary}
\begin{proof}
This is just the mirror to Lemma \ref{prop:triangle}, using Theorem \ref{thm:alpha_infinity} to see that $\pi_{X}^{\ast}$ is mirror to $\alpha_{\infty}$.
\end{proof}

\subsection{An Example}\label{subsec:example}

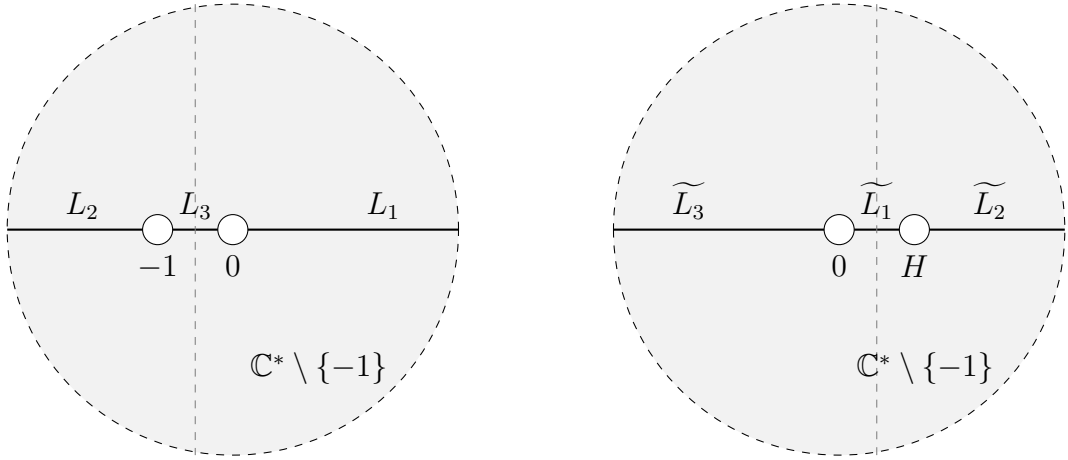
\begin{figure}
    \centering
\begin{tikzpicture}
\filldraw[dashed,fill=gray!10] (0,0) circle (3);
\draw[thick] (-3,0) to (3,0);
\draw[dashed,gray] (-0.5,3) to (-0.5,-3);
\draw[fill=white] (0,0) circle (0.2);
\draw[fill=white] (-1,0) circle (0.2);
\node[anchor=south east] at (2.2,-2.2) {$\CC^{\ast}\setminus\set{-1}$};
\node[anchor=north] at (0,-0.2) {$0$};
\node[anchor=north] at (-1,-0.2) {$-1$};
\node[anchor=south] at (-2,0) {$L_2$};
\node[anchor=south] at (-0.5,0) {$L_3$};
\node[anchor=south] at (2,0) {$L_1$};
\end{tikzpicture}
\hspace{50pt}
\begin{tikzpicture}
\begin{scope}[rotate=180]
\filldraw[dashed,fill=gray!10] (0,0) circle (3);
\draw[thick] (-3,0) to (3,0);
\draw[dashed,gray] (-0.5,3) to (-0.5,-3);
\draw[fill=white] (0,0) circle (0.2);
\draw[fill=white] (-1,0) circle (0.2);
\node[anchor=south east] at (-2.2,2.2) {$\CC^{\ast}\setminus\set{-1}$};
\node[anchor=north] at (0,0.2) {$0$};
\node[anchor=north] at (-1,0.2) {$-1$};
\node[anchor=south] at (-2,0) {$\tilde{L_2}$};
\node[anchor=south] at (-0.5,0) {$\tilde{L_1}$};
\node[anchor=south] at (2,0) {$\tilde{L_3}$};
\end{scope}
\end{tikzpicture}
    \caption{Two gluing presentations of the pair of pants from $\tilde{f_0}$ and $\tilde{f}_{\infty}$ (left to right) with generating Lagrangians labelled. The dashed vertical line represents the sectorial hypersurface dividing the two subsectors.}
    \label{fig:example2}
\end{figure}

In this section, we take up the example of $H = \set{-1}$ once again to illuminate the need for shifts by $[1]$ and different $\ZZ$-gradings. Recall that we make the following identifications, firstly, on the $B$-side:
\begin{itemize}
    \item $X = \mathrm{Spec}(\CC[x])$;
    \item $K_X = \mathrm{Spec}(\CC[x,y])$;
    \item $\AA^1[-1] = \mathrm{Spec}(\CC[t])$, $|t|=-1$.
\end{itemize}
On the $A$-side, there are two different gluing presentations of $\Pi_1 = \CC^{\ast}\setminus H$ given by Propositions \ref{prop:sector_gluing} and \ref{prop:A2}, as illustrated in Figure \ref{fig:example2}. In this case, the two choices of gradings can be easily described:
\begin{align*}
    \eta_0 & = \frac{\dd{z}}{z}, & \eta_{\infty} & = \frac{\dd{z}}{z(z+1)}
\end{align*}
For $\eta_0$, the simple clockwise Reeb orbit around $-1$ will have degree $2$, while those around the other two punctures will have degree $0$. The situation is reversed for $\eta_{\infty}$, where the simple clockwise Reeb orbit around $\infty$ will have degree $2$. This suggests that, in an informal sense, the equivalence $\scr{W}_{0}(\Pi_1) \to \scr{W}_{\infty}(\Pi_1)$ `interchanges zero and infinity'.

As illustrated in Figure \ref{fig:example2}, we use $L_1, L_2, L_3$ to denote the generating set of Lagrangians for the wrapped Fukaya category with the $0$-grading, and $\tilde{L_1}, \tilde{L_2}, \tilde{L_3}$ for those Lagrangians with the $\infty$-grading. In fact, one can show that composing with the Kn\"orrer periodicity equivalence $\mathrm{Coh}(Z) \to \mathrm{Coh}(\tilde{Z})$ interchanges $L_i$ with $\tilde{L}_i$. In particular, $j_{-}$ will correspond to $j$: in general, this will involve twisting by the monodromy. 

\begin{proposition}
Under the mirror symmetry equivalences induced by Theorems \ref{thm:alphas} and \ref{thm:alpha_infinity}, we have
\begin{itemize}
    \item $\alpha_{\infty}(L_1) \mapsto \pi_{\AA^{1}_x}^{\ast}(\scr{O}_{\AA^{1}_x})[1]$,
    \item $\alpha_{0}(\tilde{L}_2) \mapsto i_{\AA^{1}_x\ast} (\scr{O}_{\AA^{1}_x})$,
    \item $\rho(L_2) \mapsto q(\scr{O}_{\AA^{1}_y})$;
    \item $j(\CC) \mapsto \pi^{\ast}_{\AA^{1}_y}(\scr{O}_0)$.
\end{itemize}
\end{proposition}
\begin{proof}
By our definition, $\alpha_\infty$ takes the one generator $L_1$ of $\scr{W}(\CC^{\ast},z)$ to $L_1$ inside $\scr{W}_{0}(\Pi_1)$, and likewise $\alpha_0$ takes the generator of $\scr{W}(\CC^{\ast},z)$ to $\tilde{L_2}$ in $\scr{W}_{\infty}(\Pi_1)$. Under the mirror equivalence induced by gluing, the former is sent to the pushforward of the module $\CC[x]$ to $\AA^1 \cup_0 \AA^1[-1] = \tilde{Z}$. Likewise, $\tilde{L_2}$ corresponds to $\CC[x]$ as a $\CC[x,y]/(xy)$-module, which verifies the correspondence between $\alpha_0$ and $i_{X\ast}$.

Next we need to understand the Kn\"orrer periodicity equivalence
\begin{center}
\begin{tikzcd}
\mathrm{Coh}\;\CC[x,y]/(xy) \ar{r} & \mathrm{Coh}_{\mathbb{G}_m}(\AA^1 \times \AA^1 \times \mathbb{A}^1, x y t) & \ar{l} \mathrm{Coh}\;\CC[x,t]/(xt)
\end{tikzcd}
\end{center}

Under this equivalence, $\CC[x]$ in $\mathrm{Coh}(\AA^1 \cup_0 \AA^1[-1])$ is sent to the matrix factorization
\begin{center}
\begin{tikzcd}
\CC[x,y,t] \ar[bend left]{r}{t} & \CC[x,y,t], \ar[bend left]{l}{xy}
\end{tikzcd}
\end{center}
the $[1]$-shift of
\begin{center}
\begin{tikzcd}
\CC[x,y,t] \ar[bend left]{r}{xy} & \CC[x,y,t], \ar[bend left]{l}{t}
\end{tikzcd}
\end{center}
which is in turn the image of $\CC[x,y]/(xy)$ from $\mathrm{Coh}(\AA^1 \cup_0 \AA^1)$ (note the significance of $t$ having weight $1$). This module is pulled back from $\CC[x]$, the generator of $\mathrm{Coh}(X),$ under the projection $\AA^1 \cup_0 \AA^1 \to \AA^1$. Hence we have shown that under this sequence of identifications, $\alpha_{\infty}(L_1)$ corresponds to $\pi_{X}^{\ast}(\scr{O}_X)[1]$.

Verifying the mirror equivalence for the functor $\rho$ is also simple. The restriction functor $\rho$ sends $L_2$ to $\CC$, while $L_2$ corresponds to $\CC[t]$ in $\mathrm{Coh}(\AA^1 \cup_0 \AA^1[-1])$ under our gluing construction of a mirror symmetry equivalence. Note that the endomorphism algebra of $\CC[t]$ is $\CC[x,z]/(xz)$, $|z|=2$, just as for $L_2$. Passing this module through Kn\"orrer periodicity gives $\CC[y]$ as a module over $\CC[x,y]/(xy)$, which does indeed map to $\CC$ under the quotient map to $\mathrm{Sing}(\CC[x,y]/(xy))$. Similarly, the lifting functor $j$ has image $L_2$ which as we have just seen, corresponds to $\CC[y]$, the pullback of $\scr{O}_0$ under the projection $\AA^1 \cup_0 \AA^1 \to \AA^1$.
\end{proof}

\nocite{*}
\bibliography{biblio}{}
\bibliographystyle{amsalpha}

\end{document}